\documentclass{my_m2an}
\usepackage{amsmath,amsfonts,amssymb,latexsym,epic,eepic}
\usepackage{stmaryrd} 
\usepackage{dsfont} 
\usepackage{enumitem} 
\usepackage{graphics}
\usepackage{pstricks}
\usepackage{fullpage}
\usepackage{cases}
\usepackage{amsthm}
\usepackage{empheq}
\usepackage{float}
\usepackage{accents}
\usepackage{epsfig}
\usepackage{tikz}
\definecolor{bleuf}{rgb}{0.,0.,0.8}
\definecolor{bleuc}{rgb}{0.87,0.92,1.}
\definecolor{vertf}{rgb}{0,0.55,0.1}
\definecolor{vertc}{rgb}{0.6,0.95,0.75}
\definecolor{or}{rgb}{0.98,0.6,.1}
\definecolor{rougec}{rgb}{1,0.4,0.}
\newcommand{\gradi}{{\boldsymbol \nabla}}
\newcommand{\dive}{{\rm div}}
\newcommand{\bfa}{{\boldsymbol a}}
\newcommand{\bfe}{{\boldsymbol e}}
\newcommand{\bfn}{{\boldsymbol n}}
\newcommand{\bfu}{{\boldsymbol u}}
\newcommand{\bfx}{{\boldsymbol x}}
\newcommand{\bfW}{{\boldsymbol{\rm W}}}
\newcommand{\ei}{^{(i)}}
\newcommand{\eg}{\emph{e.g.}}

\newcommand{\Ds}{{\scalebox{0.6}{$D_\edge$}}}
\newcommand{\mesh}{{\mathcal M}}
\newcommand{\edge}{\sigma}
\newcommand{\edges}{{\mathcal E}}
\newcommand{\edgesint}{{\mathcal E}_{{\rm int}}}
\newcommand{\edgesext}{{\mathcal E}_{{\rm ext}}}
\newcommand{\scheme}{\mathcal{S}}
\newcommand{\edgesischeme}{\edges^{(i)}_\scheme}
\newcommand{\edgesdischeme}{\edgesd^{(i)}_\scheme}
\newcommand{\edged}{\varepsilon}
\newcommand{\edgesd}{\tilde {\mathcal E}}

\newcommand{\edgesdext}{\tilde{\mathcal E}_{{\rm ext}}}

\newcommand{\Dhi}{\Delta h_{f,i}^0}
\begin{document}

\title{A staggered pressure correction numerical scheme to compute a travelling reactive interface in a partially premixed mixture}

\author{D. Grapsas}
\address{Aix-Marseille Universit\'e, CNRS,   Centrale  Marseille,  I2M, UMR 7373, 13453 Marseille, France, \\ {dionysios.grapsas@univ-amu.fr}}
\author{R. Herbin}
\address{Aix-Marseille Universit\'e, CNRS,   Centrale  Marseille,  I2M, UMR 7373, 13453 Marseille, France, \\ {raphaele.herbin@univ-amu.fr}}
\author{J.-C. Latch\'e}
\address{Institut de Radioprotection et de S\^{u}ret\'{e} Nucl\'{e}aire (IRSN), BP 3, 13115 Saint-Paul-lez-Durance cedex, France \\ (jean-claude.latche@irsn.fr)}
\author{Y. Nasseri}
\address{Aix-Marseille Universit\'e, CNRS,   Centrale  Marseille,  I2M, UMR 7373, 13453 Marseille, France, \\ {youssouf.nasseri@univ-amu.fr}}

\begin{abstract}
We address in this paper a model for the simulation of turbulent deflagrations in industrial applications.
The flow is governed by the Euler equations for a variable composition mixture and the combustion modelling is based on a phenomenological approach: the flame propagation is represented by the transport of the characteristic function of the burnt zone, where the chemical reaction is complete; outside this zone, the atmosphere remains in its fresh state.
Numerically, we approximate this problem by a penalization-like approach, \ie\ using a finite conversion rate with a characteristic time tending to zero with the space and time steps.
The numerical scheme works on staggered, possibly unstructured, meshes.
The time-marching algorithm is of segregated type, and consists in solving in a first step the chemical species mass balances and then, in a second step, mass, momentum and energy balances.
For this latter stage of the algorithm, we use a pressure correction technique, and solve a balance equation for the so-called sensible enthalpy instead of the total energy balance, with corrective terms for consistency. 
The scheme is shown to satisfy the same stability properties as the continuous problem: the chemical species mass fractions are kept in the $[0,1]$ interval, the density and the sensible internal energy stay positive and the integral over the computational domain of a discrete total energy is conserved.
In addition, we show that the scheme is in fact conservative, \ie\  that its solution satisfy a conservative discrete total energy balance equation, with space and time discretizations which are unusual but consistent in the Lax-Wendroff sense.
Finally, we observe numerically that the penalization procedure converges, \ie\ that making the chemical time scale tend to zero allows to converge to the solution of the target (infinitely fast chemistry) continuous problem.
Tests also evidence that the scheme accuracy dramatically depends on the discretization of the convection operator in the chemical species mass balances. 
\end{abstract}

\date{\today}

\maketitle

\tableofcontents
%
%

\section{Problem position}\label{sec:int}

In this paper, we study a numerical scheme for the computation of large scale turbulent deflagrations occurring in a partially premixed atmosphere.
In usual situations, such a physical phenomena is driven by the progress in the atmosphere of a shell-shaped thin zone, where the chemical reaction occurs and which thus separates the burnt area from fresh gases; this zone is called the flame brush.
The onset of the chemical reaction is due to the temperature elevation, so the displacement of the flame brush is driven by the heat transfers inside and in the neighbour of this zone.
Modelling of deflagrations still remains a challenge, since the flame brush has a very complex structure (sometimes presented as fractal in the literature), due to thermo-convective instabilities or turbulence \cite{poi-05-the, pet-00-tur}.
Whatever the modelling strategy, the problem thus needs a multiscale approach, since the local flame brush structure is out of reach of the computations aimed at simulating the flow dynamics at the observation scale, \ie\ the whole reactive atmosphere scale.
A possible way to completely circumvent this problem is to perform an explicit computation of the flame brush location, solving a transport-like equation for a characteristic function of the burnt zone; such an approach transfers the modelling difficulty to the evaluation of the flame brush velocity (or, more precisely speaking, to the relative velocity of the flame brush with respect to the fresh gases), by an adequate closure relation, and the resulting model is generally referred to as a Turbulent Flame velocity Closure (TFC) model \cite{zim-00-gas}.
The transport equation for the characteristic function of the burnt zone is called in this context the $G$-equation, its unknown being denoted by $G$ \cite{pet-00-tur}.
Such a modelling is implemented in the in-house software P$^2$REMICS (for Partially PREMIxed Combustion Solver) developed, on the basis of the software components library CALIF$^3$S (for Components Adaptative Library For Fluid Flow Simulations, see \cite{califs}) at the French Institut de Radioprotection et S\^uret\'e Nucl\'eaire (IRSN) for safety evaluation purposes; this is the context of the work presented in the present paper.

\medskip
Usually, TFC models apply to perfectly premixed flows (\ie\ flows with constant initial composition), and the chemical state of the flow is governed by the value of $G$ only: $G \in [0,1]$, for $G \geq 0.5$, the mixture is supposed to be in its fresh (initial) state and $G < 0.5$ is supposed to correspond to the burnt state; in both cases, the composition of the gas is known (it is equal to the initial value in the fresh zones, and to the state resulting from a complete chemical reaction in the burnt zone).

\medskip
However, for partially premixed turbulent flows (\ie\ flows with non-constant initial composition), the situation is longer complex, since the composition of the mixture can no more be deduced from the value of $G$.
An extension for this situation, in the inviscid case, is proposed in \cite{bec-10-rea}.
The line followed to formulate this model is to write transport equations for the chemical species initially present in the flow, as if no chemical reaction occured, and then to compute the actual composition in the burnt zone (\ie\ the part of the physical space where $G <0.5$) as the chemical equilibrium composition, thus supposing an infinitely fast reaction.
This model is referred to in the following as the {\em ``asymptotic model"}, and is recalled in the first part of Section \ref{sec:phys}.

\medskip
We propose here an alternate extension, which consists in keeping the classical reactive formulation of the chemical species mass balance, but evaluating the reaction term as a function of $G$: it is set to zero in the fresh zone ($G \geq 0.5$), and to a finite (but possibly large) value in the burnt zone ($G <0.5$).
This model is referred to as the {\em ``relaxed model"}; it is in fact more general, as it may be readily extended to cope with diffusion terms, while the ``asymptotic model" cannot (to this purpose, a balance for the actual mass fractions is necessary).
We then build a numerical scheme, based on a staggered discretization of the unknowns, for the solution of the relaxed model; this algorithm is of fractional step type, and employs a pressure correction technique for hydrodynamics.
The balance energy solved by the scheme is the so-called (non conservative) sensible enthalpy balance, with corrective terms in order to ensure the weak consistency (in the Lax-Wendroff senses) of the scheme.
It enjoys the same stability properties as the continuous model: positivity of the density and, thanks to the choice of the enthalpy balance, the internal energy, conservation of the total energy, chemical species mass fractions lying in the interval $[0,1]$.
In addition, it is shown to be in fact conservative: indeed, its solutions satisfy a discrete conservative total energy balance whose  time and space discretization is non-standard, but weakly consistent with its continuous counterpart.
This algorithm is an extension to the reactive case of the numerical scheme for compressible Navier-Stokes equations described and tested in \cite{gra-17-unc}.

\medskip
As the reaction term  gets stiffer, the relaxed model should boil down to the asymptotic one, for which a closed form of the solution of Riemann problems is available.
Numerical tests  are performed which show that indeed this is the case.
In addition, we observe that the accuracy of the scheme (for this kind of application) is highly dependent on the numerical diffusion introduced by the scheme in the mass balance equation for the chemical species, comparing the results for three approximations of the convection operator in these equations: the standard upwind scheme, a MUSCL-like scheme introduced in \cite{pia-13-for} and a first order scheme designed to reduce diffusion proposed in \cite{dep-01-con}.  

\medskip
The presentation is structured as follows.
We first introduce the asymptotic and the relaxed models in Section \ref{sec:phys}.
Then we give an overview of the content of this paper in Section \ref{sec:gene}, writing the scheme in the time semi-discrete setting and stating its stability and consistency property.
The fully discrete setting is given in two steps, first describing the space discretization (Section \ref{sec:mesh}) and then the scheme itself (Section \ref{sec:scheme}).
The conservativity of the scheme is shown in Section \ref{sec:cons}.
Finally, numerical experiments are presented in Section \ref{sec:num}.
%
%
\section{The physical models}\label{sec:phys}

We begin with the description of the asymptotic model introduced in \cite{bec-10-rea} and then turn to the relaxed model proposed in the present work.

\bigskip
{\bf The asymptotic model} - For the sake of simplicity, only four chemical species are supposed to be present in the flow, namely the fuel (denoted by $F$), the oxydant ($O$), the product ($P$) of the reaction, and a neutral gas ($N$).
A one-step irreversible total chemical reaction is considered, which is written:
\[
\nu_F F + \nu_O O + N \rightarrow \nu_P P + N,
\]
where $\nu_F$, $\nu_O$ and $\nu_P$ are the molar stoichiometric coefficients of the reaction.
We denote by $\mathcal{I}$ the set of the subscripts used to refer to the chemical species in the flow, so $\mathcal{I}=\{F,O,N,P\}$ and the set of mass fractions of the chemical species in the flow reads $\{y_i,\ i \in \mathcal{I}\}$ (\ie\ $\{y_F,\ y_O,\ y_N,\ y_P\}$).
We now define the auxiliary unknowns $\{\tilde y_i,\ i \in \mathcal{I}\}$ as the result of the (inert) transport by the flow of the initial state, which means that the $\{\tilde y_i,\ i \in \mathcal{I}\}$ are the solutions to the following system of equation:
\begin{equation}
\partial_t(\rho \tilde y_i)+\dive(\rho\tilde y_i\bfu)= 0,\quad \tilde y_i(\bfx,0)=y_{i,0}(\bfx) \qquad \mbox{for } i \in \mathcal{I},
\end{equation}
where $\rho$ stands for the fluid density, $\bfu$ for the velocity, and $y_{i,0}(\bfx)$ is the initial mass fraction of the chemical species $i$ in the flow.
These equations are supposed to be posed over a bounded domain $\Omega$ of $\xR^d$, $d\in\{1,2,3\}$ and a finite time interval $(0,T)$.
The initial conditions are supposed to verify $\sum_{i\in \mathcal{I}} y_{i,0} =1$ everywhere in $\Omega$, and this property is assumed to be valid for any $t \in (0,T)$, which is equivalent with the mixture mass balance, given below.
The characteristic function $G$ is supposed to obey the following equation:
\begin{equation}\label{eq:G}
\partial_t(\rho G) + \dive (\rho G \bfu) + \rho_u u_f |\gradi G| =0,
\end{equation}
associated to the initial conditions $G=0$ at the location where the flame starts and $G=1$ elsewhere.
The quantity $\rho_u$ is a constant density, which, from a physical point of view, stands for a characteristic value for the unburnt gases density.
The chemical mass fractions are now computed as:
\begin{equation}\label{eq:y_asymp}
\left| \begin{array}{ll}
\mbox{if } G > 0.5,
& \displaystyle
y_i = \tilde y_i \quad \mbox{for } i \in \mathcal{I},
\\[1ex]
\mbox{if } G \leq 0.5,
& \displaystyle
y_F=\nu_F W_F \tilde z^+,\ y_O=\nu_O W_O \tilde z^-,\ y_N=\tilde y_N,
\mbox{ with } \tilde z = \frac 1 {\nu_F W_F} \tilde y_F -\frac 1 {\nu_O W_O} \tilde y_O.
\end{array} \right.
\end{equation}
In these relation, $\tilde z^+$ and $\tilde z^-$ stand for the positive and negative part of $\tilde z$, respectively, \ie\ $\tilde z^+=\max(\tilde z,0)$ and $\tilde z^-=-\min(\tilde z,0)$, and, for $i \in \mathcal{I}$, $W_i$ is the molar mass of the chemical species $i$.
The physical meaning of Relation \eqref{eq:y_asymp} is that the chemical reaction is supposed to be infinitely fast, and thus that the flow composition is stuck to the chemical equilibrium composition in the so-called burnt zone, which explains why the model is qualified as ``asymptotic".
The product mass fraction is given by $y_P = 1 - (y_F+ y_O+y_N)$.
The flow is governed by the Euler equations:
\begin{subequations}
\begin{align}\label{eq:mass} &
\partial_t \rho + \dive( \rho \bfu) = 0,
\\[1ex] \label{eq:mom} &
\partial_t (\rho u_i) + \dive(\rho u_i \bfu) + \partial_i p = 0,
\quad i=1,d,
\\[1ex] \label{eq:etot} &
\partial_t (\rho E) + \dive(\rho E \bfu) + \dive (p \bfu)= 0,
\\[1ex] \label{eq:state} &
 p=(\gamma-1)\, \rho e_s, \qquad E=\frac 1 2|\bfu|^2+ e, \quad e=e_s + \sum_{i\in\mathcal{I}} y_i \Delta h_{f,i}^0\,
\end{align}\label{cvns:eq:pb} \end{subequations}
where $p$ stands for the pressure, $E$ for the total energy, $e$ for the internal energy, $e_s$ for the so-called sensible internal energy and, for $i \in \mathcal{I}$, $\Delta h_{f,i}^0$ is the formation enthalpy of the chemical species $i$.
The equation of state \eqref{eq:state} supposes that the fluid is a perfect mixture of ideal gases, with the same iso-pressure to iso-volume specific heat ratio $\gamma >1$.
This set of equations is complemented by homogeneous Neumann boundary conditions for the velocity:
\begin{equation}\label{eq:bc}
\bfu \cdot \bfn =0\quad \mbox{a.e. on }\partial \Omega,
\end{equation}
where $\partial\Omega$ stands for the boundary of $\Omega$ and $\bfn$ its outward normal vector.

\bigskip
{\bf The ``relaxed" model} -- This model retains the original form governing equations for reactive flows: a  a transport/reaction equation is written for each of the chemical species mass fractions;   the value of $G$ controls the reaction rate $\dot \omega$, which is set to zero when $G \geq 0.5$, and takes non-zero (and possibly large) values otherwise.
The unknowns $\{y_i,\ i \in \mathcal{I}\}$ are thus now solution to the following balance equations:
\begin{equation}
\partial_t(\rho y_i)+\dive(\rho y_i \bfu)= \dot \omega_i ,\quad \tilde y_i(\bfx,0)=y_{i,0}(\bfx) \qquad \mbox{for } i \in \mathcal{I},
\end{equation}
where the reactive term $\dot \omega_i$ is given by:
\begin{equation}\label{eq:reactive_term}
\dot \omega_i= \frac 1 {\varepsilon}\ \zeta_i\, \nu_i W_i\, \dot \omega,
\quad \dot \omega= \eta(y_F,y_O)\ (G-0.5)^-,\quad \eta(y_F,y_O)=\min (\frac{y_F}{\nu_F W_F}, \frac{y_O}{\nu_O W_O}),
\end{equation}
with $\zeta_F=\zeta_O=-1$, $\zeta_P=1$ and $\zeta_N=0$.
Note that, since $\nu_F W_F + \nu_O W_0 = \nu_P W_P$, we have $\sum_{i \in \mathcal{I}} \dot \omega_i=0$, which, summing on $i \in \mathcal{I}$ the species mass balance, allows to recover the equivalence between the mass balance and the fact that $\sum_{i \in \mathcal{I}} y_i=1$.
The factor $\eta(y_F,y_O)$ is a cut-off function, which prevents the chemical species mass fractions from taking negative values (and, consequently, values greater than 1, since their sum is equal to 1).

\medskip
The rest of the model is left unchanged.
%
%
\section{General description of the scheme and main results}\label{sec:gene}

\paragraph{\textbf{Time semi-discrete algorithm}}
Instead of the total energy balance equation, the scheme solves a balance equation for the sensible enthalpy $h_s = e_s+p/\rho$, which is formally derived as follows.
The first step is to establish the kinetic energy balance formally and subtract from \eqref{eq:etot} to obtain a balance equation for the internal energy.
Thanks to the mass balance equation, for any regular function $\psi$
\[
\partial_t(\rho \psi) + \dive(\rho \psi \bfu) = \rho \partial_t \psi + \rho \bfu \cdot \gradi \psi.
\]
Using twice this identity and then the momentum balance equation, we have for $1 \leq i \leq d$:
\[
\frac 1 2 \partial_t(\rho u_i^2) + \frac 1 2 \dive(\rho u_i^2\, \bfu)
= \rho\, u_i \partial_t u_i + \rho u_i \bfu \cdot \gradi u_i
= u_i \bigl[ \partial_t(\rho u_i) + \dive(\rho u_i \bfu) \bigr]
= -u_i \partial_i p,
\]
and, summing for $i=1$ to $d$, we obtain the kinetic energy balance:
\[
\frac 1 2 \partial_t (\rho |\bfu|^2)+\frac 1 2 \dive(\rho |\bfu|^2 \bfu )
= \bfu \cdot \bigl[\partial_t(\rho \bfu) + \dive(\rho \bfu \otimes \bfu) \bigr]
= -\bfu \cdot \gradi p.
\]
Substituting the expression of the total energy in \eqref{eq:etot}, yields
\[
\partial_t(\rho e) + \dive(\rho e \bfu) + \frac 1 2 \partial_t(\rho |\bfu|^2) + \frac 1 2 \dive(\rho |\bfu|^2) + \bfu\cdot\gradi p + p \dive(\bfu) = 0,
\]
which, using the kinetic energy balance, gives the total internal energy balance:
\begin{equation} \label{eq:e}
\partial_t(\rho e) + \dive(\rho e \bfu) + p \dive(\bfu) = 0.
\end{equation}
Using the linearity of the mass balance of the chemical species $i$, for any $i\in\mathcal I$, we derive the reactive energy balance:
\begin{equation} \label{eq:er}
\partial_t \bigl[\rho \big(\sum_{i \in \mathcal{I}} \Dhi y_i\big) \bigr] + \dive\bigl[ \rho \big(\sum_{i \in \mathcal{I}} \Dhi y_i\big) \bfu \bigr]
= \sum_{i \in \mathcal{I}} \Dhi \dot\omega_i = -\dot\omega_\theta.
\end{equation}
Subtracting \eqref{eq:er} from \eqref{eq:e} yields the sensible internal energy balance:
\begin{equation} \label{eq:es}
\partial_t(\rho e_s) + \dive(\rho e_s \bfu) + p \dive(\bfu) = \dot\omega_\theta.
\end{equation}
Finally, using the relation between the sensible energy and the sensible enthalpy, we obtain the sensible enthalpy balance:
\begin{equation} \label{eq:hs}
 \partial_t(\rho h_s) + \dive(\rho h_s\bfu) - \partial_t p - \bfu\cdot\gradi p = \dot\omega_\theta.
\end{equation}

The numerical resolution of the mathematical model is realized by a fractional step algorithm, which implements a pressure correction technique for hydrodynamics in order to separate the resolution of the momentum balance from the other equations of the Euler system.
Supposing that the time interval $(0,T)$ is split in $N$ sub-intervals, of constant length $\delta t=T/N$, the semi-discrete algorithm is given by:
\begin{subequations}\label{eq:sd_scheme}
\begin{align} \nonumber &
\mbox{Reactive step:}
\\ \label{eq:Gsd} & \quad G^{n+1}:
&&
\frac{1}{\delta t}(\rho^n G^{n+1}-\rho^{n-1}G^{n})+\dive(\rho^n G^k \bfu^n) + \rho_u u_f\, |\gradi G^{n+1}| = 0,
\\[0.5ex]\label{eq:yN} & \quad Y_N^{n+1}:
&&
\frac{1}{\delta t} (\rho^n y_N^{n+1}-\rho^{n-1} y_N^n) + \dive(\rho^n y_N^k \bfu^n) = 0.
\\[0.5ex]\label{eq:z} & \quad z^{n+1}:
&&
\frac{1}{\delta t}(\rho^n z^{n+1}-\rho^{n-1} z^n)+\dive(\rho^n z^k \bfu^n) = 0.
\\[0.5ex]\label{eq:yF} & \quad Y_F^{n+1}:
&&
\frac{1}{\delta t} (\rho^n y_F^{n+1}-\rho^{n-1} y_F^n) + \dive(\rho^n y_F^k \bfu^n) =
-\frac 1 {\varepsilon} \nu_F W_F\, \dot{\omega}(y_F^{n+1},z^{n+1}),
\\[0.5ex]\label{eq:yP} & \quad Y_P^{n+1}:
&&
y_F^{n+1}+y_O^{n+1}+y_N^{n+1}+y_P^{n+1}=1.
\displaybreak[1] \\[3ex] \nonumber &
\mbox{Euler step:}
\\[0.5ex]\label{eq:pred} & \quad \tilde \bfu^{n+1}:
&&
\begin{array}{l} \displaystyle
\frac{1}{\delta t}(\rho^n  \tilde u^{n+1}_i-\rho^{n-1} u_i^n) + \dive(\rho^n  \tilde u_i^{n+1} \bfu^n)
\\ \displaystyle \hspace{32ex}
+ \Big( \frac{\rho^n}{\rho^{n-1}} \Big)^{1/2} \partial_i p^n = 0, \quad  i=1,\dots,d,
\end{array}
\displaybreak[1] \\[3ex]\label{eq:cor} & \quad
\begin{array}{l} \bfu^{n+1},\ \rho^{n+1},\\ h_s^{n+1},\ p^{n+1}:\end{array}
&&
\left| \begin{array}{l} \displaystyle
\dfrac{1}{\delta t}\ \rho^n (u^{n+1}_i-\tilde u_i^{n+1}) + \partial_i p^{n+1}
- \Big( \frac{\rho^n}{\rho^{n-1}} \Big)^{1/2} \partial_i p^n = 0, \quad i=1,\dots,d,
\\[2ex] \displaystyle
\dfrac{1}{\delta t}(\rho^{n+1}-\rho^n) + \dive(\rho^{n+1} \bfu^{n+1})=0, 
\\[2ex] \displaystyle
\frac{1}{\delta t}\,(\rho^{n+1} h_s^{n+1}-\rho^n h_s^n)
+ \dive( \rho^{n+1} h_s^{n+1} \bfu^{n+1})
-\dfrac{1}{\delta t}\,(p^{n+1}- p^n)
\\[2ex] \displaystyle \hfill
- u^{n+1} \cdot \gradi p^{n+1}
= \dot{\omega}_\theta^{n+1} + S^{n+1},
\\[2ex] \displaystyle
p^{n+1} = \frac {\gamma -1}\gamma\ \rho^{n+1}\, h_s^{n+1}.
\end{array} \right.
\end{align}
\end{subequations}
Equations \eqref{eq:Gsd}-\eqref{eq:cor} are solved successively, and the unknown for each equation is specified before each equation.
In the convection term of the equations of the reactive step, the index $k$ may take the value $n$ (so the scheme is explicit) or $n+1$ (so the scheme is implicit).
The unknown $z$ is an affine combination of $y_F$ and $y_O$, defined so that the reactive term cancels:
\begin{equation}\label{eq:def_z}
z = \frac 1 {\nu_F W_F} y_F -\frac 1 {\nu_O W_O} y_O.
\end{equation}
Thus the value of $y_O^{n+1}$ is deduced from $y_F^{n+1}$ and $z^{n+1}$, which allows to express $\dot{\omega}$ in \eqref{eq:yF} as a function of $y_F^{n+1}$ and $z^{n+1}$, instead of $y_F^{n+1}$ and $y_O^{n+1}$ as suggested by Relation \eqref{eq:reactive_term}.
In addition, we have:
\[
\eta(y_F^{n+1},y_O^{n+1})=\min (\frac{y_F^{n+1}}{\nu_F W_F}, \frac{y_O^{n+1}}{\nu_O W_O}) = 
\left| \begin{array}{ll} \displaystyle
\frac 1 {\nu_F W_F}\, y_F^{n+1}
&
\mbox{if } z^{n+1} \leq 0,
\\[4ex] \displaystyle
\frac 1 {\nu_O W_O}\, y_O^{n+1} = \frac 1 {\nu_F W_F}\, y_F^{n+1} - z^{n+1}
&
\mbox{otherwise.}
\end{array} \right.
\]
Hence, because of the specific form of the function $\eta$, the right hand side of \eqref{eq:yF} boils down to an affine term, even if $\eta$ vanishes when $y_F$ or $y_O$ vanishes, and the scheme is fully implicit in time with respect to the reaction term.
This is the motivation for the choice of the form of $\eta$.
It is fundamental to remark that Equations \eqref{eq:yN}-\eqref{eq:yP} are equivalent to the following system:
\begin{equation} \label{eq:chem_sys}
\frac{1}{\delta t} (\rho^n y_i^{n+1}-\rho^{n-1} y_i^n) + \dive(\rho^n y_i^k \bfu^n) =
\frac 1 {\varepsilon} \zeta_i \nu_i W_i\, \dot{\omega}(y_F^{n+1},y_O^{n+1}), \quad \mbox{for } i \in \mathcal{I},
\end{equation}
where we recall that $\zeta_F=\zeta_O=-1$, $\zeta_P=1$ and $\zeta_N=0$.
Indeed, dividing the fuel mass balance equation \eqref{eq:yF} by $\nu_F W_F$, substracting Equation \eqref{eq:z} and finally multiplying by $\nu_O W_O$ yields the desired mass balance equation for the oxydant chemical species.
Finally, we suppose that the product mass balance holds:
\begin{equation} \label{eq:yP_bis}
\frac{1}{\delta t} (\rho^n y_P^{n+1}-\rho^{n-1} y_P^n) + \dive(\rho^n y_P^k \bfu^n) =
\frac 1 {\varepsilon} \nu_P W_P\, \dot{\omega}(y_F^{n+1},y_O^{n+1}).
\end{equation}
Since the sum of the chemical reaction terms vanishes, we have for $\Sigma= y_F + y_O + y_P +y_N$, summing all the chemical species mass balances,
\begin{equation} \label{eq:Sigma}
\frac{1}{\delta t} (\rho^n \Sigma^{n+1}-\rho^{n-1} \Sigma^n) + \dive(\rho^n \Sigma^k \bfu^n) =0,
\end{equation}
and this equation may equivalently replace the product mass balance equation \eqref{eq:yP_bis}.
Thanks to the mixture balance, we see that, provided that $\Sigma^n$ satisfies $\Sigma^n=1$ everywhere in $\Omega$, the solution to Equation \eqref{eq:Sigma} is $\Sigma^{n+1}=1$ everywhere in $\Omega$.
Since the initialization yields $\Sigma^0=1$, this last equality is indeed true, and \eqref{eq:yP_bis} is equivalent to \eqref{eq:yP}.
Finally, note that, when the chemical step is performed, the mass balance at step $n+1$ is not yet solved; hence the (unusual) backward time shift for the densities and for the mass fluxes in the equations of this step.

\medskip
Equations \eqref{eq:pred}-\eqref{eq:cor} implement a pressure correction technique, where the correction step couples the velocity correction equation, the mass balance and the sensible enthalpy balance.
This coupling ensures that the pressure and velocity are kept constant through the contact discontinuity associated to compositional non-reactive Euler equations (precisely speaking, the usual contact discontinuity, already present in 1D equations, but not slip lines); for this property to hold, it is necessary that all chemical species share the same heat capacity ratio $\gamma$.
The term $S_K^{n+1}$ in the sensible enthalpy balance equation is a corrective term which is necessary for consistency; schematically speaking, it compensates the numerical dissipation which appears in a discrete kinetic energy balance that is obtained from the discrete momentum balance.
Its expression is given in Section \ref{sec:scheme}, and its derivation is explained in Section \ref{sec:cons}, where the conservativity of the scheme is discussed.

\bigskip
\paragraph{\textbf{Space discretization}}
The space dicretization is performed by a finite volume technique, using a staggered arrangement of the unknowns (the scalar variables are approximated at the cell centers and the velocity components at the face centers), using either a MAC scheme (for structured discretizations) or the degrees of freedom of low-order non-conforming finite elements: Crouzeix-Raviart \cite{cro-73-con} for simplicial cells and Rannacher-Turek \cite{ran-92-sim} for quadrangles ($d=2$) or hexahedra ($d=3$).
For the Euler equations (\ie\ Steps \eqref{eq:pred}- \eqref{eq:cor}), upwinding is performed by building positivity-preserving convection operators, in the spirit of the so-called Flux-Splitting methods, and only first-order upwinding is implemented.
The pressure gradient is built as the transpose (with respect to the $L^2$ inner product) of the natural velocity divergence operator.
For the balance equations for the other scalar unknowns, the time discretization is implicit when first-order upwinding is used in the convection operator (in other words, $k=n+1$ in \eqref{eq:Gsd}-\eqref{eq:yF}) or explicit ($k=n$ in \eqref{eq:Gsd}-\eqref{eq:yF}) when a higher order (of MUSCL type, \cf\ Section \ref{an:MUSCL}) flux or an anti-diffusive flux (\cf\ Section \ref{an:AD}) is used.

\bigskip
\paragraph{\textbf{Properties of the scheme}}
First, the positivity of the density is ensured by construction of the discrete mass balance equation,\ie\ by the use of a first order upwind scheme.
In addition, the physical bounds of the mass fractions are preserved thanks to the following (rather standard) arguments: first, building a discrete convection operator which vanishes when the convected unknown is constant thanks to the discrete mass balance equation ensures a positivity-preservation property \cite{lar-91-how}, under a CFL condition if an explicit time approximation is used; second, the discretization of the chemical reaction rate ensures either that it vanishes when the unknown of the equation vanishes (for $y_F$ and $y_O$), or that it is non-negative (for $y_P$).
Consequently, mass fractions are non-negative and, since their sum is equal to $1$ (see above), they are also bounded by $1$.

\medskip
The positivity of the sensible energy stems from two essential arguments: first, a discrete analog of the internal energy equation \eqref{eq:e} may be obtained from the discrete sensible enthalpy balance, by mimicking the continuous computation; second, this discrete relation may be shown to have only positive solutions, once again thanks to the consistency of the discrete convection operator and the mass balance.
This holds provided that the equation is exothermic ($\dot{\omega}_\theta \geq 0$) and thanks to the non-negativity of $S^{n+1}$ (see below).

\medskip
In order to calculate correct shocks, it is crucial for the scheme to be consistent with the following weak formulation of the problem:
\begin{equation}
\label{eq:weak_form}
\!
\begin{array}{ll}
& \displaystyle
\forall\phi\in C_c^\infty( \Omega\times[0,T) \big),
\\[2ex] & \displaystyle \quad 
\int_0^T\!\int_\Omega\!\big[\rho\partial_t\phi + \rho\bfu\cdot\gradi\phi\big]{\rm d}\bfx\,{\rm d}t 
+ \int_\Omega\!\rho_0(\bfx)\phi(\bfx,0){\rm d}\bfx = 0,
\\[2ex] & \displaystyle \quad 
\int_0^T\!\int_\Omega\!\big[\rho u_i\partial_t\phi + (\rho\bfu u_i)\cdot\gradi\phi + p\partial_i\phi\big]{\rm d}\bfx\,{\rm d}t 
+ \int_\Omega\!\rho_0(\bfx)(u_i)_0(\bfx)\phi(\bfx,0){\rm d}\bfx = 0, \quad 1\leq i\leq d,
\\[2ex] & \displaystyle \quad 
\int_0^T\!\int_\Omega\!\big[ \rho E\partial_t\phi + (\rho E+p)\bfu\cdot\gradi\phi \big]{\rm d}\bfx\,{\rm d}t 
+ \int_\Omega\! \rho_0(\bfx) E_0(\bfx) \phi(\bfx,0) {\rm d}\bfx = 0,
\\[2ex] & \displaystyle \quad 
\int_0^T\!\int_\Omega\!\big[ \rho y_i \partial_t \phi + \rho y_i \bfu \cdot \gradi\phi \big] {\rm d}\bfx\,{\rm d}t + \int_0^T\!\int_\Omega\! \rho_0(\bfx) y_{i,0}(\bfx) \phi(\bfx,0) {\rm d}\bfx = - \int_0^T\!\int_\Omega\! \dot\omega_i \phi\,{\rm d}\bfx\,{\rm d}t, \quad 1\leq i\leq d,
\\[2ex] & \displaystyle \quad 
p = (\gamma-1) \rho e_s.
\end{array}
\end{equation}

\medskip
Remark that this system features the total energy balance equation and not the sensible enthalpy balance equation, which is actually solved here.
However, we show in Section \ref{sec:cons} that the solutions of the scheme satisfy a discrete total energy balance, with a time and space dicretization which is unusual but allows however to prove the consistency in the Lax-Wendroff sense.
Finally, the integral of the total energy over the domain is conserved, which yields a stability result for the scheme (irrespectively of the time and space step, for this relation; recall however that the overall stability of the scheme needs a CFL condition if an explicit version of the convection operator for chemical species is used).
%
%
\section{Meshes and unknowns}\label{sec:mesh}

Let the computational domain $\Omega$ be an open polygonal subset of $\xR^d$, $1 \leq d \leq 3$, with boundary $\partial \Omega$ and  let $\mesh$ be a decomposition of $\Omega$, supposed to be regular in the usual sense of the finite element literature (\eg\ \cite{cia-91-bas}).
The cells may be:
\begin{list}{-}{\itemsep=0.5ex \topsep=0.5ex \leftmargin=1.cm \labelwidth=0.3cm \labelsep=0.5cm \itemindent=0.cm}
\item for a general domain $\Omega$, either convex quadrilaterals ($d=2$) or hexahedra ($d=3$) or simplices, both type of cells being possibly combined in a same mesh,
\item for a domain the boundaries of which are hyperplanes normal to a coordinate axis, rectangles ($d=2$) or rectangular parallelepipeds ($d=3$) (the faces of which, of course, are then also necessarily normal to a coordinate axis).
\end{list}
By $\edges$ and $\edges(K)$ we denote the set of all $(d-1)$-faces $\edge$ of the mesh and of the element $K \in \mesh$ respectively.
The set of faces included in the boundary of $\Omega$ is denoted by $\edgesext$ and the set of internal edges (\ie\ $\edges \setminus \edgesext$) is denoted by $\edgesint$; a face $\edge \in \edgesint$ separating the cells $K$ and $L$ is denoted by $\edge=K|L$.
The outward normal vector to a face $\edge$ of $K$ is denoted by $\bfn_{K,\edge}$.
For $K \in \mesh$ and $\edge \in \edges$, we denote by $|K|$ the measure of $K$ and by $|\edge|$ the $(d-1)$-measure of the face $\edge$.
For any $K\in\mesh$ and $\edge\in\edges(K)$, we denote by $d_{K,\edge}$ the Euclidean distance between the center $x_K$ of the mesh and  the edge $\edge$.
For any $\edge\in\edges$, we define $d_{\edge}=d_{K,\edge}+d_{L,\edge}$, if $\edge \in \edgesint$ and $d_{\edge}=d_{K,\edge}$ if $\edge\in\edgesext$.
The size of the mesh is denoted by $h$.
For $1 \leq i \leq d$, we denote by $\edges \ei \subset \edges$ and $\edgesext \ei \subset \edgesext$ the subset of the faces of $\edges$ and $\edgesext$ respectively which are perpendicular to the $i^{th}$ unit vector of the canonical basis of $\xR^d$.

\medskip
The space discretization is staggered, using either the Marker-And Cell (MAC) scheme \cite{har-65-num,har-71-num}, or nonconforming low-order finite element approximations, namely the Rannacher and Turek (RT) element \cite{ran-92-sim} for quadrilateral or hexahedric meshes, or the lowest degree Crouzeix-Raviart (CR) element \cite{cro-73-con} for simplicial meshes.

\medskip
For all these space discretizations, the degrees of freedom for the pressure, the density, the enthalpy, the mixture, fuel and neutral gas mass fractions and the flame indicator are associated to the cells of the mesh $\mesh$
 and are denoted by:
\[
 \bigl\{ p_K,\ \rho_K,\ h_K,\ y_{F,K},\ y_{N,K},\ z_K,\ G_K,\ K \in \mesh \bigr\}.
\]

\medskip
Let us then turn to the degrees of freedom for the velocity (\ie\ the discrete velocity unknowns).
\begin{list}{-}{\itemsep=0.5ex \topsep=0.5ex \leftmargin=1.cm \labelwidth=0.3cm \labelsep=0.5cm \itemindent=0.cm}
\item {\bf Rannacher-Turek} or {\bf Crouzeix-Raviart} discretizations --
The degrees of freedom for the velocity components are located at the center of the faces of the mesh, and we choose the version of the element where they represent the average of the velocity through a face.
The set of degrees of freedom reads:
\[
\lbrace \bfu_{\edge},\ \edge \in \edges\rbrace,\mbox{ of components }\lbrace u_{\edge,i},\ \edge \in \edges,\ 1 \leq i \leq d \rbrace.
\]
\item {\bf MAC} discretization -- The degrees of freedom for the $i^{th}$ component of the velocity are defined at the centre of the faces of $\edges \ei$, so the whole set of discrete velocity unknowns reads:
\[
\big\{ u_{\edge,i},\ \edge \in \edges\ei, \ 1 \leq i \leq d \big\}.
\]
\end{list}

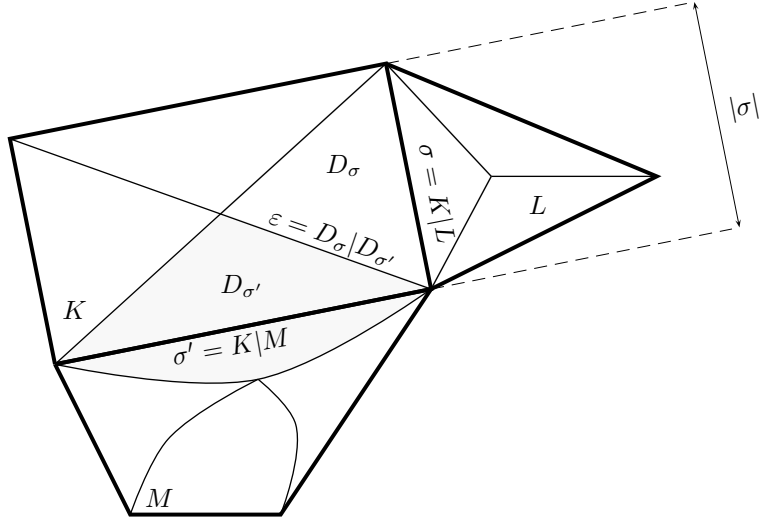
\begin{figure}[!t]
\begin{center}
\newgray{grayml}{.9}
\newgray{grayc}{.97}
\psset{unit=1cm}
\begin{pspicture}(0,0)(12,7)
\rput[bl](0,1){
   \pspolygon[linecolor=grayml](3.2,3)(6,2)(6.8,3.5)(5.4,5)
   \rput[bl](4.6,3.5){{$D_\edge$}}
   \psccurve[fillstyle=solid, fillcolor=grayc, linecolor=grayc](3.2,3)(3.2,3)(6,2)(6,2)(3.7,0.8)(1,1)(1,1)
   \rput[bl](3.2,1.9){{$D_{\edge'}$}}
   \rput[bl]{10}(2.6,0.9){$\edge'=K|M$}
   \psline[linecolor=black, linewidth=1.5pt]{-}(1,1)(6,2)(5.4,5)(0.4,4)(1,1)
   \psline[linecolor=black, linewidth=0.5pt]{-}(1,1)(5.4,5)
   \psline[linecolor=black, linewidth=0.5pt]{-}(6,2)(0.4,4)
   \rput[bl](1.1,1.6){{$K$}}
   \psline[linecolor=black, linewidth=1.5pt]{-}(6,2)(9,3.5)(5.4,5)
   \psline[linecolor=black, linewidth=0.5pt]{-}(6,2)(6.8,3.5)
   \psline[linecolor=black, linewidth=0.5pt]{-}(9,3.5)(6.8,3.5)
   \psline[linecolor=black, linewidth=0.5pt]{-}(5.4,5)(6.8,3.5)
   \rput[bl](7.3,3){{$L$}}
   \psline[linecolor=black, linewidth=1.5pt]{-}(1,1)(6,2)(4,-1)(2,-1)(1,1)
   \pscurve[linecolor=black, linewidth=0.5pt]{-}(4,-1)(4.2,0.2)(3.7,0.8)
   \pscurve[linecolor=black, linewidth=0.5pt]{-}(2,-1)(2.5,0.)(3.7,0.8)
   \pscurve[linecolor=black, linewidth=0.5pt]{-}(1,1)(3.7,0.8)(6,2)
   \rput[bl](2.2,-0.9){{$M$}}
   \psline[linecolor=black, linewidth=0.3pt, linestyle=dashed]{-}(6,2)(10,2.8)
   \psline[linecolor=black, linewidth=0.3pt, linestyle=dashed]{-}(5.4,5)(9.4,5.8)
   \psline[linecolor=black, linewidth=0.3pt]{<->}(10.1,2.82)(9.5,5.82) \rput[bl]{11}(10,4.2){$|\edge|$}
   \rput[bl]{-79}(5.75,3.9){$\edge=K|L$}
   \rput[bl]{-19}(3.8,2.8){$\edged=D_\edge|D_{\edge'}$}
}
\end{pspicture}
\caption{Primal and dual meshes for the Rannacher-Turek and Crouzeix-Raviart elements.}
\label{fig:mesh}
\end{center}
\end{figure}

\medskip
For the definition of the schemes, we need a dual mesh which is defined as follows.
\begin{list}{-}{\itemsep=0.5ex \topsep=0.5ex \leftmargin=1.cm \labelwidth=0.3cm \labelsep=0.5cm \itemindent=0.cm}
\item {\bf Rannacher-Turek} or {\bf Crouzeix-Raviart} discretizations --
For the RT or CR discretizations, the dual mesh is the same for all the velocity components.
When $K\in\mesh$ is a simplex, a rectangle or a rectangular cuboid, for $\edge \in \edges(K)$, we define $D_{K,\edge}$ as the cone with basis $\edge$ and with vertex the mass center of $K$ (see Figure \ref{fig:mesh}).
We thus obtain a partition of $K$ in $m$ sub-volumes, where $m$ is the number of faces of the mesh, each sub-volume having the same measure $| D_{K,\edge}|= |K|/m$.
We extend this definition to general quadrangles and hexahedra, by supposing that we have built a partition still of equal-volume sub-cells, and with the same connectivities; note that this is of course always possible, but that such a volume $D_{K,\edge}$ may be no longer a cone; indeed, if $K$ is far from a parallelogram, it may not be possible to build a cone having $\edge$ as basis, the opposite vertex lying in $K$ and a volume equal to $|K|/m$ (note that these dual cells do not need to be constructed in the implementation of the scheme, only their volume is needed).
The volume $D_{K,\edge}$ is referred to as the half-diamond cell associated to $K$ and $\edge$.\\
For $\edge \in \edgesint$, $\edge=K|L$, we now define the diamond cell $D_\edge$ associated to $\edge$ by $D_\edge=D_{K,\edge} \cup D_{L,\edge}$; for an external face $\edge \in \edgesext \cap \edges(K)$, $D_\edge$ is just the same volume as $D_{K,\edge}$.
\item {\bf MAC} discretization --
For the MAC scheme, the dual mesh depends on the component of the velocity.
For each component, the MAC dual mesh only differs from the RT or CR dual mesh by the choice of the half-diamond cell, which, for $K \in \mesh$ and $\edge \in \edges(K)$, is now the rectangle or rectangular parallelepiped of basis $\edge$ and of measure $| D_{K,\edge}|= |K|/2$.
\end{list}

\medskip
We denote by $|D_\edge|$ the measure of the dual cell $D_\edge$, and by $\edged=D_\edge|D_{\edge'}$ the dual face separating two diamond cells $D_\edge$ and $D_{\edge'}$.

\medskip
In order to be able to write a unique expression of the discrete equations for both MAC and CR/RT schemes, we introduce the set of faces $\edgesischeme$ associated with the degrees of freedom of each component of the velocity ($\scheme$ stands for ``scheme''):
\[
\edgesischeme= \left| \begin{array}{ll}
\edges\ei \setminus \edgesext\ei \mbox{ for the MAC scheme},
\\
\edges \setminus \edgesext\ei \mbox{ for the CR or RT schemes.}
\end{array} \right.
\]
Similarly, we unify the notation for the set of dual faces for both schemes by defining:
\[
\edgesdischeme= \left| \begin{array}{ll}
\edgesd\ei \setminus \edgesdext\ei \mbox{ for the MAC scheme},
\\
\edgesd \setminus \edgesdext\ei \mbox{ for the CR or RT schemes,}
\end{array} \right.
\]
where the symbol $\tilde ~$ refers to the dual mesh; for instance, $\edgesd\ei$ is thus the set of faces of the dual mesh associated with the $i^{th}$ component of the velocity, and $\edgesdext\ei$ stands for the subset of these dual faces included in the boundary.
Note that, for the MAC scheme, the faces of $\edgesd\ei$ are perpendicular to a unit vector of the canonical basis of $\xR^d$, but not necessarily to the $i^{th}$ one.

%
%
\section{The scheme}\label{sec:scheme}

In this section, we give the fully discrete form of the scheme.
Even if it corresponds to the reverse order with respect to the semi-discrete scheme given in \eqref{eq:sd_scheme}, we begin with hydrodynamics (Section \ref{subsec:euler}) and then turn to the mass balance step for chemical species and the transport of the characteristic function for the burnt zone (Section \ref{subsec:chem}).
This choice is due to the fact that the definition of the convection operators for scalar variables necessitates to introduce first the discretization of the mixture mass balance equation.

\subsection{Euler step} \label{subsec:euler}

For $0 \leq n <N$, the step $n+1$ of the algorithm for the resolution of the Euler equations reads:
\begin{subequations} \label{scheme:euler}
\begin{align}
& \nonumber \quad \mbox{\textit{Pressure gradient scaling step} -- Solve for $(\widetilde{\gradi p})^{n+1}$:} \\ 
\label{scheme:gradp} & \qquad \forall \edge \in \edges, 
\qquad (\widetilde {\gradi p})^{n+1}_\edge = \Big( \frac{\rho^n_\Ds}{\rho^{n-1}_\Ds} \Big)^{1/2} (\gradi p)^{n}_\edge. \\[2ex]
& \nonumber \quad \mbox{\textit{Prediction step} -- Solve for $\tilde \bfu^{n+1}$:} \\
& \nonumber \qquad \mbox{For } 1 \leq i \leq d,\ \forall \edge \in \edgesischeme, \\
& \hspace{2.2cm} \label{scheme:mom} \qquad \dfrac{1}{\delta t}(\rho^n_\Ds \tilde u^{n+1}_{\edge,i}-\rho^{n-1}_\Ds u_{\edge,i}^n) + \dive_\edge(\rho^n \tilde u_i^{n+1} \bfu^n) + (\widetilde {\gradi p})^{n+1}_{\edge,i} = 0. \\[2ex]
& \nonumber \quad \mbox{\textit{Correction step} -- Solve for $\rho^{n+1}$, $p^{n+1}$ and $\bfu^{n+1}$:} \\
& \nonumber \qquad \mbox{For } 1 \leq i \leq d,\ \forall \edge \in \edgesischeme, \\
& \hspace{2.2cm} \label{scheme:cor} \qquad \dfrac{1}{\delta t}\ \rho^n_\Ds\ (u^{n+1}_{\edge,i}-\tilde u_{\edge,i}^{n+1}) + (\gradi p)_{\edge,i}^{n+1} - (\widetilde{\gradi p})_{\edge,i}^{n+1} = 0, \\ 
\label{scheme:mass} & \qquad \forall K \in \mesh, \qquad
\dfrac{1}{\delta t}(\rho^{n+1}_K-\rho^n_K) + \dive_K(\rho \bfu)^{n+1}=0, \\
& \begin{matrix} 
\qquad \forall K \in \mesh, 
\label{scheme:hs} \qquad 
\displaystyle \frac{1}{\delta t}\,\left[\rho^{n+1}_K\,(h_s)^{n+1}_K-\rho^{n}_K\,(h_s)^n_K\right] + \dive_K(\rho h_s \bfu)^{n+1} \hfill \\
\displaystyle \hspace{7cm} -\frac{1}{\delta t}\,(p^{n+1}_K- p^n_K) - \big(u \cdot \gradi p\big)_K^{n+1} = (\dot{\omega}_\theta)_K^{n+1} + S_K^{n+1},
\end{matrix}
\\ \label{scheme:eos} & \qquad
\forall K \in \mesh, \qquad
p_K^{n+1}=\frac{\gamma-1}{\gamma}(h_s)_K^{n+1}\rho_K^{n+1}. 
\end{align}
\end{subequations}

\medskip
The initial approximations for $\rho^{-1}$, $h_s^0$ and $\bfu^0$ are given by the mean values of the initial conditions over the primal and dual cells:
\[
\begin{array}{l}
\displaystyle
\forall K\in\mesh, \quad \rho_K^{-1} = \frac1{|K|} \int_K\! \rho_0(\bfx) {\rm d}\bfx \quad \mbox{and} \quad (h_s)_K^0 = \frac1{|K|} \int_K\! (h_s)_0(\bfx),  
\\[3ex] \displaystyle
\forall\edge\in\edgesischeme,\ 1\leq i\leq d,\quad u_{\edge,i}^0 = \frac1{|D_\edge|} \int_{D\edge}\! (\bfu_0(\bfx))_i {\rm d}\bfx.
\end{array}
\]
Then, $\rho^0$ is computed by the mass balance equation \eqref{scheme:mass} and $p^0$ is computed by the equation of state \eqref{scheme:eos}.

\medskip
We now define each of the discrete operators featured in System \eqref{scheme:euler}.
%
%

\medskip
\paragraph{\textbf{Mass balance equation}} Equation \eqref{scheme:mass} is a finite volume discretisation of the mass balance \eqref{eq:mass} over the primal mesh.
For a discrete density field $\rho$ and a discrete velocity field $\bfu$, the discrete divergence is defined by:
\[
 \dive_K(\rho\bfu) = \frac{1}{|K|}\sum_{\edge\in\edges(K)} F_{K,\edge}, \quad F_{K,\edge} = |\edge|\ \rho_\edge u_{K,\edge},
\]
where $u_{K,\edge}$ is an approximation of the normal velocity to the face $\edge$ outward $K$.
The definition of this latter quantity depends on the discretization: in the MAC case, $u_{K,\edge} = u_{\edge,i}\ \bfe^{(i)}\cdot\bfn_{K,\edge}$ for a face $\edge$ of $K$ perpendicular to $\bfe^{(i)}$, with $\bfe^{(i)}$ the $i$-th vector of the orthonormal basis of $\mathbb R^d$, and, in the CR and RT cases, $u_{K,\edge} =\bfu_\edge\cdot\bfn_{K,\edge}$ for any face $\edge$ of $K$.
The density at the face $\edge=K|L$ is approximated by the upwind technique, so $\rho_\edge=\rho_K$ if $u_{K,\edge} \geq 0$ and $\rho_\edge=\rho_L$ otherwise.
Since we assume that the normal velocity vanishes on the boundary faces, the definition is complete.
%
%

\medskip
\paragraph{\textbf{Convection operators associated to the primal mesh}} We may now define the discrete convection operator of any discrete field $z$ defined on the primal cell by
\[
\dive_K(\rho z \bfu) = \frac{1}{|K|} \sum_{\edge\in\edges(K)} F_{K,\edge}\ z_\edge,
\]
where $z_\edge$ is the upwind approximation with respect to the mass flux $F_{K,\edge}$ at the face $\edge$.
%
%

\medskip
\paragraph{\textbf{Momentum balance equation and pressure gradient scaling}}
We now turn to the discrete momentum balance \eqref{scheme:mom}.
For the MAC discretization, but also for the RT and CR discretizations, the time derivative and convection terms are approximated in \eqref{scheme:mom} by a finite volume technique over the dual cells, so the convection term reads:
\[
\dive_\edge(\rho \tilde u_i \bfu) = \dive_\edge \bigl( \tilde u_i (\rho  \bfu) \bigr)=
\frac 1 {|D_\edge|} \sum_{\edged \in\edgesd(D_\edge)} F_{\edge,\edged} \tilde u_{\edged,i},
\]
where $F_{\edge,\edged}$ stands for a mass flux through the dual face $\edged$, and $\tilde u_{\edged,i}$ is a centered approximation of the $i^{th}$ component of the velocity $\tilde \bfu$ on $\edged$.
The density at the dual cell $\rho_\Ds$ is obtained by a weighted average of the density in the neighbour cells:
$ |D_\edge|\, \rho_\Ds = |D_{K,\edge}|\, \rho_K + |D_{L,\edge}|\, \rho_L$ for $\edge = K|L \in \edgesint$, and $\rho_\Ds = \rho_K$ for an external face of a cell $K$. 
The mass fluxes $(F_{\edge,\edged})_{\edged \in \edges(D_\edge)}$ are evaluated as linear combinations, with constant coefficients, of the primal mass fluxes at the neighbouring faces, in such a way that the following discrete mass balance over the dual cells is implied by the discrete mass balance \eqref{scheme:mass}:
\begin{equation}\label{eq:mass_D}
\forall \edge\in\edges \mbox{ and } n\in\mathbb N,\qquad
\frac{|D_\edge|}{\delta t} \ (\rho^{n+1}_\Ds-\rho^n_\Ds)
+ \sum_{\edged\in\edges(D_\edge)} F_{\edge,\edged}^{n+1}=0.
\end{equation}
This relation is critical to derive a discrete kinetic energy balance (see Section \ref{sec:cons} below).
The computation of the dual mass fluxes is such that the flux through a dual face lying on the boundary, which is then also a primal face, is the same as the primal flux, that is zero.
For the expression of these densities and fluxes, we refer to \cite{gas-11-sta, her-14-som, her-10-kin}.
Since the mass balance is not yet solved at the velocity prediction stage, they have to be built from the mass balance at the previous time step: hence the backward time shift for the densities in the time-derivative term.

\medskip
The term $(\gradi p)_{\edge,i}$ stands for the $i$-th component of the discrete pressure gradient at the face $\edge$.
This gradient operator is built as the transpose of the discrete operator for the divergence of the velocity, \ie\ in such a way that the following duality relation with respect to the L$^2$ inner product holds:
\[
 \sum_{K\in\mesh}|K|p_K\dive_K(\bfu) + \sum_{i=1}^d \sum_{\edge\in\edgesischeme} |D_\edge|u_{\edge,i}(\gradi p)_{\edge,i} = 0.
\]
This leads to the following expression:
\[
 \forall\edge=K|L\in\edgesint,\qquad (\gradi p)_{\edge,i} = \frac{|\edge|}{|D_\edge|}(p_L-p_K)\bfn_{K,\edge}\cdot\bfe^{(i)}.
\]
The scaling of the pressure gradient \eqref{scheme:gradp} is necessary for the solution to the scheme to satisfy a local discrete (finite volume) kinetic energy balance \cite[Lemma 4.1]{gra-17-unc}.
%
%

\medskip
\paragraph{\textbf{Sensible enthalpy equation}}
The equation is discretised in such a way that the present enthalpy formulation is strictly equivalent to the internal energy formulation of the energy balance equation used in \cite{gra-17-unc}.
Consequently, the term $-\big(u \cdot \gradi p \big)_K$ reads:
\[
-\big(u \cdot \gradi p \big)_K=\frac{1}{|K|}\sum_{\edge\in\edges(K)}|\edge|\,u_{K,\edge}\,(p_K-p_\edge),
\]
where $p_\edge$ is the upwind approximation of $p$ at the face $\edge$ with respect to $u_{K,\edge}$.
The reaction heat, $\displaystyle(\dot{\omega}_\theta)_K$, is written in the following way:
\[
(\dot{\omega}_\theta)_K=-\sum_{i=1}^{N_s} \Delta h_{f,i}^0 \, (\dot{\omega}_i)_K=\left(\nu_F\,W_F\,\Delta h_{f,F}^0+\nu_O\,W_O\,\Delta h_{f,O}^0-\nu_P\,W_P\,\Delta h_{f,P}^0\right)\,\dot{\omega}_K.
\]
The definition of $\dot\omega_K$ is given in Section \ref{subsec:chem}, and the definition of the corrective term $S_K^{n+1}$ is given in Section \ref{sec:cons} (see Equation \eqref{eq:def_SK} and Remark \ref{rmrk:consist} below).
%
%
\subsection{Chemistry step}  \label{subsec:chem}

For $0 \leq n <N$, the step $n+1$ for the solution of the transport of the characteristic function of the burnt zone and the chemical species mass balance equations reads:
\begin{subequations}
\label{scheme:chem}
\begin{align}
& \nonumber \quad \mbox{\textit{Computation of the burnt zone characteristic function} -- Solve for $G^{n+1}$:}  & & & & & & & & & & & & & & \\
\label{scheme:G} & \qquad \forall K\in\mesh,
\qquad \frac{1}{\delta t}(\rho_K^nG_K^{n+1}-\rho_K^{n-1}G_K^{n})+\dive_K(\rho^n G^{n+1} \bfu^n) + (\rho_u^n u_f^n\, |\gradi G|)_K = 0. \\[2ex]
& \nonumber \quad \mbox{\textit{Computation of the variable $z$} -- Solve for $z^{n+1}$:} \\
\label{scheme:z} & \qquad \forall K\in\mesh,
\qquad \frac{1}{\delta t}(\rho_K^nz_K^{n+1}-\rho_K^{n-1}z_K^{n})+\dive_K(\rho^n z^{n+1} \bfu^n) = 0. \\[2ex]
& \nonumber \quad \mbox{\textit{Neutral gas mass fraction computation} -- Solve for $y_N^{n+1}$:} \\
\label{scheme:yN} & \qquad \forall K\in\mesh,
\qquad \frac{1}{\delta t}\big[\rho_K^n(y_N)_K^{n+1}-\rho_K^{n-1}(y_N)_K^{n}\big]+\dive_K(\rho^n y_N^{n+1} \bfu^n) = 0. \\[2ex]
& \nonumber \quad \mbox{\textit{Fuel mass fraction computation} -- Solve for $y_F^{n+1}$:} \\ 
\label{scheme:yF} & \qquad \forall K\in\mesh,
\qquad \frac{1}{\delta t}\big[\rho_K^n(y_F)_K^{n+1}-\rho_K^{n-1}(y_F)_K^{n}\big]+\dive_K(\rho^n y_F^{n+1} \bfu^n) = -\frac 1 {\varepsilon} \nu_F W_F\,\dot{\omega}^{n+1}_K. \\[2ex]
& \nonumber \quad \mbox{\textit{Product mass fraction computation} -- Compute $y_P^{n+1}$ given by:} \\ 
\label{scheme:yP} & \qquad \forall K\in\mesh,
\qquad (y_P)_K^{n+1} = 1 - (y_F)_K^{n+1} - (y_O)_K^{n+1} - (y_N)_K^{n+1}.
\end{align}
\end{subequations}

\medskip
The initial value of the chemical variables is the mean value of the initial condition over the primal cells:
\[
\forall K\in\mesh,\ G_K^0 = \frac1{|K|} \int_K\! G_0(\bfx)\,{\rm d}\bfx, \quad z_K^0 = \frac1{|K|} \int_K\! z_0(\bfx)\,{\rm d}\bfx, \quad (y_i)_K^0 = \frac1{|K|} \int_K\! (y_i)_0(\bfx)\,{\rm d}\bfx,\ i=N,F,
\]
where the reduced variable $z$ is the linear combination of $y_F$ and $y_O$ given by Equation \eqref{eq:def_z}.
According to the developments of Section \ref{sec:gene}, the chemical reaction term reads $\dot{\omega}^{n+1}_K= \eta((y_F)_K^{n+1}, z_K^{n+1})\ (G_K^{n+1}-0.5)^-$ with
\[
\eta((y_F)_K^{n+1},z_K^{n+1})=
\left| \begin{array}{ll} \displaystyle
\frac 1 {\nu_F W_F}\, (y_F)_K^{n+1}
&
\mbox{if } z^{n+1} \leq 0,
\\[3ex] \displaystyle
\frac 1 {\nu_F W_F}\, (y_F)_K^{n+1} - z_K^{n+1}
&
\mbox{otherwise,}
\end{array} \right.
\]
and the chemical species mass fractions satisfy the following system, which is equivalent to \eqref{scheme:z}-\eqref{scheme:yP}:
\begin{equation} \label{eq:chem_sys_d}
\frac{1}{\delta t} (\rho_K^n (y_i)_K^{n+1}-\rho_K^{n-1} (y_i)_K^n) + \dive_K(\rho^n y_i^k \bfu^n) =
\frac 1 {\varepsilon} \zeta_i \nu_i W_i\, \dot{\omega}^{n+1}_K, \quad \mbox{for } i \in \mathcal{I} \mbox{ and } K \in \mesh.
\end{equation}

\medskip
At the continuous level, the last term of equation \eqref{scheme:G} may be written:
\[
\rho_u \, u_f\, |\gradi G| = \bfa \cdot \gradi G = \dive(G\, \bfa) - G\, \dive(\bfa),\quad \mbox{with }
\bfa = \rho_u \, u_f\, \frac{\gradi G}{|\gradi G|}.
\]
Using an upwind finite volume discretization of both divergence terms in this relation, we get:
\[
|K|\ (\rho_u^n\,u_f^n \ |\gradi G|)_K=\sum_{\edge \in \edges(K)} |\edge|\ (G_\edge^{n+1}-G_K^{n+1})\, \bfa_\edge^n \cdot \bfn_{K,\edge},
\]
where $G_\edge^{n+1}$ stands for the upwind approximation of $G^{n+1}$ on $\edge$ with respect to $\bfa^n \cdot \bfn_{K,\edge}$.
The flame velocity on $\edge$, $\bfa_\edge^n$, is evaluated as
\[
\bfa_\edge^n = (\rho_u \,u_f)_\edge^n \ \frac{(\gradi G)^n_\edge}{|(\gradi G)^n_\edge|},
\]
where $(\rho_u \,u_f)_\edge^n$ stands for an approximation of the product $\rho_u \,u_f$ on the face $\edge$ at $t^n$ (this product being often constant in applications), and the gradient of $G$ on $\edge=K|L$ is computed as:
\[
(\gradi G)_\edge=\frac{1}{|K\cup L|} \Bigl[
\sum_{\edge' \in \edges(K)} |\edge'| \ \hat G_{\edge'}\ \bfn_{K,\edge'}
+ \sum_{\edge' \in \edges(L)} |\edge'| \ \hat G_{\edge'}\ \bfn_{L,\edge'}
\Bigr],
\]
where $\hat G_{\edge'}$ is a second order approximation of $G$ at the barycenter of the face $\edge'$.

%
%
\section{Scheme conservativity}\label{sec:cons}

Let the discrete sensible internal energy be defined by $p_K^n= (\gamma -1)\, \rho_K^n (e_s)_K^n$ for $K \in \mesh$ and $0 \leq n \leq N$.
In view of the equation of state \eqref{scheme:eos}, this definition implies $\rho_K^n (h_s)_K^n=\rho_K^n (e_s)_K^n+p_K^n$, for $K \in \mesh$ and $0 \leq n \leq N$.
The following lemma states that the discrete solutions satisfy a local internal energy balance.

\begin{lmm}[Discrete internal energy balance]\ \\
A solution to \eqref{scheme:euler}-\eqref{scheme:chem} satisfies the following equality, for any $K\in\mesh$ and $0 \leq n < N$:
\begin{equation} \label{scheme:e}
\frac{1}{\delta t}\bigl[ (\rho e)_K^{n+1} - (\rho e)_K^n \bigr] + \widetilde\dive_K(\rho e \bfu)^{n+1} + p_K^{n+1}\dive_K(\bfu)^{n+1} = S_K^{n+1},
\end{equation}
where
\[ \begin{array}{l} \displaystyle
(\rho e)_K^{n+1} = \rho_K^{n+1}(e_s)_K^{n+1} + \rho_K^n\sum_{i \in \mathcal{I}} \Dhi (y_i)_K^{n+1},
\\[3ex] \displaystyle
\widetilde\dive_K(\rho e\bfu)^{n+1} = \dive_K\Bigl[(\rho e_s)^{n+1} \bfu^{n+1} + \rho^n \bigr[ \sum_{i\in\mathcal{I}} \Dhi y_i^{n+1} \bigl]\bfu^n\Bigr].
\end{array}
\]
\end{lmm}

\begin{proof}
We begin with deriving a local sensible internal energy balance, starting from the sensible enthalpy balance \eqref{scheme:hs} and mimicking the formal passage between these two equations at the continuous level given previously (\ie\ the passage from Equation \eqref{eq:hs} to Equation \eqref{eq:es}).
To this purpose, let us write \eqref{scheme:hs} as $T_1+T_2= T_3$ with
\[
\begin{array}{l} \displaystyle
T_1 = \frac{1}{\delta t}\,\bigl[\rho^{n+1}_K\,(h_s)^{n+1}_K-\rho^{n}_K\,(h_s)^n_K \bigr] -\frac{1}{\delta t}\,(p^{n+1}_K- p^n_K),
\\[2ex] \displaystyle
T_2 = \dive_K(\rho h_s \bfu)^{n+1} - \big(\bfu \cdot \gradi p\big)_K^{n+1},
\\[2ex] \displaystyle
T_3 = (\dot{\omega}_\theta)_K^{n+1} + S_K^{n+1}.
\end{array}
\]
Using $\rho_K^\ell (h_s)_K^\ell=\rho_K^\ell (e_s)_K^\ell+p_K^\ell$ for $\ell=n$ and $\ell=n+1$, we easily get
\[
T_1 = \frac{1}{\delta t}\,\bigl[\rho^{n+1}_K\,(e_s)^{n+1}_K-\rho^{n}_K\,(e_s)^n_K \bigr].
\]
The term $T_2$ reads:
\[
|K|\ T_2= \sum_{\edge \in \edges(K)} |\edge|\ \bigl[\rho_\edge^{n+1} (h_s)_\edge^{n+1} - p_\edge^{n+1} + p_K^{n+1} \bigr] u_{K,\edge}^{n+1}.
\]
If $u_{K,\edge}^{n+1} > 0$, by definition, $\rho_\edge^{n+1} (h_s)_\edge^{n+1}=\rho_K^{n+1} (h_s)_K^{n+1}$ and $p_\edge^{n+1}=p_K^{n+1}$; otherwise, thanks to the assumptions on the boundary conditions, $\edge$ is an internal face and, denoting by $L$ the adjacent cell to $K$ such that $\edge=K|L$, $\rho_\edge^{n+1} (h_s)_\edge^{n+1}=\rho_L^{n+1} (h_s)_L^{n+1}$ and $p_\edge^{n+1}=p_L^{n+1}$.
In both cases, denoting by $(e_s)_\edge^{n+1}$ the upwind choice for $(e_s)^{n+1}$ at the face $\edge$, we get
\[
\rho_\edge^{n+1} (h_s)_\edge^{n+1} - p_\edge^{n+1} = \rho_\edge^{n+1} (e_s)_\edge^{n+1},
\]
so, finally
\[
|K|\ T_2= \sum_{\edge \in \edges(K)} F_{K,\edge}^{n+1} (e_s)_\edge^{n+1}  + p_K^{n+1} \sum_{\edge \in \edges(K)}|\edge|\  u_{K,\edge}^{n+1}.
\]
We thus get the following sensible internal energy balance:
\begin{equation}\label{eq:es_long}
\frac{|K|}{\delta t}\,\bigl[\rho^{n+1}_K\,(e_s)^{n+1}_K-\rho^{n}_K\,(e_s)^n_K \bigr]
+\sum_{\edge \in \edges(K)} F_{K,\edge}^{n+1} (e_s)_\edge^{n+1}  + p_K^{n+1} \sum_{\edge \in \edges(K)}|\edge|\  u_{K,\edge}^{n+1}
= |K|\ \bigl[ (\dot{\omega}_\theta)_K^{n+1} + S_K^{n+1} \bigr],
\end{equation}
or, using the discrete differential operator formalism,
\begin{equation}\label{eq:esd}
\frac{1}{\delta t}\,\bigl[\rho^{n+1}_K\,(e_s)^{n+1}_K-\rho^{n}_K\,(e_s)^n_K \bigr]
+\dive_K(\rho e_s \bfu)^{n+1} + p_K^{n+1} \dive_K \bfu^{n+1}
= (\dot{\omega}_\theta)_K^{n+1} + S_K^{n+1}.
\end{equation}
We now derive from this relation a discrete (sensible and chemical) internal energy balance.
Multiplying the mass fraction balance equations by the corresponding formation enthalpy $(\Dhi)_{i\in\mathcal I}$ and suming over $i\in\mathcal I$ yields:
\[
\frac 1 \delta t \sum_{i\in\mathcal{I}} \Dhi \bigl[ \rho_K^n (y_i)_K^{n+1} - \rho_K^{n+1} (y_i)_K^n \bigr] +
\sum_{\edge \in \edges(K)} F_{K,\edge}^n \sum_{i\in\mathcal{I}} \Dhi\ (y_i)_\edge^{n+1}
= \sum_{i\in\mathcal{I}} \Dhi\ (\dot\omega_i)_K^{n+1} = (\dot\omega_\theta)_K^{n+1}.
\]
Adding this relation to \eqref{eq:es_long} yields the balance equation which we are looking for.
\end{proof}

\medskip
\begin{rmrk}[Positivity of the sensible internal energy]
Equation \eqref{eq:esd} implies that the sensible internal energy remains positive, provided that the right-hand side is non-negative, which is true if $\dot{\omega}_\theta \geq 0$, \ie\ if the chemical reaction is exothermic.
The proof of this property may be found in \cite[Lemma 4.3]{gra-17-unc}, and relies on two arguments: first, the convection operator may be recast as a discrete positivity-preserving transport operator thanks to the mass balance, and, second, the pressure $p_K^{n+1}$ vanishes when $e_K^{n+1}$, by the equation of state.
\end{rmrk}

\medskip 
The following local discrete kinetic energy balance holds on the dual mesh (see \cite[Lemma 4.1]{gra-17-unc} for a proof).

\begin{lmm}[Discrete kinetic energy balance on the dual mesh]\ \\
A solution to \eqref{scheme:euler}-\eqref{scheme:chem} satisfies the following equality, for $1\leq i \leq d$, $\edge \in \edgesischeme$ and $0 \leq n < N$:
\begin{equation} \label{eq:ec_face}
\frac{|D_\edge|}{\delta t}\bigl[ (e_k)_{\edge,i}^{n+1} - (e_k)_{\edge,i}^n \bigr] 
+ \sum_{\edged\in\edgesd(D_\edge)} F_{\edge,\edged}^n (e_k)_{\edged,i}^{n+1}
+ |D_\edge|\ (\gradi p)_{\edge,i}^{n+1}u_{\edge,i}^{n+1} = -R_{\edge,i}^{n+1},
\end{equation}
where
\[
\begin{array}{l} \displaystyle
(e_k)_{\edge,i}^{n+1} = \frac 1 2 \rho_\Ds^n (u_{\edge,i}^{n+1})^2 + \delta t^2  \frac{|D_\edge|}{2 \rho_\Ds^n} \big( (\gradi p)_{\edge,i}^{n+1} \big)^2,
\\[2ex] \displaystyle
(e_k)_{\edged,i}^{n+1} = \frac 1 2 \tilde u_{\edge,i}^{n+1} \tilde u_{\edge',i}^{n+1},
\\[2ex] \displaystyle
R_{\edge,i}^{n+1} = \frac{|D_\edge|\ \rho_\Ds^{n-1}}{2\delta t}(\tilde u_{\edge,i}^{n+1} - u_{\edge,i}^n)^2.
\end{array}
\]
\end{lmm}

\medskip
We now derive a kinetic energy balance equation on the primal cells from Relation \eqref{eq:ec_face}.
For the sake of clarity, we make a separate exposition for the Rannacher-Turek case and the MAC case.
The case of simplicial discretizations, with the degrees of freedom of the Crouzeix-Raviart element, is an easy extension of the Rannacher-Turek case.

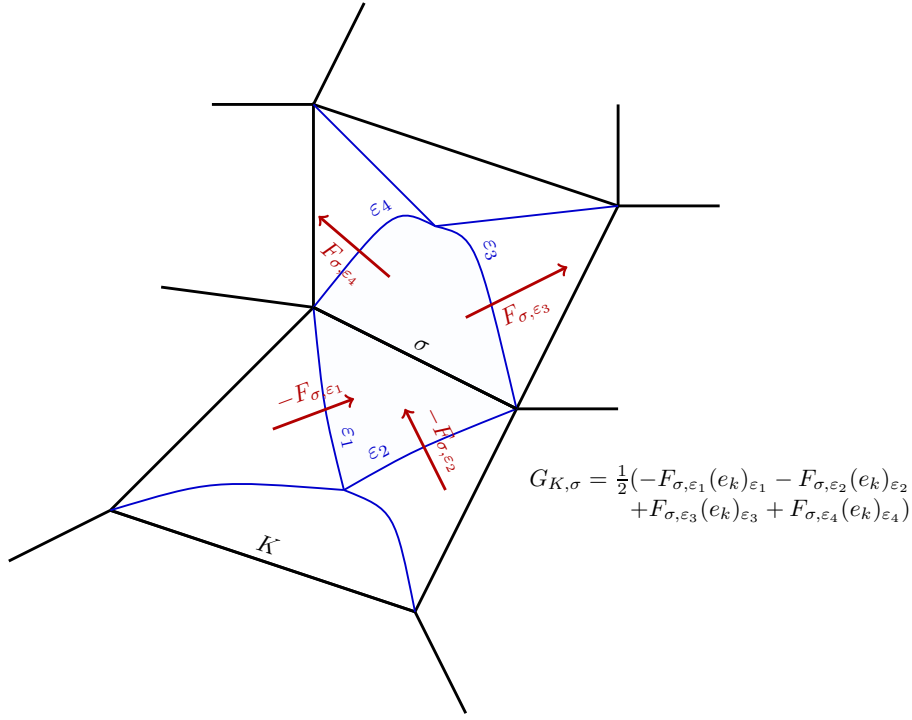
\begin{figure}[!t]
\begin{center}
\scalebox{0.9}{
\begin{tikzpicture}[scale=1.5]
\fill[color=bleuc!30!white,opacity=0.3] (3.,4.) .. controls (3.1,3.) .. (3.3,2.2) .. controls (4.,2.6) .. (5.,3.);
\fill[color=bleuc!30!white,opacity=0.3] (3.,4.) .. controls (3.8,5.) .. (4.2,4.8) .. controls (4.6,4.7) .. (5.,3.);
\draw[very thick, color=black] (1.,2.) -- (4.,1.) -- (5.,3.) -- (3.,4.) -- (1.,2.);
\draw[very thick, color=black] (1.,2.) -- (4.,1.) node[midway,sloped,above]{$K$};
\draw[very thick, color=black] (5.,3.) -- (6.,5.) -- (3.,6.) -- (3.,4.);
\draw[very thick, color=black] (1.,2.) -- (0.,1.5);
\draw[very thick, color=black] (4.,1.) -- (4.5,0.);
\draw[very thick, color=black] (3.,4.) -- (1.5,4.2);
\draw[very thick, color=black] (5.,3.) -- (6.,3);
\draw[very thick, color=black] (7.,5.) -- (6.,5.);
\draw[very thick, color=black] (6.,5.) -- (6.,6.);
\draw[very thick, color=black] (3.,6.) -- (3.5,7.);
\draw[very thick, color=black] (3.,6.) -- (2.,6.);
\draw[thick, color=bleuf] (3.,4.) .. controls (3.1,3.) .. (3.3,2.2) node[near end,sloped,above]{$\edged_1$};
\draw[very thick, color=red!70!black, <-] (3.4,3.1) -- (2.6,2.8) node[midway,sloped,above]{$\ -F_{\edge,\edged_1}$};
\draw[thick, color=bleuf] (3.3,2.2) .. controls (4.,2.6) .. (5.,3.) node[near start,sloped,above]{$\edged_2$};
\draw[very thick, color=red!70!black, <-] (3.9,3.) -- (4.3,2.2) node[midway,sloped,above]{$-F_{\edge,\edged_2}$};
\draw[thick, color=bleuf] (4.2,4.8) .. controls (4.6,4.7) .. (5.,3.) node[midway,sloped,above]{$\edged_3$};
\draw[very thick, color=red!70!black, <-] (5.5,4.4) -- (4.5,3.9) node[midway,sloped,below]{$F_{\edge,\edged_3}$};
\draw[thick, color=bleuf] (3.,4.) .. controls (3.8,5.) .. (4.2,4.8)  node[midway,sloped,above]{$\edged_4$};
\draw[very thick, color=red!70!black, <-] (3.05,4.9) -- (3.75,4.3) node[midway,sloped,below]{$\,F_{\edge,\edged_4}$};
\draw[thick, color=bleuf] (1.,2.) .. controls (2.,2.3) .. (3.3,2.2);
\draw[thick, color=bleuf] (4.,1.) .. controls (3.8,2.) .. (3.3,2.2);
\draw[thick, color=bleuf] (6.,5.) -- (4.2,4.8);
\draw[thick, color=bleuf] (3.,6.) -- (4.2,4.8);
\draw[very thick, color=black] (5.,3.) -- (3.,4.);
\draw[very thick, color=black] (5.,3.) -- (3.,4.) node[midway,sloped,above]{$\edge$};
\draw (7,2.3) node{$G_{K,\edge}=\frac 1 2 (-F_{\edge,\edged_1} (e_k)_{\edged_1}-F_{\edge,\edged_2} (e_k)_{\edged_2}$};
\draw (7.5,2.) node{$+F_{\edge,\edged_3} (e_k)_{\edged_3}+F_{\edge,\edged_4} (e_k)_{\edged_4})$};
\end{tikzpicture}
}
\end{center}
\caption{From fluxes at dual faces to fluxes at primal faces, for the Rannacher-Turek discretization.}
\label{fig:convRT}
\end{figure}

\medskip
\paragraph{\textbf{The Rannacher-Turek case}}
Since the dual meshes are the same for all the velocity components in this case, we may sum up Equation \eqref{eq:ec_face} over $i=1,\dots d$ to obtain, for $\edge \in \edges$ and $0 \leq n <N$:
\begin{equation} \label{eq:ec_faceRT}
\frac{|D_\edge|}{\delta t}\bigl[ (e_k)_\edge^{n+1} - (e_k)_\edge^n \bigr] 
+ \sum_{\edged\in\edgesd(D_\edge)} F_{\edge,\edged}^n (e_k)_\edged^{n+1}
+ |D_\edge|\ (\gradi p)_\edge^{n+1} \cdot \bfu_\edge^{n+1} = -R_\edge^{n+1},
\end{equation}
with
\[
(e_k)_\edge^\ell = \sum_{i=1}^d (e_k)_{\edge,i}^\ell,\mbox{ for } \ell=n \mbox{ or } \ell=n+1,\quad
(e_k)_\edged^{n+1} = \sum_{i=1}^d (e_k)_{\edged,i}^{n+1},\quad
R_\edge^{n+1} = \sum_{i=1}^d R_{\edge,i}^{n+1}.
\]
For $K\in\mesh$, let us define a kinetic energy associated to $K$ and the flux $G_{K,\edge}^{n+1}$ as follows (see Figure \ref{fig:convRT}):
\[
\begin{array}{l} \displaystyle
(e_k)_K^\ell = \frac 1 {|K|} \sum_{\edge \in \edges(K)} |D_\edge|\ (e_k)_\edge^\ell,\ \ell=n \mbox{ or } \ell=n+1,
\\[2ex] \displaystyle
G_{K,\edge}^{n+1} = -\frac 1 2 \sum_{\edged \in \edges(D_\edge), \edged \subset K} F_{\edge,\edged}^n\ (e_k)_\edged^{n+1}
+ \frac 1 2\sum_{\edged \in \edges(D_\edge), \edged \not \subset K} F_{\edge,\edged}^n\ (e_k)_\edged^{n+1}.
\end{array}
\]
We easily check that the fluxes $G_{K,\edge}^{n+1}$ are conservative, in the sense that, for $\edge=K|L$, $G_{K,\edge}^{n+1}=-G_{L,\edge}^{n+1}$.
Let us now divide Equation \eqref{eq:ec_faceRT} by $2$ and sum over the faces of $K$.
A reordering of the summations, using the conservativity of the mass fluxes through the dual edges and the expression of the discrete pressure gradient, yields:
\begin{multline} \label{eq:ec_meshRT}
\frac{|K|}{\delta t}\bigl[ (e_k)_K^{n+1} - (e_k)_K^n \bigr] 
+ \sum_{\edge \in \edges(K)} G_{K,\edge}^{n+1}
+ \sum_{\edge=K|L} |\edge|\ (p_L^{n+1}-p_K^{n+1})\ u_{K,\edge}^{n+1} 
= -R_K^{n+1},
\\
\mbox{with } R_K^{n+1}=\frac 1 2 \sum_{\edge \in \edges(K)} R_\edge^{n+1}. \hspace{10ex}
\end{multline}

\begin{figure}[!t]
\begin{center}
\scalebox{0.9}{
\begin{tikzpicture}[scale=1.3]
\fill[color=bleuc!30!white,opacity=0.3] (2.75,5.) -- (5.75,5.) -- (5.75,7.) -- (2.75,7.) -- (2.75,5.) ;
\draw[very thick, color=black] (1.,5.) -- (8.,5.) node[midway,sloped,above]{$\hspace{30ex} K$} ;
\draw[very thick, color=black] (1.,7.) -- (8.,7.) ;
\draw[very thick, color=black] (1.5,4.5) -- (1.5,7.5) ;
\draw[very thick, color=black] (4.,4.5) -- (4.,7.5)  node[midway,sloped,above]{$\edge$} ;
\draw[very thick, color=black] (7.5,4.5) -- (7.5,7.5) ;
\draw[thick, color=bleuf] (2.75,5.) -- (2.75,7.) node[near start,sloped,above]{$\edged_1$};
\draw[very thick, color=red!70!black, <-] (2.25,6.) -- (3.25,6.) node[near end,sloped,above]{$F_{\edge,\edged_1}$};
\draw[thick, color=bleuf] (5.75,5.) -- (5.75,7.) node[near start,sloped,above]{$\edged_2$};
\draw[very thick, color=red!70!black, ->] (5.25,6.) -- (6.25,6.) node[near start,sloped,above]{$-F_{\edge,\edged_2}$};
\draw (10.,5.9) node{$G_{K,\edge,1}=\frac 1 2 \bigl[ F_{\edge,\edged_1} (e_k)_{\edged_1,1}-F_{\edge,\edged_2} (e_k)_{\edged_2,1} \bigr]$};
\fill[color=bleuc!30!white,opacity=0.3] (1.2,1.) -- (4.24,1.) -- (4.24,3.) -- (1.2,3.) -- (1.2,1.) ;
\fill[color=vertc!30!white,opacity=0.3] (4.26,1.) -- (7,1.) -- (7,3.) -- (4.26,3.) -- (4.26,1.) ;
\draw[very thick, color=black] (1.,1.) -- (8.,1.) node[midway,sloped,above]{$K$} ;
\draw[very thick, color=black] (1.,3.) -- (8.,3.) node[midway,sloped,above]{$\tau$} ;
\draw[very thick, color=black] (2.5,0.5) -- (2.5,3.5)  node[midway,sloped,below]{$\edge$} ;
\draw[very thick, color=black] (6,0.5) -- (6,3.5) node[midway,sloped,below]{$\edge'$} ;
\draw[thick, color=bleuf] ((1.2,2.95) -- (4.2,2.95) node[near start,sloped,above]{$\edged$};
\draw[very thick, color=red!70!black, ->] (2.7,2.5) -- (2.7,3.5) node[near end,sloped,below]{$F_{\edge,\edged}$};
\draw[thick, color=bleuf] (4.3,2.95) -- (7,2.95) node[near start,sloped,above]{$\edged'$};
\draw[very thick, color=red!70!black, ->] (5.8,2.5) -- (5.8,3.5) node[near end,sloped,above]{$\ F_{\edge',\edged'}$};
\draw (10.,1.9) node{$G_{K,\tau,1}=\frac 1 2 \bigl[ F_{\edge,\edged} (e_k)_{\edged,1}+F_{\edge',\edged'} (e_k)_{\edged',1} \bigr]$};
\end{tikzpicture}
}
\end{center}
\caption{From fluxes at dual faces to fluxes at primal faces, for the MAC discretization -- First component of the velocity.}
\label{fig:convMAC}
\end{figure}
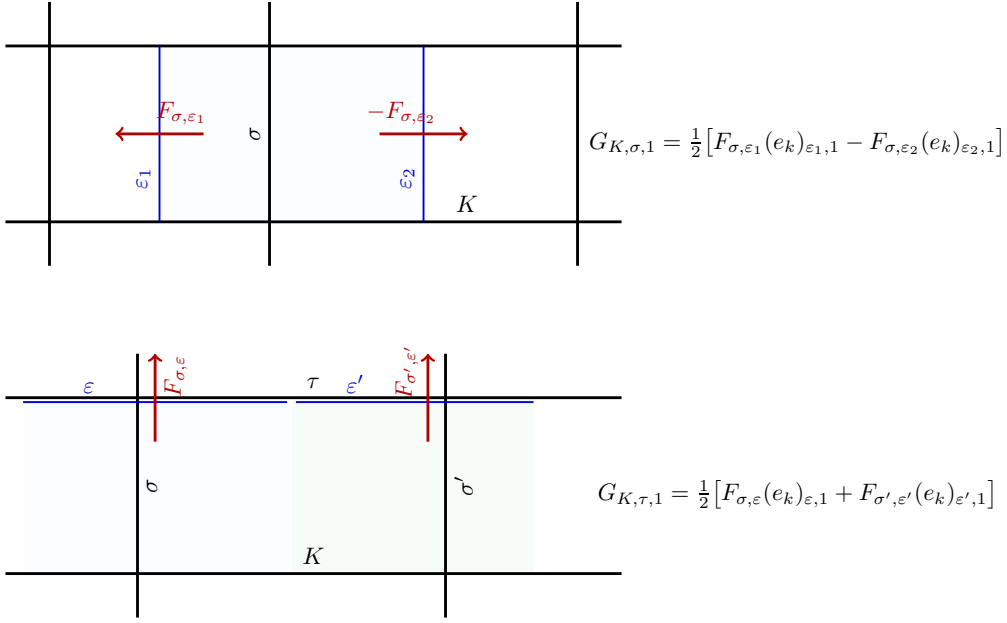

\medskip
\paragraph{\textbf{The MAC case}}
Let $1 \leq i \leq d$, let $K \in \mesh$, let us denote by $\edge$ and $\edge'$ the two faces of $\edges\ei(K)$, and let us define:
\[
(e_k)_{K,i}^\ell = \frac 1 {|K|} \Bigl[ |D_\edge|\ (e_k)_{\edge,i}^\ell+ |D_\edge|\ (e_k)_{\edge,i}^\ell \Bigr],\mbox{ for } \ell=n \mbox{ or } \ell=n+1.
\]
Let $\tau=\edge$ or $\tau=\edge'$, let $\edged$ and $\edged'$ be the two faces of $D_\tau$ perpendicular to $\bfe\ei$, and let $\edge'$ be included in $K$ (see Figure \ref{fig:convMAC}, top).
Then we define
\[
G_{K,\tau,i}^{n+1} = \frac 1 2 \bigl[ F_{\tau,\edged} (e_k)_{\edged,i}^{n+1} - F_{\tau,\edged'} (e_k)_{\edged',i}^{n+1} \bigr].
\]
For $\tau \in \edges(K) \setminus \{\edge,\edge'\}$, let $\edged$ and $\edged'$ be such that $\tau \subset (\bar \edged \cup \bar \edged')$ with $\edged$ a face of $D_\edge$ and $\edged'$ a face of $D_\edged'$ (see Figure \ref{fig:convMAC}, bottom).
Then we define
\[
G_{K,\tau,i}^{n+1} = F_{\sigma,\edged} (e_k)_{\edged,i}^{n+1} - F_{\sigma',\edged'} (e_k)_{\edged',i}^{n+1}.
\]
Let us now divide by $2$ Equation \eqref{eq:ec_face} written for $\edge$ and for $\edge'$ and sum.
We get:
\begin{equation} \label{eq:ec_meshMACi}
\frac{|K|}{\delta t}\bigl[ (e_k)_{K,i}^{n+1} - (e_k)_{K,i}^n \bigr] 
+ \sum_{\edge \in \edges(K)} G_{K,\edge,i}^{n+1}
+ \sum_{\edge \in \edges\ei(K),\ \edge=K|L} |\edge|\ (p_L^{n+1}-p_K^{n+1}) u_{K,\edge}^{n+1} 
= -\frac 1 2 \bigl( R_{\edge,i}^{n+1} + R_{\edge,i}^{n+1} \bigr).
\end{equation}
Let now
\[
(e_k)_K^\ell = \sum_{i=1}^d (e_k)_{K,i}^\ell,\mbox{ for } \ell=n \mbox{ or } \ell=n+1,\quad \mbox{and} \quad
G_{K,\sigma}^{n+1} = \sum_{i=1}^d G_{K,\sigma,i}^{n+1}, \mbox{ for } \edge \in \edges(K).
\]
Since only one equation is written for a given face $\edge$ of the mesh (for the velocity component $i$ with $i$ such that the normal vectorto $\edge$ is parallel to $\bfe\ei$), we may define in the MAC case $R_\edge^{n+1}=R_{\edge,i}^{n+1}$.
Summing Equation \ref{eq:ec_meshMACi} over the space dimension, we finally get
\begin{multline} \label{eq:ec_meshMAC}
\frac{|K|}{\delta t}\bigl[ (e_k)_K^{n+1} - (e_k)_K^n \bigr] 
+ \sum_{\edge \in \edges(K)} G_{K,\edge}^{n+1}
+ \sum_{\edge=K|L} |\edge|\ (p_L^{n+1}-p_K^{n+1})\ u_{K,\edge}^{n+1} 
= -R_K^{n+1},
\\
\mbox{with } R_K^{n+1}=\frac 1 2 \sum_{\edge \in \edges(K)} R_\edge^{n+1}, \hspace{10ex}
\end{multline}
which is formally the same equation as Relation \eqref{eq:ec_meshRT} (although with a different definition of all the terms in the equation except the pressure gradient).

\medskip
\begin{rmrk}[On the definition of the cell kinetic energy]
Note that, both in the Rannanacher-Turek and the MAC case, the cell kinetic energy is not a convex combination of the face kinetic energies, since, on a non-uniform mesh, the equality $|K| = \frac 1 2 \sum_{\edge \in \edges(K)} |D_\edge|$ is generally false.
Consequently, the cell kinetic energy may for instance oscillate from cell to cell while the face kinetic energy does not.
Anyway, the discrete time derivative of the cell kinetic energy is consistent in the Lax-Wendroff sense.
\end{rmrk}

\medskip
Equations \eqref{eq:ec_meshRT} and \eqref{eq:ec_meshMAC} suggest a choice for the term $S_K^{n+1}$, the purpose of which is to compensate the numerical dissipation terms appearing in the kinetic energy balance:
\begin{equation}\label{eq:def_SK}
S_K^{n+1} = R_K^{n+1}, \mbox{ for } K \in \mesh \mbox{ and } 0 \leq n < N.
\end{equation}
This expression yields a conservative scheme, in the sense that the discrete solutions satisfy a discrete total energy balance without any remainder term (see Equation \eqref{eq:etot} below); as a consequence, the scheme can be proven to be consistent in the Lax-Wendroff sense.
However, different definitions are possible (and this latitude may be useful in explicit variants of the scheme, to ensure the positivity of $S_K^{n+1}$, see Remark \ref{rmrk:consist} below.

\medskip
We are now in position to state a total energy balance for the scheme.

\begin{thrm}[Discrete total energy and stability of the scheme]\ \\
A solution to \eqref{scheme:euler}-\eqref{scheme:chem} satisfies the following equality, for any $K\in\mesh$ and $0 \leq n < N$:
\begin{equation} \label{eq:etotd}
 \frac{1}{\delta t}\bigl[(\rho E)_K^{n+1} - (\rho E)_K^n\bigr] + \widetilde\dive_K((\rho E +p) \bfu)^{n+1} = 0,
\end{equation}
where
\[ \begin{array}{l} \displaystyle
(\rho E)_K^\ell = (e_k)_K^\ell + \rho_K^\ell (e_s)_K^\ell + \rho_K^{l-1}\sum_{i \in \mathcal{I}} \Dhi (y_i)_K^\ell,\mbox{ for } \ell=n \mbox{ and } \ell=n+1,
\\[3ex] \displaystyle
\widetilde\dive_K((\rho E+p)\,\bfu)^{n+1} = \dive_K\Bigl[(\rho e_s)^{n+1} \bfu^{n+1} + \rho^n \bigr[ \sum_{i\in\mathcal{I}} \Dhi y_i^{n+1} \bigl]\bfu^n\Bigr]
+ \sum_{\edge=K|L} |\edge|\ (p^{n+1}_K+p^{n+1}_L)\, u_{K,\edge}^{n+1}.
\end{array}
\]
Let us suppose that $e_s^0,\ \rho^0$ and $\rho^{-1}$ are positive. 
Then, a solution to \eqref{scheme:euler}-\eqref{scheme:chem} satisfies $\rho^{n+1}>0$, $e^{n+1}>0$ and the following stability result:
\begin{equation*}
E^n = E^0,
\end{equation*}
where, for $0\leq n \leq N$,
\[
E^n = \sum_{K\in\mesh} |K|(\rho e)_K^n + \frac12 \sum_{i=1}^d \sum_{\edge \in \edgesischeme} |D_\edge| (u_{\edge,i}^n)^2 
+ \delta t^2\sum_{\edge \in \edgesint} \frac{|D_\edge|}{\rho_\Ds^{n-1}}|(\gradi p)_\edge^n|^2.
\]
\end{thrm}

\begin{proof}
The discrete total energy balance equation \eqref{eq:etotd} is obtained by summing the internal energy balance \eqref{scheme:e} and the kinetic energy balance, \ie\ Equation \eqref{eq:ec_meshRT} in the Rannacher-Turek case and Equation \eqref{eq:ec_meshMAC} for the MAC scheme, and remarking that the numerical dissipation terms in the kinetic energy balance $R_K^{n+1}$ exactly compensate with the corrective terms $S_K^{n+1}$ in the internal energy balance.
Then the stability result is obtained by summation over the time steps.
\end{proof}

\medskip
\begin{rmrk}[Consistency of the scheme] \label{rmrk:consist}
The consistency in the Lax-Wendroff sense follows from the conservativity of the scheme (for all balance equations) so, in particular, from the fact that the discrete solutions satisfy the discrete total energy balance \eqref{eq:etotd}, thanks to standard (but technical) arguments.\\
Note however that the consistency of the scheme does not require a strict conservativity, and in particular, variants for the choice \eqref{eq:def_SK} of the compensation term in the sensible enthalpy balance are possible; indeed, what is really needed is only that the difference between the dissipation in the kinetic energy balance and its compensation tend to zero in a distributional sense.
In practice, this allows a different redistribution of the face residuals to the neighbour primal cells, and this can help to preserve the non-negativity of the compensation term for explicit versions of the scheme.
\end{rmrk}
%
%
\section{Higher order convection schemes}

\subsection{The MUSCL scheme} \label{an:MUSCL}

The MUSCL discretization of the convection operators of the chemical species balance and $G$-equation is inspired by the discretisation proposed in \cite{pia-13-for}.
Let us use the following system of equations,
\[
\begin{array}{l}
\displaystyle \partial_t\rho + \dive(\rho\bfu) = 0 \\[1ex]
\displaystyle \partial_t(\rho\bfu) + \dive(\rho\bfu y) = 0,
\end{array}
\]
in order to explain here the MUSCL discretization of the convection operator in the transport equation of $y$.

\medskip
The discretization of the above system reads:
\[
\begin{array}{ll}
\forall K\in\mesh, \quad & \displaystyle \frac{\rho_K^{n+1}-\rho_K^n}{\delta t} + \frac1{|K|} \sum_{\edge\in\edges(K)} F_{K,\edge}^{n+1} = 0,
\\[4ex] & \displaystyle
\frac{\rho_K^{n+1} y_K^{n+1} - \rho_K^n y_K^n}{\delta t} + \frac1{|K|} \sum_{\edge\in\edges(K)} F_{K,\edge}^{n+1} y_\edge^n = 0.
\end{array}
\]
For any $\edge\in\edges$, the procedure consists in three steps:
\begin{list}{-}{\itemsep=0.5ex \topsep=0.5ex \leftmargin=1.cm \labelwidth=0.3cm \labelsep=0.5cm \itemindent=0.cm}
\item calculate a tentative value for $y_\edge$ as a linear interpolate of nearby values,
\item calculate an interval for $y_\edge$ which guarantees the stability of the scheme,
\item project the tentative value $y_\edge$ to this stability interval.
\end{list}

\medskip
For the tentative value of $y_\edge$, let us choose some real coefficients $(\alpha_K^\edge)_{K\in\mesh}$ such that
\[
\bfx_\edge = \sum_{K\in\mesh} \alpha_K^\edge \bfx_K, \qquad\qquad \sum_{K\in\mesh} \alpha_K^\edge = 1.
\]
The coefficients used in this interpolation are chosen in such a way that as few as possible cells, to be picked up in the closest cells to $\edge$, take part.
For example, for $\edge=K|L$ and if $\bfx_K,\ \bfx_\edge,\ \bfx_L$ are aligned, only two non-zero coefficients exist in the family $(\alpha_K^\edge)_{K\in\mesh}$, namely $\alpha_K^\edge$ and $\alpha_K^\edge$.
Then, these coefficients are used to calculate the tentative value of $y_\edge$ by
\[
y_\edge = \sum_{K\in\mesh} \alpha_K^\edge y_K.
\]

\medskip
The construction of the stability interval must be such that the following property holds:
\begin{equation}
\label{eq:MUSCL-hyp}
\begin{array}{c}
\forall K\in\mesh,\ \forall\edge\in\edges(K)\cap\edgesint,\ \exists\beta_K^\edge\in[0,1] \text{ and } M_K^\edge\in\mesh \text{ such that} \\[2ex]
y_\edge - y_K = 
\left|
\begin{array}{l}
\beta_K^\edge(y_K - y_{M_K^\edge}), \text{ if } F_{K,\edge} \geq 0,\\[2ex]
\beta_K^\edge(y_{M_K^\edge} - y_K), \text{ otherwise.}
\end{array}
\right.
\end{array}
\end{equation}
Indeed, under this latter hypothesis and a CFL condition, the scheme preserves the initial bounds of $y$.

\medskip
\begin{rmrk}
Note that in this work the presence of Neumann homogeneous boundary conditions limits our study to the internal faces, but what follows may be naturally generalized for non-homogeneous Dirichlet or/and Neumann boundary conditions.
\end{rmrk}

\begin{dfntn}
The so-called CFL number reads for any $0 \leq n \leq N$:
\[
 {\rm CFL}^n = \max_{K\in\mesh}\Big\{ \frac{\delta t}{\rho_K^{n+1}\ |K|} \sum_{\edge \in \edges(K)} \big|F_{K,\edge}^{n+1}\big| \Big\}.
\]
\end{dfntn}

\begin{lmm}
Let us suppose that ${\rm CFL}^{n+1}\leq1$.
For $K\in\mesh$, let us note by $\mathcal V(K)$ the union of the set of cells $M_K^\edge,\ \edge\in\edges(K)$ such that \eqref{eq:MUSCL-hyp} holds.
Then $\forall K\in\mesh$, the value of $y_K^{n+1}$ is a convex combination of $\{y_K^n, (y_M^n)_{M\in\mathcal V(K)} \}$.
\end{lmm}

\begin{proof}
The discrete mass balance equation yields:
\[
\rho_K^n = \rho_K^{n+1} + \frac{\delta t}{|K|} \sum_{\edge \in \edges(K)} F_{K,\edge}^{n+1}.
\]
Replacing this expression of $\rho_K^n$ in the discrete balance equation of $y$ and using the relations provided by \eqref{eq:MUSCL-hyp}, we obtain:
\[
\begin{array}{ll}
\rho_K^{n+1} y_K^{n+1} 
& \displaystyle \hspace{-2ex}
= \rho_K^n y_K^n - \frac{\delta t}{|K|} \sum_{\edge \in \edges(K)} F_{K,\edge}^{n+1} y_\edge^n 
\\[2ex] & \displaystyle \hspace{-2ex}
= \rho_K^{n+1} y_K^n - \frac{\delta t}{|K|} \sum_{\edge \in \edges(K)} F_{K,\edge}^{n+1} (y_\edge^n - y_K^n)
\\[2ex] & \displaystyle \hspace{-2ex}
= \rho_K^{n+1} y_K^n - \frac{\delta t}{|K|} \sum_{\edge \in \edges(K)} \big(F_{K,\edge}^{n+1}\big)^+ (y_\edge^n - y_K^n) 
+ \frac{\delta t}{|K|} \sum_{\edge \in \edges(K)} \big(F_{K,\edge}^{n+1}\big)^- (y_\edge^n - y_K^n)
\\[2ex] & \displaystyle \hspace{-2ex}
= \rho_K^{n+1} y_K^n - \frac{\delta t}{|K|} \sum_{\edge \in \edges(K)} \big(F_{K,\edge}^{n+1}\big)^+ \beta_K^\edge (y_K^n - y_{M_K^\edge}^n)
+ \frac{\delta t}{|K|} \sum_{\edge \in \edges(K)} \big(F_{K,\edge}^{n+1}\big)^- (y_{M_K^\edge}^n - y_K^n).
\end{array}
\]
This relation yields
\[
y_K^{n+1} = y_K^n \Big( 1 - \frac{\delta t}{\rho_K^{n+1}\ |K|} \sum_{\edge \in \edges(K)} \beta_K^\edge \big|F_{K,\edge}^{n+1}\big| \Big)
+ \frac{\delta t}{|K|} \sum_{\edge \in \edges(K)} y_{M_K^\edge}^n \beta_K^\edge \big|F_{K,\edge}^{n+1}\big|,
\]
which concludes the proof under the hypothesis that ${\rm CFL}\leq1$.
\end{proof}

\medskip
In practice, in order to construct the stability interval, we use a stronger version of \eqref{eq:MUSCL-hyp}, which allows us to be more precise in the choice of the control volumes $M_K^\edge,\ \forall K\in\mesh \text{ and } \forall\edge\in\edges(K)$.
Let $\edge\in\edges$, let us denote by $V^-$ and $V^+$ the upstream and downstream cell separated by $\edge$, and by $\mathcal V_\edge(V^-)$ and $\mathcal V_\edge(V^+)$ two sets of neighbouring cells of $V^-$ and $V^+$ respectively, and let us suppose:
\[
\begin{array}{l}
\text{(H1)} -\quad \text{there exists } M\in\mathcal V_\edge(V^+) \text{ such that } u_\edge^n\in |[u_M^n, u_M^n + \dfrac{\zeta^+}{2}(u_{V^+}^n - u_M^n) ]|,\\[2ex]
\text{(H2)} -\quad \text{there exists } M\in\mathcal V_\edge(V^-) \text{ such that } u_\edge^n\in |[u_{V^-}^n, u_{V^-}^n + \dfrac{\zeta^-}{2}(u_{V^-}^n - u_M^n) ]|,
\end{array}
\]
where for $a,\,b\in\mathbb R$, we denote by $|[a,b]|$ the interval $\{\alpha a + (1-\alpha)b,\ \alpha\in[0,1] \}$, and $\zeta^+$ and $\zeta^-$ are two numerical parameters lying in the interval $[0,2]$.

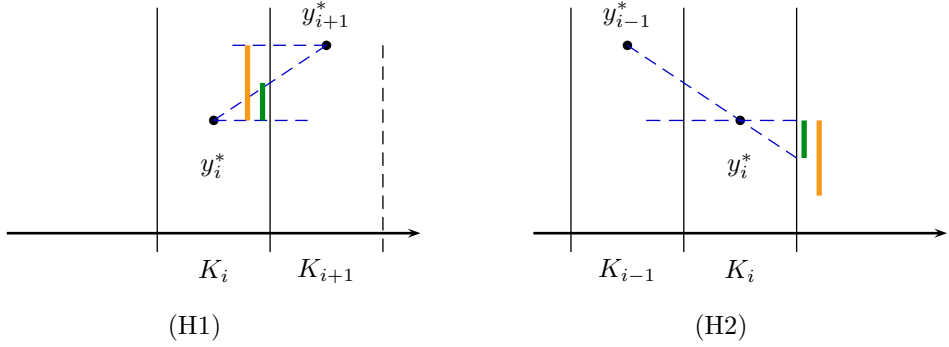
\begin{figure}[htb]
\psset{unit=1cm}
\begin{center}\begin{pspicture}(0,-0.5)(13,4)
\rput[bl](7,0){
%
   \psline[linecolor=black, linewidth=1pt]{->}(0.5,0.5)(6,0.5)
%
   \psline[linecolor=black, linewidth=0.5pt]{-}(1,0.25)(1,3.5)
   \psline[linecolor=black, linewidth=0.5pt]{-}(2.5,0.25)(2.5,3.5)
   \psline[linecolor=black, linewidth=0.5pt]{-}(4,0.25)(4,3.5)
   \rput(1.75,0){$K_{i-1}$}
   \rput(3.25,0){$K_i$}
%
   \psdots*(1.75,3) \rput(1.75,3.4){$y^\ast_{i-1}$}
   \psdots*(3.25,2.)\rput(3.25,1.4){$y^\ast_i$}
%
   \psline[linecolor=bleuf, linewidth=0.5pt, linestyle=dashed]{-}(1.75,3)(4,1.5)
   \psline[linecolor=bleuf, linewidth=0.5pt, linestyle=dashed]{-}(2,2)(4,2)
   \psline[linecolor=vertf, linewidth=2pt]{-}(4.1,1.5)(4.1,2)
   \psline[linecolor=or, linewidth=2pt]{-}(4.3,1)(4.3,2)
   \rput(3,-0.75){(H2)}
}
\rput[bl](0,0){
%
   \psline[linecolor=black, linewidth=1pt]{->}(0.5,0.5)(6,0.5)
%
   \psline[linecolor=black, linewidth=0.5pt]{-}(2.5,0.25)(2.5,3.5)
   \psline[linecolor=black, linewidth=0.5pt]{-}(4,0.25)(4,3.5)
   \psline[linecolor=black, linewidth=0.5pt, linestyle=dashed]{-}(5.5,0.25)(5.5,3.)
   \rput(3.25,0){$K_i$}
   \rput(4.75,0){$K_{i+1}$}
%
   \psdots*(3.25,2.) \rput(3.25,1.4){$y^\ast_i$}
   \psdots*(4.75,3) \rput(4.75,3.4){$y^\ast_{i+1}$}
%
   \psline[linecolor=bleuf, linewidth=0.5pt, linestyle=dashed]{-}(4.75,3)(3.25,2)
   \psline[linecolor=bleuf, linewidth=0.5pt, linestyle=dashed]{-}(3.25,2)(4.5,2)
   \psline[linecolor=bleuf, linewidth=0.5pt, linestyle=dashed]{-}(3.5,3)(4.75,3)
   \psline[linecolor=vertf, linewidth=2pt]{-}(3.9,2)(3.9,2.5)
   \psline[linecolor=or, linewidth=2pt]{-}(3.7,2)(3.7,3)
   \rput(3,-0.75){(H1)}
}
\end{pspicture}\end{center}
\caption{
Conditions (H1) and (H2) in 1D.
}
\label{fig:muscl1D}
\end{figure}

\begin{rmrk}[1D case]
Let us take the example of an interface $\edge$ separating $K_i$ and $K_{i+1}$ in a 1D case (see Figure \ref{fig:muscl1D} for the notations), with a uniform meshing and a positive advection velocity, so that $V^-=K_i$ and $V^+=K_{i+1}$.
In 1D, a natural choice is $\mathcal V_\edge(K_i)=\{K_{i-1}\}$ and $\mathcal V_\edge(K_{i+1})=\{K_i \}$.
On Figure \ref{fig:muscl1D}, we sketch: on the left, the admissible interval given by (H1) with $\zeta^+=1$ (green) and $\zeta^+=2$ (orange);
on the right, the admissible interval given by (H2) with $\zeta^-=1$ (green) and $\zeta^-=2$ (orange).
The parameters $\zeta^-$ and $\zeta^+$ may be seen as limiting the admissible slope between $(\bfx_i,y^n_i)$ and $(\bfx_\edge, y^n_\edge)$ (with $\bfx_i$ the abscissa of the mass centre of $K_i$ and $\bfx_\edge$ the abscissa of $\edge$), with respect to a left and right slope, respectively.
For $\zeta^-=\zeta^+=1$, one recognizes the usual minmod limiter (\eg\ \cite[Chapter III]{god-96-num}).
Note that, since, on the example depicted on Figure \ref{fig:muscl1D}, the discrete function $y^n$ has an extremum in $K_i$, the combination of the conditions (H1) and (H2) imposes that, as usual, the only admissible value for $y^n_\edge$ is the upwind one.
\label{rm:MUSCL-1D}
\end{rmrk}

\medskip
(H1)-(H2) and \eqref{eq:MUSCL-hyp} are linked in the following way: let $K\in\mesh$ and $\edge\in\edges(K)$. 
If $F_{K,\edge}^n\leq0,\ \ie\ K$ is the downstream cell for $\edge$, denoted above by $V^+$, since $\zeta^+\in[0,2]$, condition (H1) yields that there exists $M\in\mesh$ such that $u_\edge^n\in|[u_K^n,u_M^n]|$, which is \eqref{eq:MUSCL-hyp}.
Otherwise, \ie\ if $F_{K,\edge}^n\geq0$ and $K$ is the upstream cell for $\edge$, denoted above by $V^-$, condition (H2) yields that there exists $M\in\mesh$ such that $y_\edge^n\in|[y_K^n,2y_K^n-y_M^n]|$, so $y_\edge^n-y_K^n\in|[0,y_K^n-y_M^n]|$, which is once again \eqref{eq:MUSCL-hyp}.

\begin{rmrk}
For $\edge\in\edges$, if $V^-\in\mathcal V_\edge(V^+)$, the upstream choice $y_\edge^n=y_{V^-}^n$ always satisfies the conditions (H1)-(H2), and is the only one to satisfy them if we choose $\zeta^-=\zeta^+=0$.
\end{rmrk}

\begin{figure}[tb]
\psset{unit=1cm}
\begin{center}\scalebox{0.9}{\begin{pspicture}(0,-1)(13,3.5)
\rput[bl](0,0){
%
   \pspolygon*[linecolor=bleuc](2,1)(3.5,1.1)(3.4,2.4)(2.2,2)
   \rput[bl](2.2,1.1){{$V^-$}}
%
   \pspolygon*[linecolor=vertc](3.5,1.1)(5,1.3)(5.2,2.5)(3.4,2.4)
   \rput[bl](4.5,1.4){{$V^+$}}
%
   \pspolygon*[linecolor=or](0.4,1.7)(2.1,1.5)(2.2,2)(0.4,2.2)
   \pspolygon*[linecolor=or](0.5,1.1)(2,1)(2.1,1.5)(0.4,1.7)
   \pspolygon*[linecolor=or](2.2,2)(3.4,2.4)(3,3)(2,2.9)
%
   \psline[linecolor=black, linewidth=1pt]{-}(2,1)(3.5,1.1)(3.4,2.4)(2.2,2)(2,1)
   \psline[linecolor=black, linewidth=1pt]{-}(3.5,1.1)(5,1.3)(5.2,2.5)(3.4,2.4)(3.5,1.1)
   \psline[linecolor=black, linewidth=1pt]{-}(0.4,1.7)(2.1,1.5)(2.2,2)(0.4,2.2)(0.4,1.7)
   \psline[linecolor=black, linewidth=1pt]{-}(0.5,1.1)(2,1)(2.1,1.5)(0.4,1.7)(0.5,1.1)
   \psline[linecolor=black, linewidth=1pt]{-}(2,0)(3.5,0.1)(3.5,1.1)(2,1)(2,0)
   \psline[linecolor=black, linewidth=1pt]{-}(2.2,2)(3.4,2.4)(3,3)(2,2.9)(2.2,2)
   \psline[linecolor=black, linewidth=1pt]{-}(3.4,2.4)(5.2,2.5)(5.1,3.1)(3,3)(3.4,2.4)
   \psline[linecolor=black, linewidth=1pt]{-}(0.4,2.2)(2.2,2)(2,2.9)(0,3.3)(0.4,2.2)
   \psline[linecolor=black, linewidth=1pt]{-}(0.5,1.1)(2,0)
   \psline[linecolor=black, linewidth=1pt]{-}(3.5,0.1)(5,1.3)
%
   \psline[linecolor=black, linewidth=1pt]{->}(0,0.6)(0.8,0.2) \rput[bl](0,0){${\bf F}$}
   \rput[bl](2.7,-1){(a)}
}
\rput[bl](7,0){
%
   \pspolygon*[linecolor=bleuc](2,1)(3.5,1.1)(3.4,2.4)(2.2,2)
   \rput[bl](2.2,1.1){{$V^-$}}
%
   \pspolygon*[linecolor=vertc](3.5,1.1)(5,1.3)(5.2,2.5)(3.4,2.4)
   \rput[bl](4.5,1.4){{$V^+$}}
%
   \pspolygon*[linecolor=or](0.4,1.7)(2.1,1.5)(2.2,2)(0.4,2.2)
   \pspolygon*[linecolor=or](0.5,1.1)(2,1)(2.1,1.5)(0.4,1.7)
%
   \psline[linecolor=black, linewidth=1pt]{-}(2,1)(3.5,1.1)(3.4,2.4)(2.2,2)(2,1)
   \psline[linecolor=black, linewidth=1pt]{-}(3.5,1.1)(5,1.3)(5.2,2.5)(3.4,2.4)(3.5,1.1)
   \psline[linecolor=black, linewidth=1pt]{-}(0.4,1.7)(2.1,1.5)(2.2,2)(0.4,2.2)(0.4,1.7)
   \psline[linecolor=black, linewidth=1pt]{-}(0.5,1.1)(2,1)(2.1,1.5)(0.4,1.7)(0.5,1.1)
   \psline[linecolor=black, linewidth=1pt]{-}(2,0)(3.5,0.1)(3.5,1.1)(2,1)(2,0)
   \psline[linecolor=black, linewidth=1pt]{-}(2.2,2)(3.4,2.4)(3,3)(2,2.9)(2.2,2)
   \psline[linecolor=black, linewidth=1pt]{-}(3.4,2.4)(5.2,2.5)(5.1,3.1)(3,3)(3.4,2.4)
   \psline[linecolor=black, linewidth=1pt]{-}(0.4,2.2)(2.2,2)(2,2.9)(0,3.3)(0.4,2.2)
   \psline[linecolor=black, linewidth=1pt]{-}(0.5,1.1)(2,0)
   \psline[linecolor=black, linewidth=1pt]{-}(3.5,0.1)(5,1.3)
%
   \psline[linecolor=black, linewidth=1pt]{->}(0,0.6)(0.8,0.2) \rput[bl](0,0){${\bf F}$}
   \rput[bl](2.7,-1){(b)}
}
\end{pspicture}}\end{center}
\caption{Notations for the definition of the limitation process.
In orange, control volumes of the set $\mathcal V_\edge(V^-)$ for $\edge=V^-|V^+$, with a constant advection field ${\bf F}$: upwind cells (a) or opposite cells (b).}
\label{fig:MUSCL2}
\end{figure}
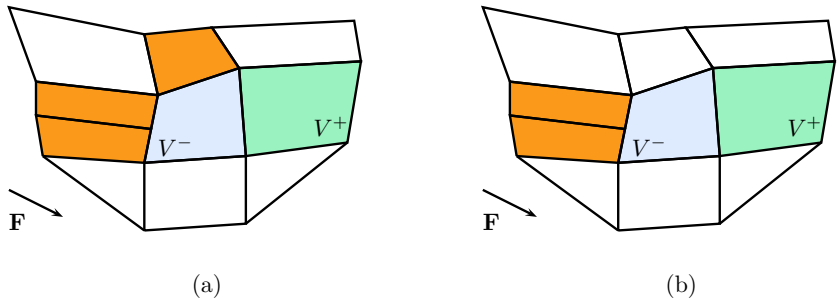

\medskip
Finally, we need to specify the choice of the sets $\mathcal V_\edge(V^-)$ and $\mathcal V_\edge(V^+)$.
Here, we just set $\mathcal V_\edge(V^+) = \{V^-\}$; such a choice guarantees that at least the upstream choice is in the intersection of the intervals defined by (H1) and (H2), as explained in Remark \ref{rm:MUSCL-1D}.
$\mathcal V_\edge(V^-)$ may be defined in two different ways (\cf\ Figure \ref{fig:MUSCL2}):
\begin{itemize}
\item[--] as the ``upstream cells'' to $V^-$, $\ie\ \mathcal V_\edge(V^-) = \{ L\in\mesh,\ L \text{ shares a face $\edge$ with $V^-$ and } F_{V^-,\edge}<0 \}$,
\item[--] when this makes sense (\ie\ with a mesh obtained by $Q_1$ mappings from the $(0,1)^d$ reference element), the opposite cells to $\edge$ in $V^-$ are chosen.
Note that for a structured mesh, this choice allows to recover the usual minmod limiter.
\end{itemize}

%
%
\subsection{An anti-diffusive scheme} \label{an:AD}

The scheme of Despr\'es-Lagouti\`ere \cite{dep-01-con} for the constant velocity advection problem presents some interesting proporties in one dimensional  (or structured multi-dimensional) space; in particular, it limits the numerical diffusion notably.
In this work, in order to treat convection operators of the chemical variables, which read in a simplified case
\[
\partial_t(\rho y) + \dive(\rho y \bfu) = 0,
\]
we used the following generalization, which allows us to work on untructured meshes:
\[
\forall K\in\mesh,\qquad y_K^{n+1} = \frac{\rho_K^{n-1}}{\rho_K^n} y_K^n + \frac{\delta t}{|K|}\frac1{\rho_K^n}\sum_{\edge \in \edges(K)} F_{K,\edge}^n y_\edge^n,
\]
where for $\edge=K|L$ and given that $F_{K,\edge}^n\geq0$, at first we estimate $y_\edge^n= y_L^n$ and second, to ensure stability, we project to the interval
\[
I= [y_K^n,\,y_K^n+\frac1d\frac{1-\nu}{\nu}(y_K-y_M),\quad \nu=\frac{F_{K,\edge}^{n+1}\,\delta t}{\rho_K^{n+1}\ |K|},
\]
where $M\in\mesh$ is the control volume which stands at the opposite side of $K$ with respect to $L$.
The scheme presented in \cite{dep-01-con} is recovered by this formulation for the one-dimensional constant velocity convection equation.
%
%
\section{Numerical tests}\label{sec:num}

At the continuous level, the boundedness of the chemical mass fractions formally implies that, when $\varepsilon \rightarrow 0$, the relaxed model converges to the asymptotic one.
Indeed, integrating any of the reactive species mass balance equations with respect to time and space, we observe that $||\dot\omega||_{L^1(\Omega\times(0,T))}$ tends to zero as $\varepsilon$, and thus two separate zones appear: a zone characterized by $G < 0.5$ where the reaction is complete, and a zone corresponding to $G \geq 0.5$, where no reaction has occured.

\medskip
A closed form of the solution of the Riemann problem for the asymptotic model is available \cite{bec-10-rea}.
In order to perform numerical tests, a Riemann problem with initial conditions such that the analytic solution has the profile presented in Figure \ref{fig:ansol} is chosen.
Moreover, the selected configuration imposes zero amplitude for the contact discontinuity and the left non linear wave, thus the solution consists of three different constant states: $\bfW_R^*, \bfW^{**}$ and $\bfW_R$.
The right state corresponds to a stoechiometric mixture of hydrogen and air (so the molar fractions of Hydrogen, Oxygen and Nitrogen are $2/7$, $1/7$ and $4/7$ respectively) at rest, at the pressure $p=9.9\,10^4$ Pa and the temperature $T=283^\circ$ K.
The velocity is supposed to be zero in the left state, which is sufficient to determine the solution.
Physically, speaking, supposing that the initial discontinuity lies at $x=0$, this situation corresponds to the left part of a (symmetrical) constant velocity plane deflagration starting at $x=0.$.
The flame velocity is $u_f= 63$ m/s  and the formation enthalpies are zero except for the product (\ie\ steam), with $\Delta h_{f,O}^0=-13.255\, 10^6$ J (Kg K)$^{-1}$.
The quantity $\rho_u$ is the analytical density in the intermediate state (so the total velocity of the flame brush is equal to the sum of $u_f$ and the material velocity on the right side of the reactive shock, see \cite{bec-10-rea}).
The computation is initialized by the analytical solution at $t=0.002$ and the final time is $t=0.005$.
The computational domain is the interval $(0,4.5)$.

\begin{center}
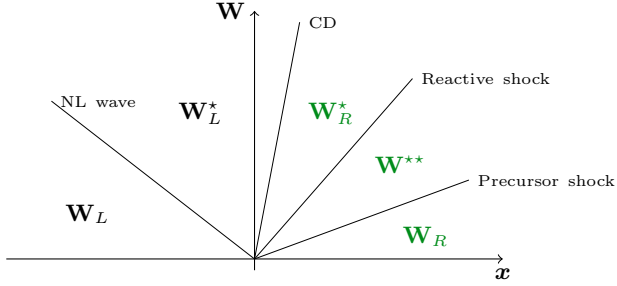
\begin{figure}[H]
\begin{tikzpicture}[scale=3]
    \draw[thin, ->] (-1.1,0) -- (1.1,0) node(xline)[below]
        {\small $\bfx$};
    \draw[thin, ->] (0,-0.05) -- (0,1.1) node(yline)[left] {\small $\bfW$};
    \draw[very thin] (0,0) coordinate (A) -- (0.95,0.35)
        coordinate (B) node[right, text width=6em] {\tiny Precursor shock};
    \draw[very thin] (0,0) coordinate (C) -- (0.7,0.8)
        coordinate (D) node[right, text width=6em] {\tiny Reactive shock};
    \draw[very thin] (0,0) coordinate (E) -- (0.2,1.05)
        coordinate (F) node[right, text width=6em] {\tiny CD};    
    \draw[very thin] (0,0) coordinate (G) -- (-0.9,0.7)
        coordinate (H) node[right, text width=6em] {\tiny NL wave};
    \node[left] at (0.48,0.65) {\small \textcolor{vertf}{$\bfW_R^\star$}};
    \node[left] at (0.80,0.42) {\small \textcolor{vertf}{$\bfW^{\star\star}$}};
    \node[left] at (-0.1,0.65) {\small $\bfW_L^\star$};
    \node[left] at (0.9,0.1) {\small \textcolor{vertf}{$\bfW_R$}};
    \node[left] at (-0.6,0.2) {\small $\bfW_L$};
\end{tikzpicture}
\vspace{-2ex}
\caption{The analytic solution of the numerical test configuration.} 
\label{fig:ansol}
\end{figure}
\end{center}

\medskip
The numerical tests performed aim at checking the convergence of the scheme to such a solution, which in fact may result from two different properties: the convergence of the relaxed model to the asymptotic model when $\varepsilon$ tends to zero, and the convergence of the scheme towards a numerical solution when the time and space steps tend to zero.
To this purpose, we choose $\varepsilon$ proportional to the space step and make it tend to zero, with a constant CFL number.
We test the scheme behaviour with three different discretizations of the convection operator in the chemical mass species balances: the standard upwind scheme, a MUSCL-like discretization which is an extension to variable density flows of the scheme proposed in \cite{pia-13-for} and is described in Section \ref{an:MUSCL}, and a first-order anti-diffusive scheme proposed in \cite{dep-01-con} and given in Section \ref{an:AD} for the sake of completeness.

\begin{figure}[tb]
\begin{center}
\scalebox{0.85}{\begin{minipage}{\textwidth}
\epsfig{file=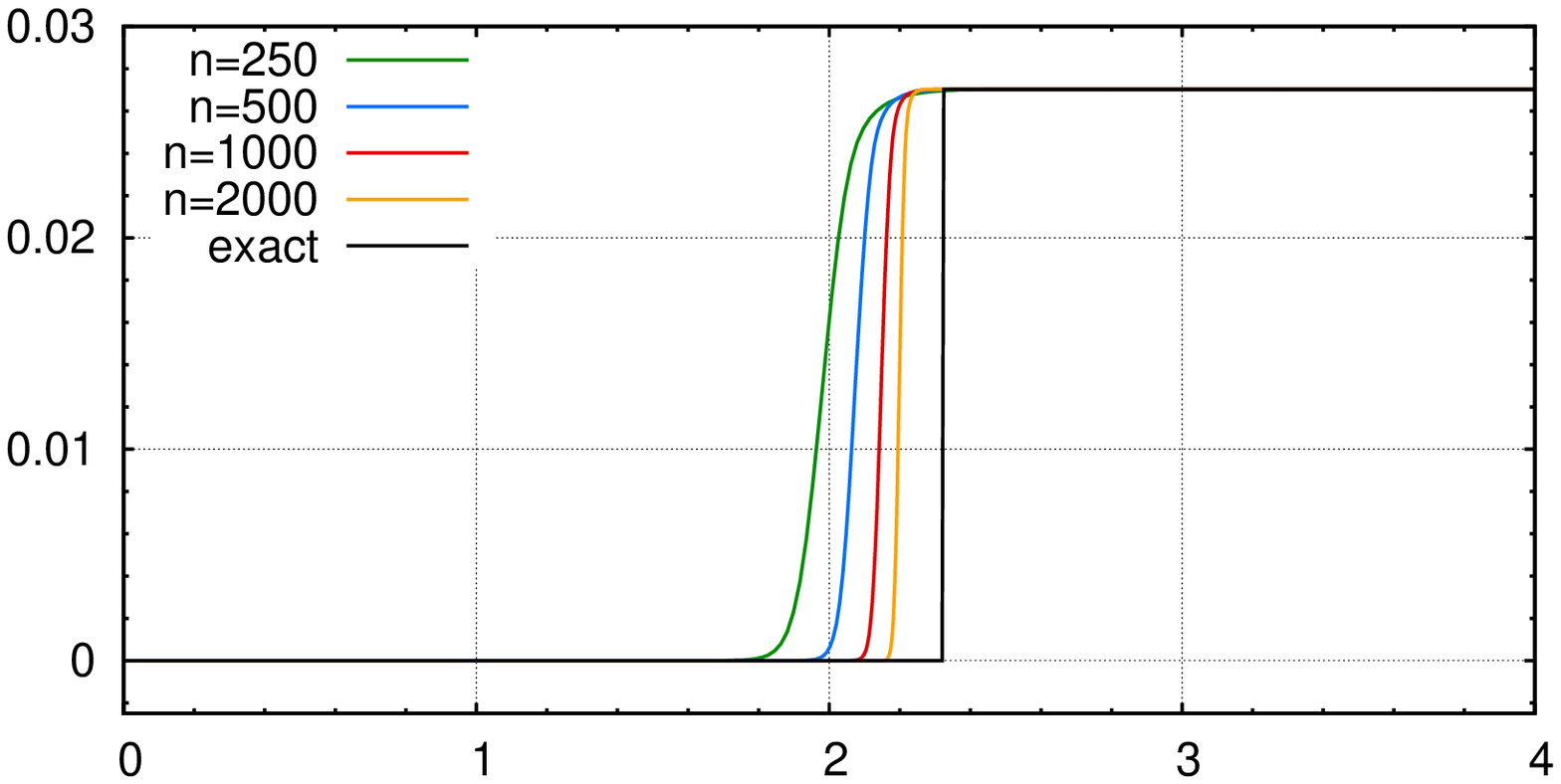,width=0.48\linewidth}
\epsfig{file=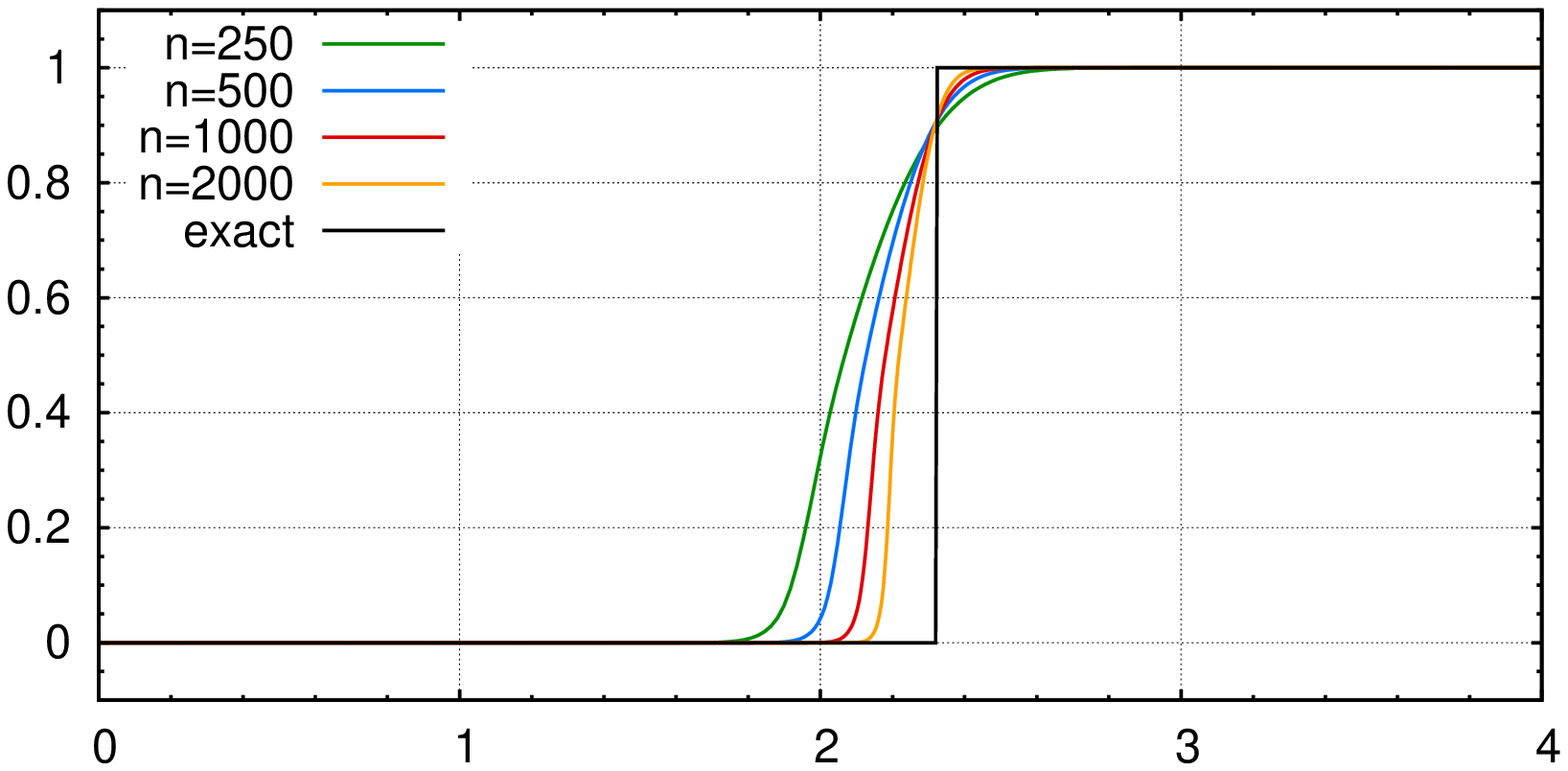,width=0.48\linewidth}\\
\epsfig{file=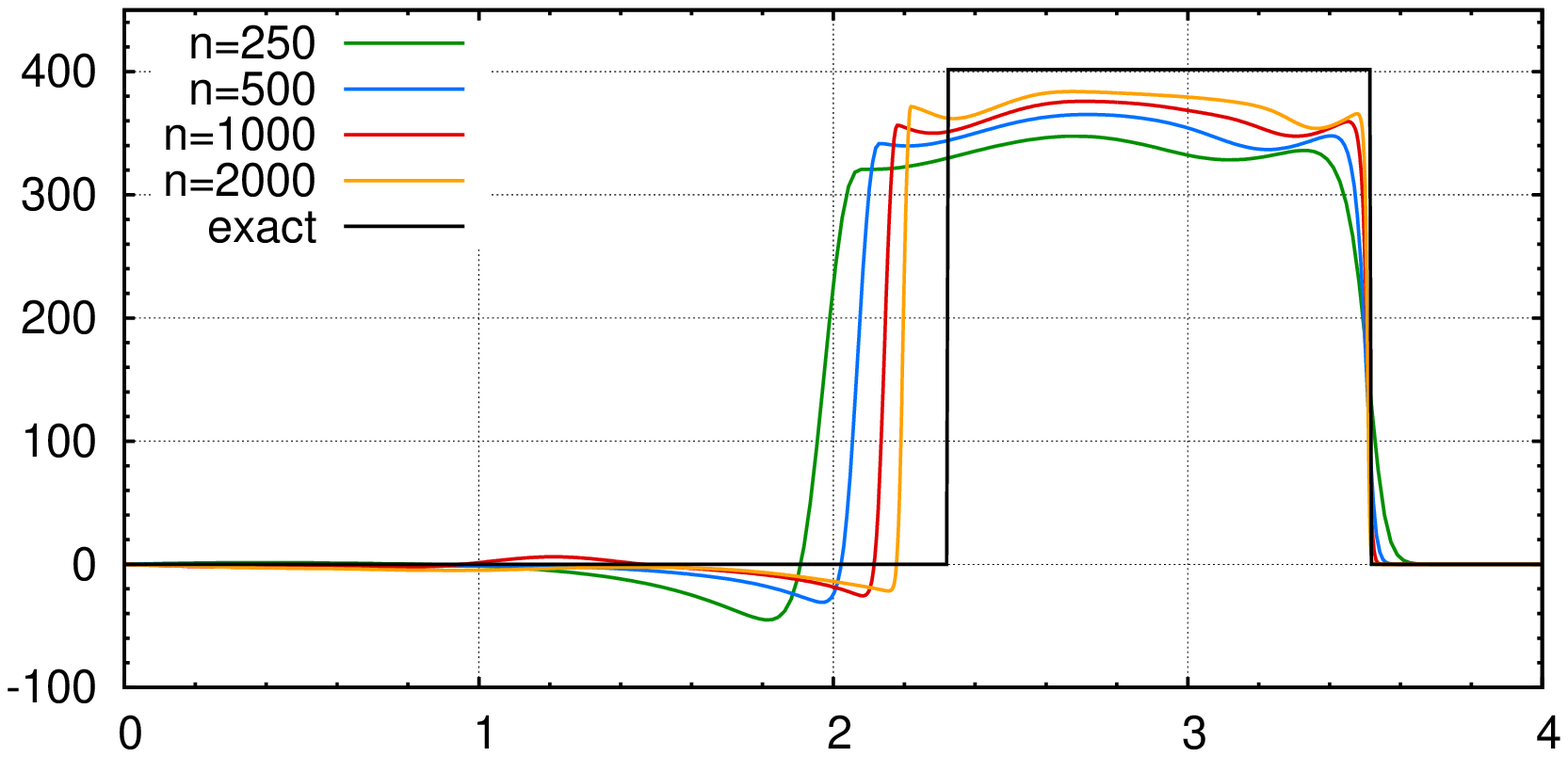,width=0.48\linewidth}
\epsfig{file=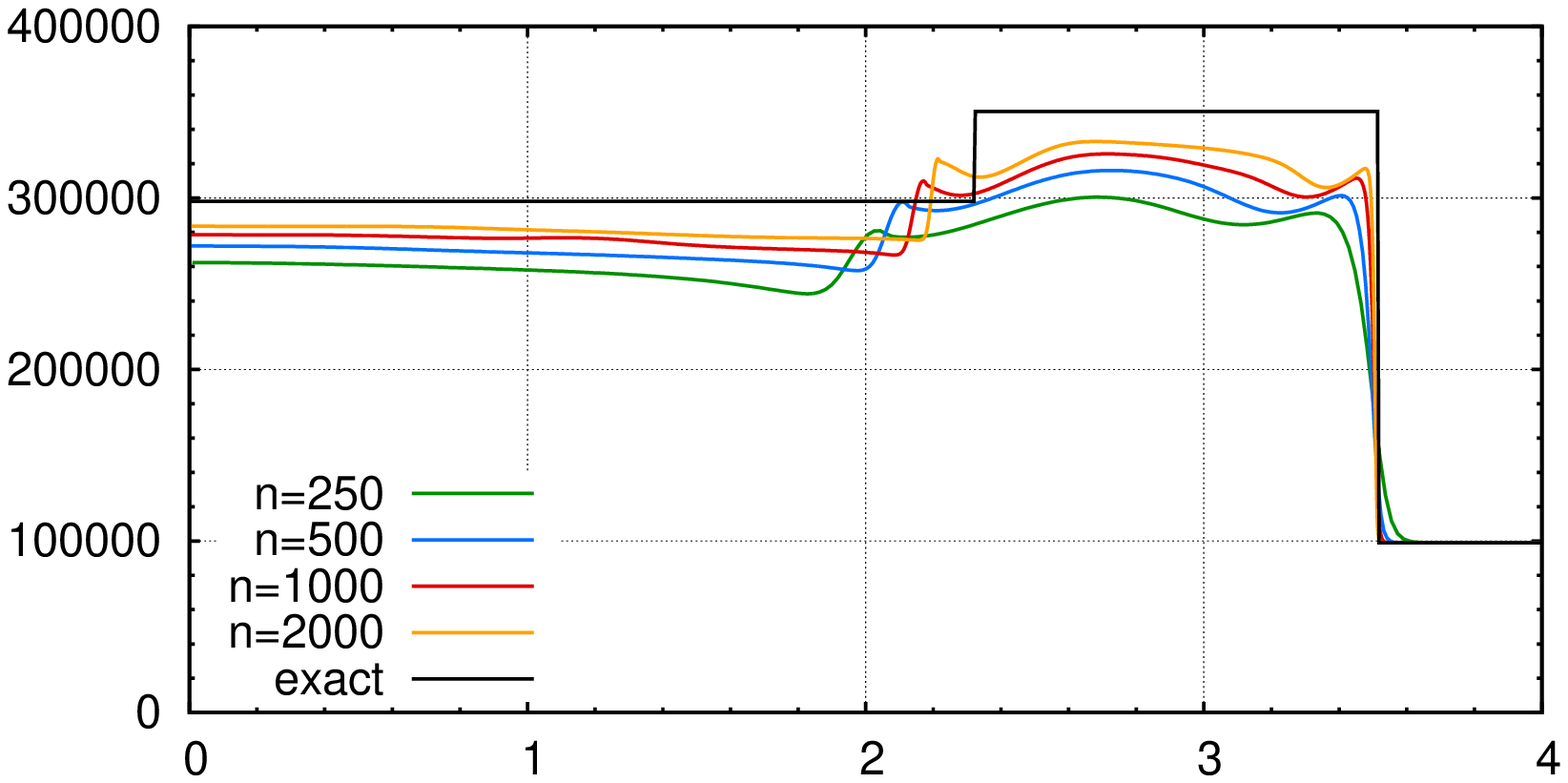,width=0.48\linewidth}\\
\epsfig{file=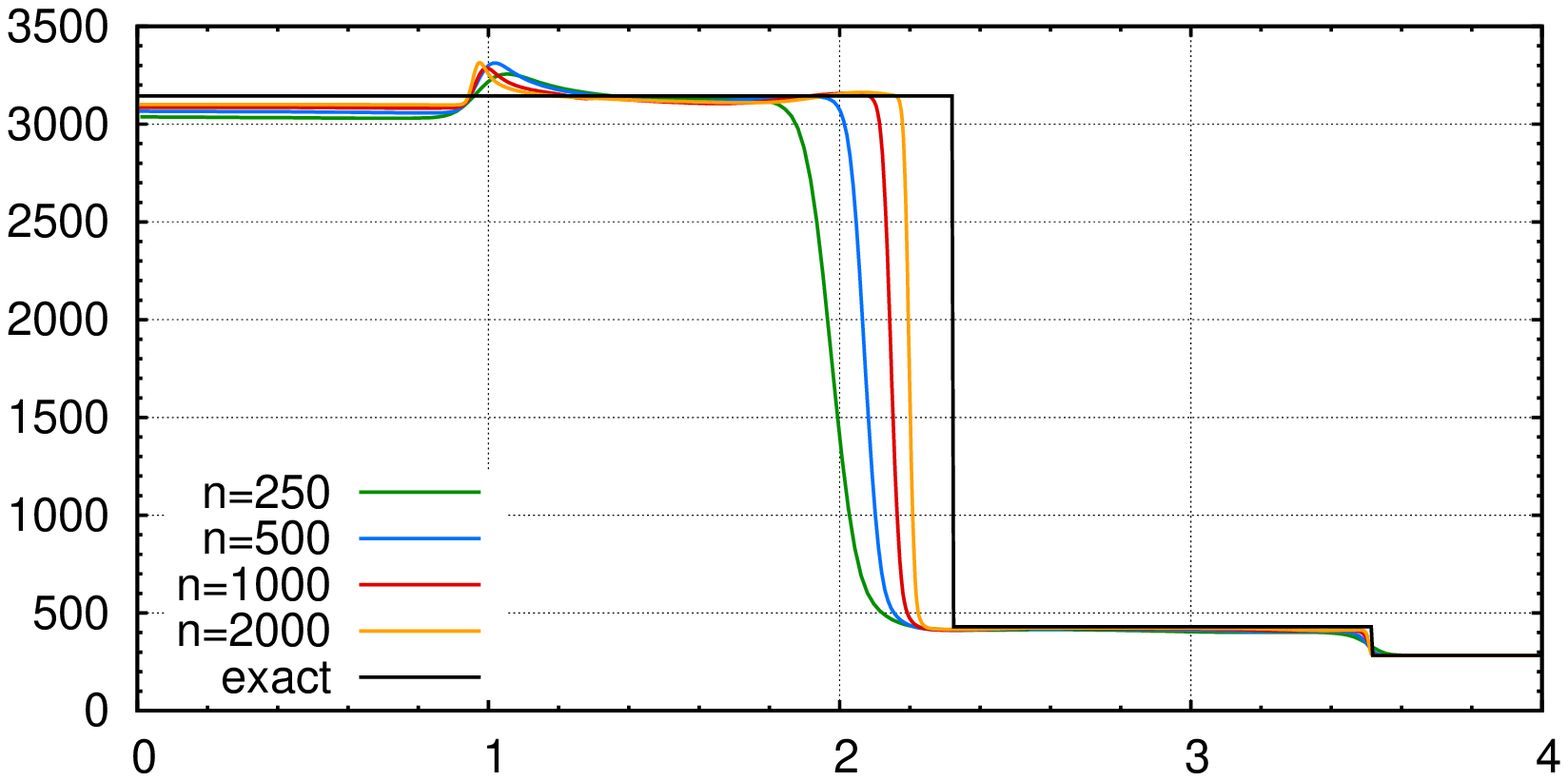,width=0.48\linewidth}
\epsfig{file=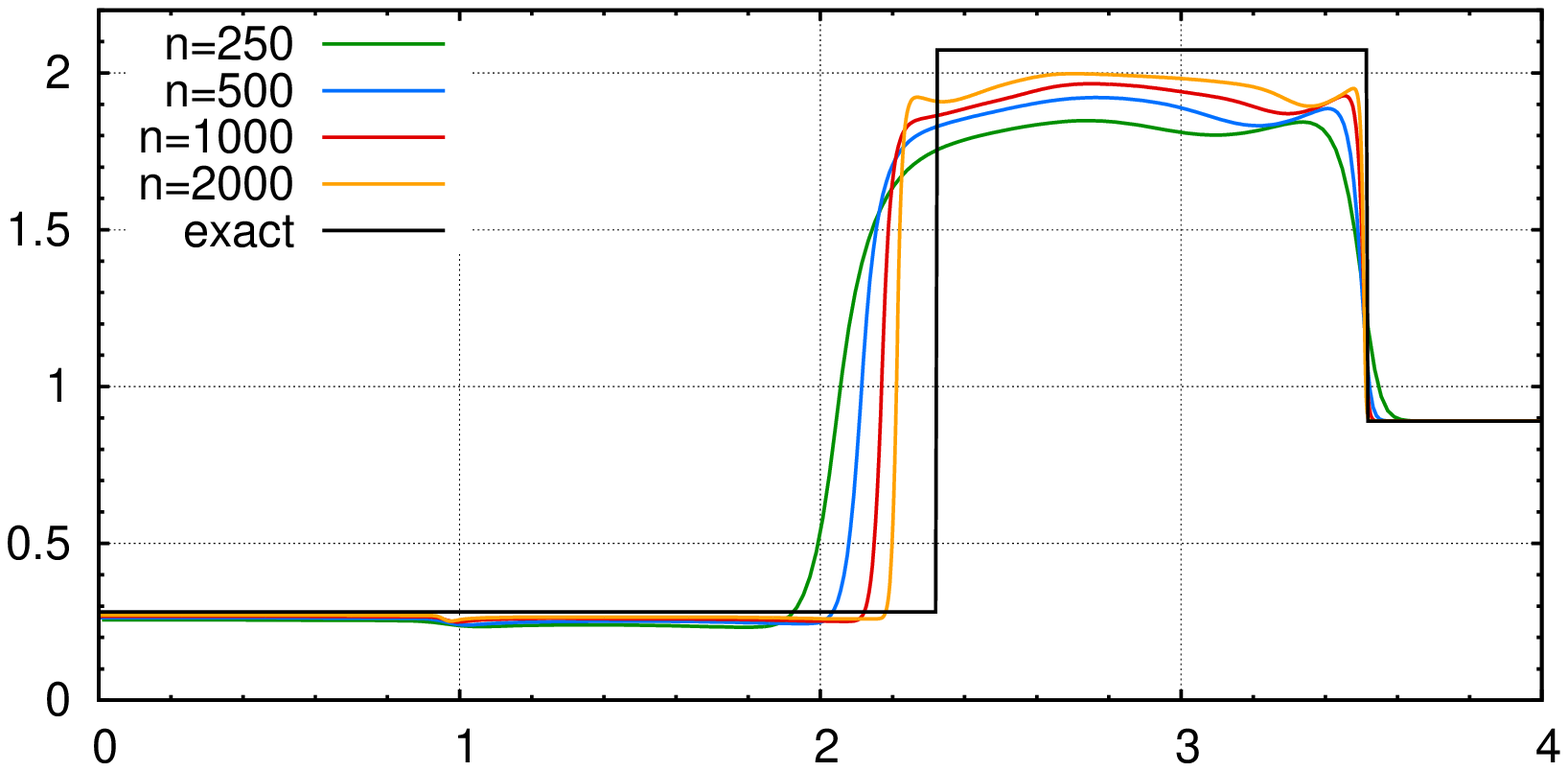,width=0.48\linewidth}
\end{minipage}}
\end{center}
\caption{Upwind scheme -- From top left to bottom right, fuel mass fraction, $G$, velocity, temperature and density at $t=0.005$, as a function of the space variable.
\label{fig:upwind}}
\end{figure}

\begin{figure}[tb]
\begin{center}
\scalebox{0.85}{\begin{minipage}{\textwidth}
\epsfig{file=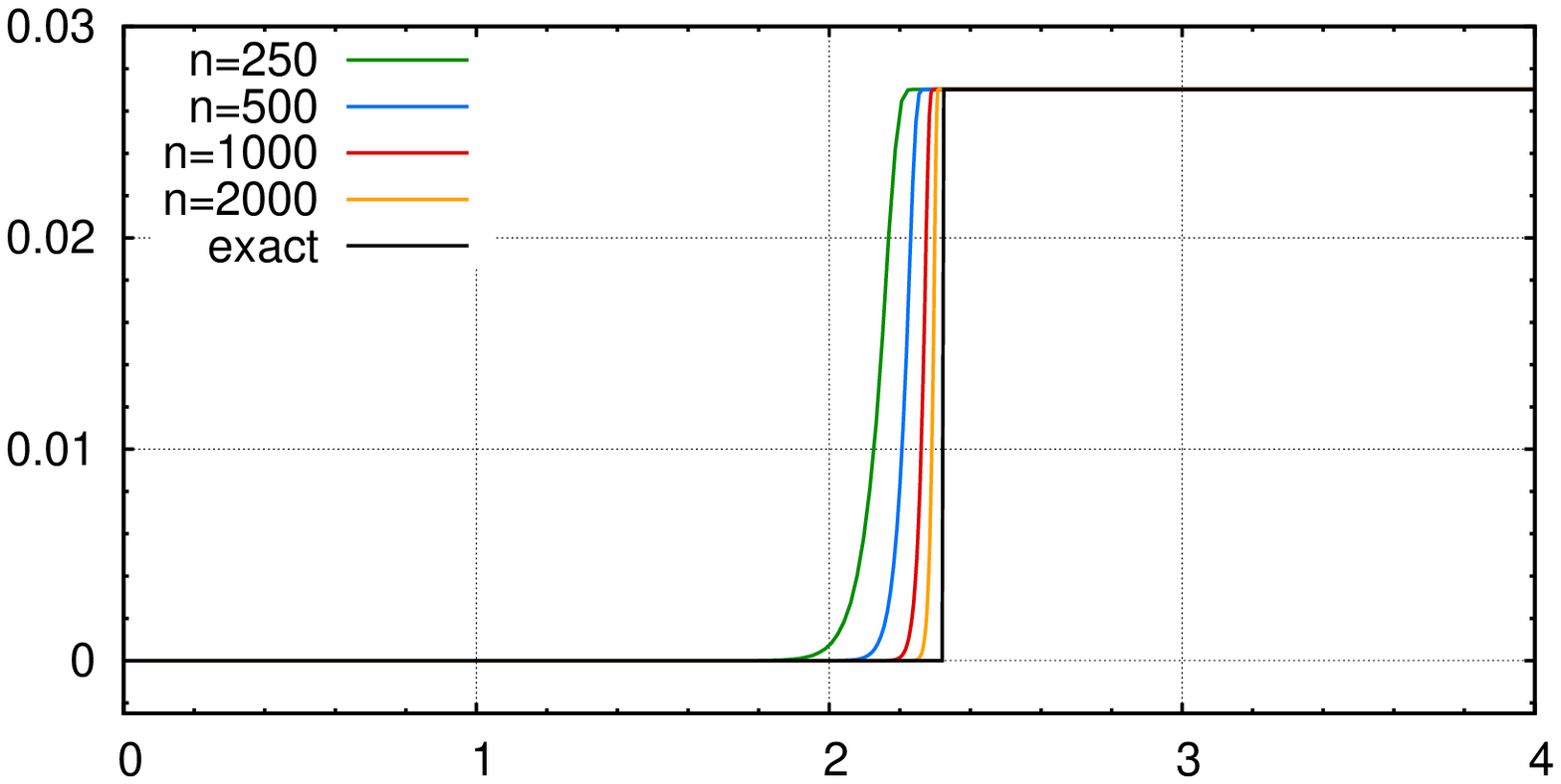,width=0.48\linewidth}
\epsfig{file=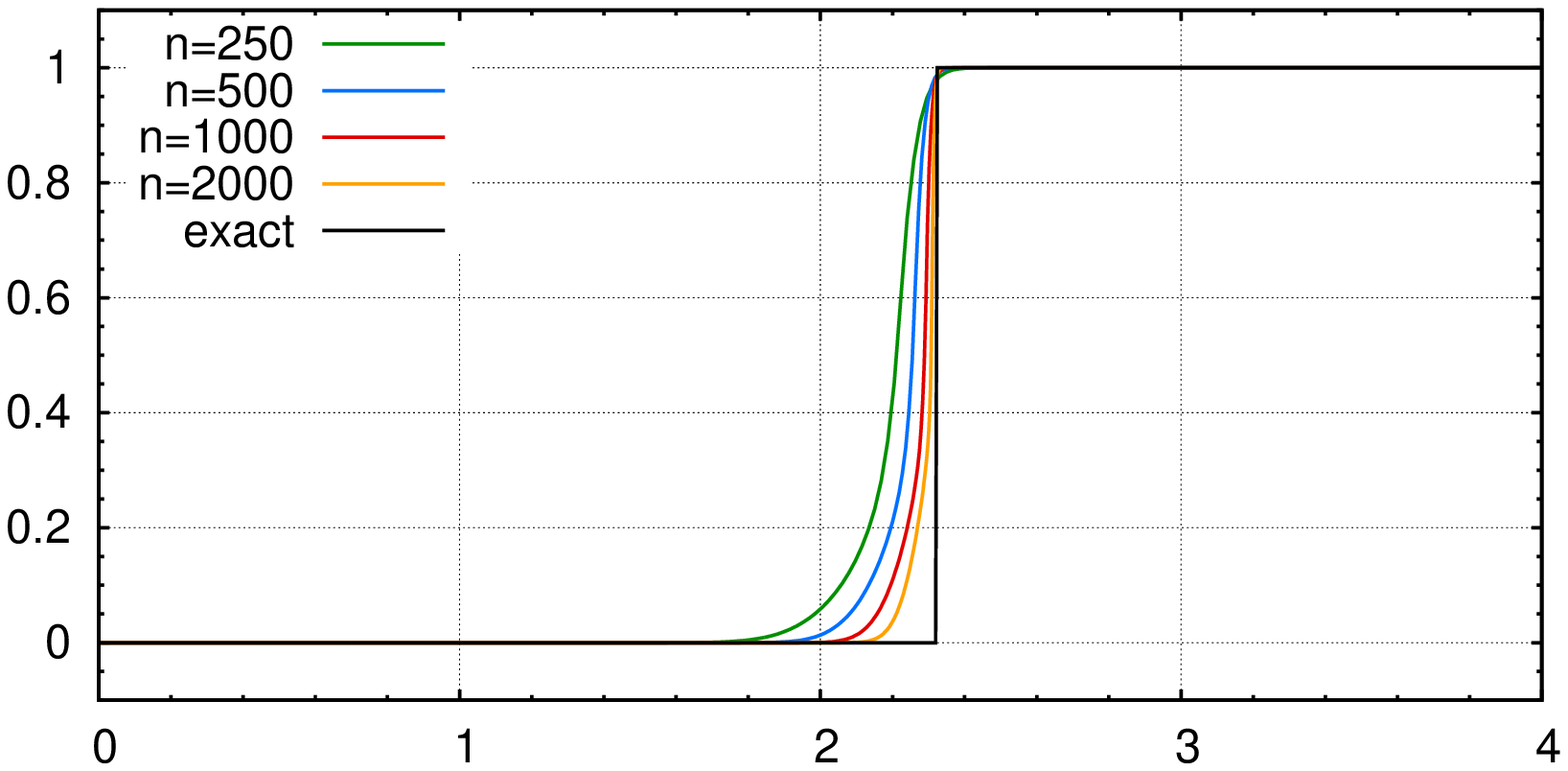,width=0.48\linewidth}\\
\epsfig{file=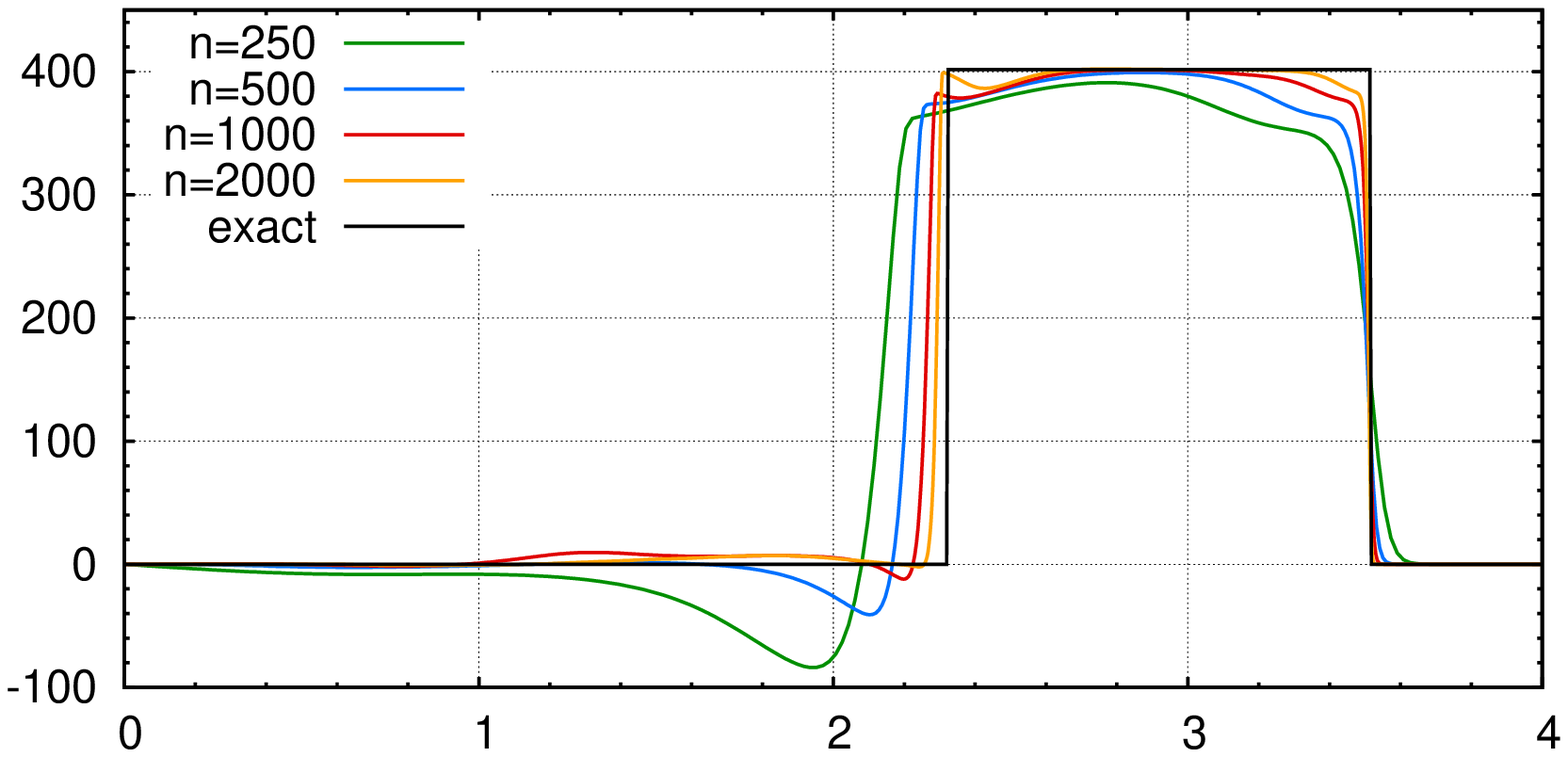,width=0.48\linewidth}
\epsfig{file=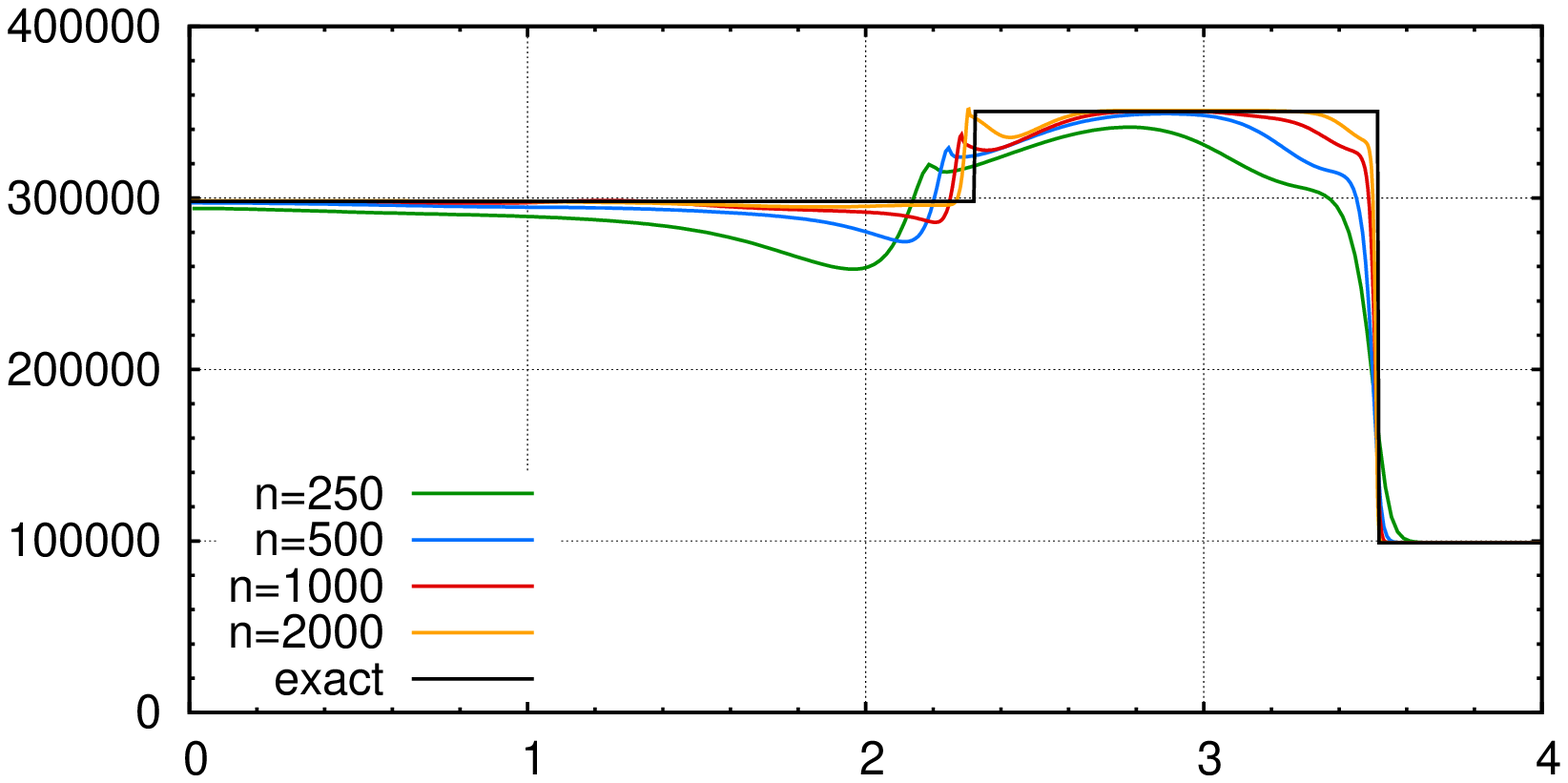,width=0.48\linewidth}\\
\epsfig{file=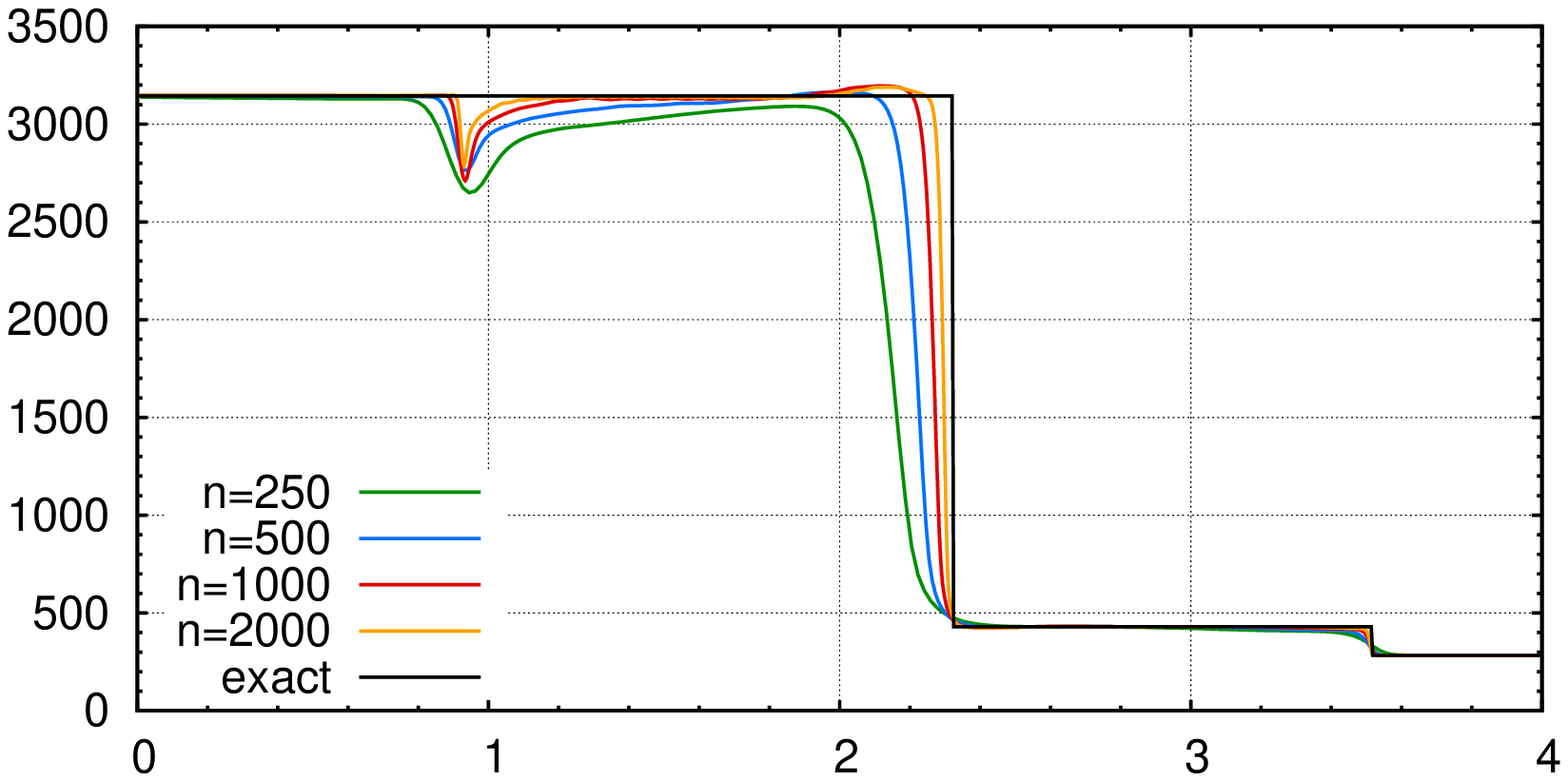,width=0.48\linewidth}
\epsfig{file=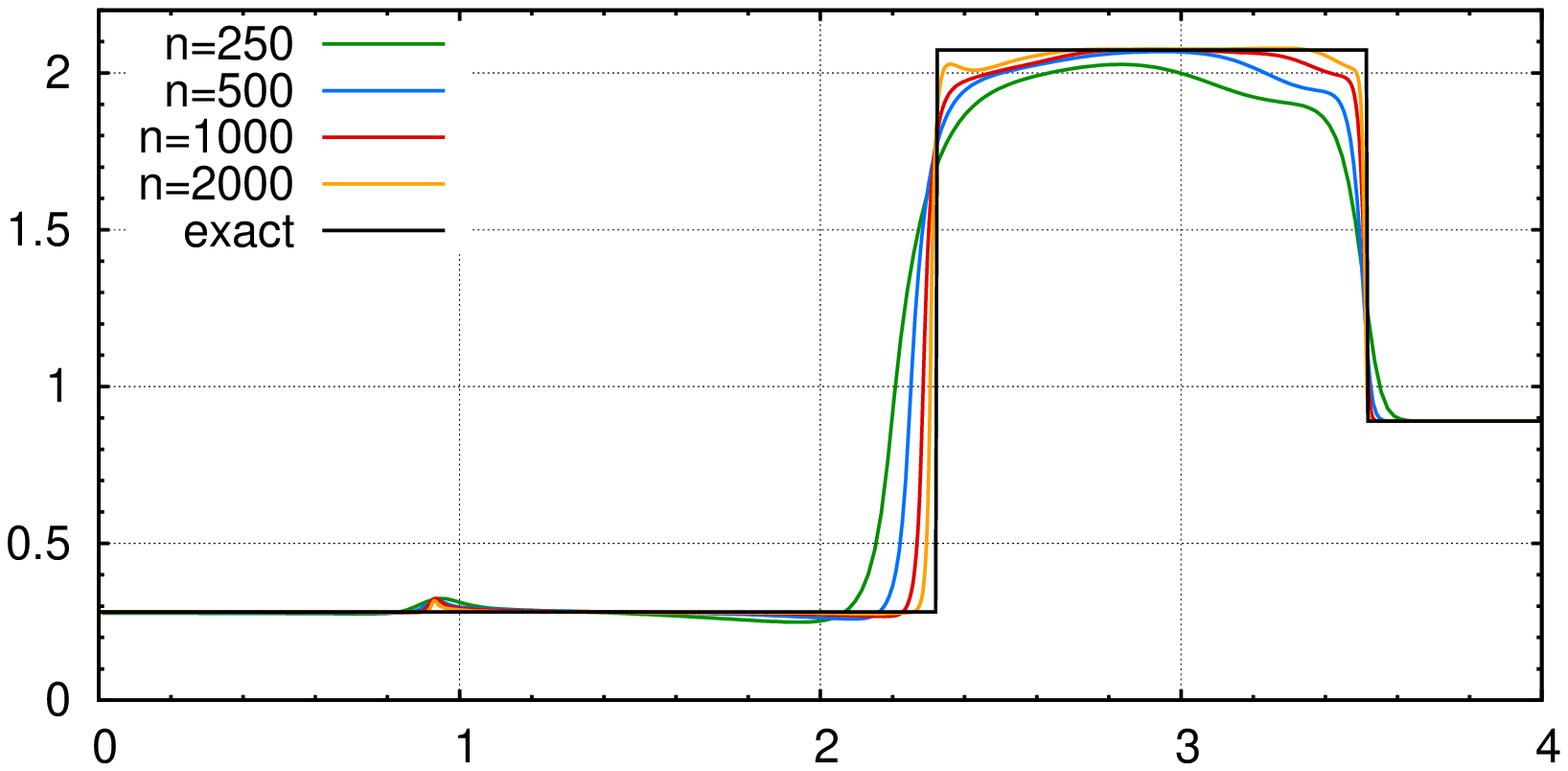,width=0.48\linewidth}
\end{minipage}}
\end{center}
\caption{MUSCL scheme -- From top left to bottom right, fuel mass fraction, $G$, velocity, temperature and density at $t=0.005$, as a function of the space variable.
\label{fig:muscl}}
\end{figure}

\begin{figure}[tb]
\begin{center}
\scalebox{0.85}{\begin{minipage}{\textwidth}
\epsfig{file=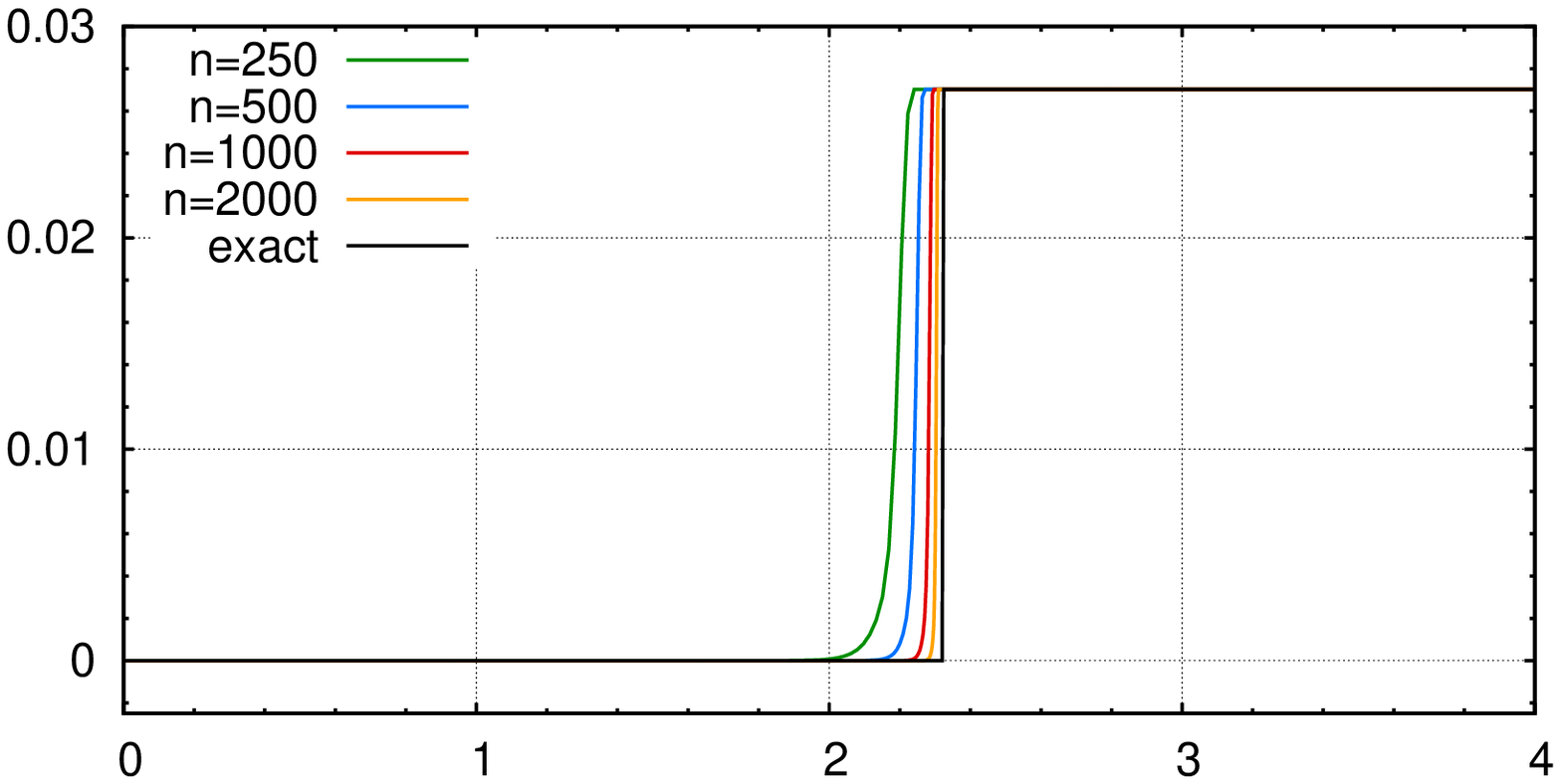,width=0.48\linewidth}
\epsfig{file=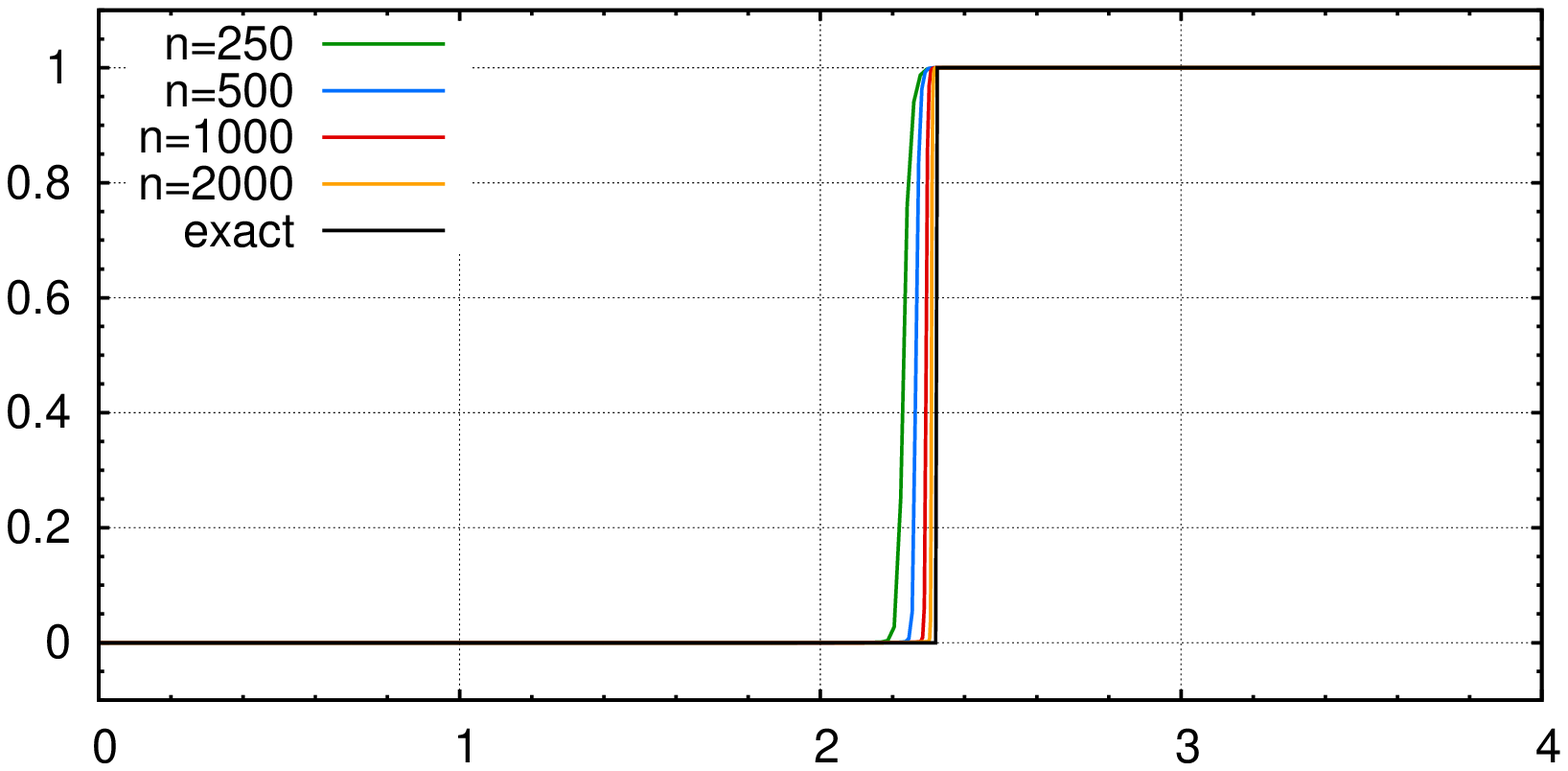,width=0.48\linewidth}\\
\epsfig{file=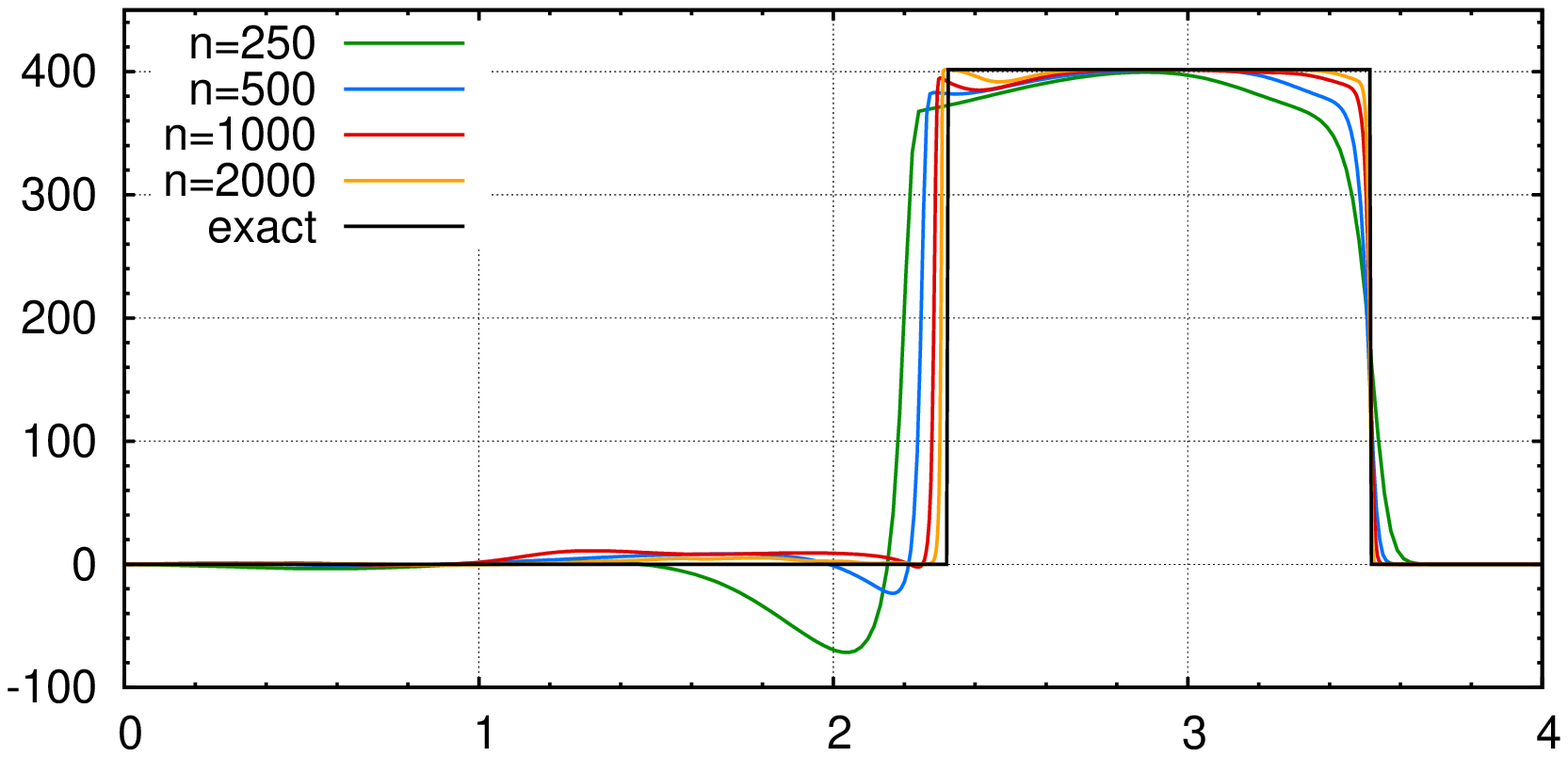,width=0.48\linewidth}
\epsfig{file=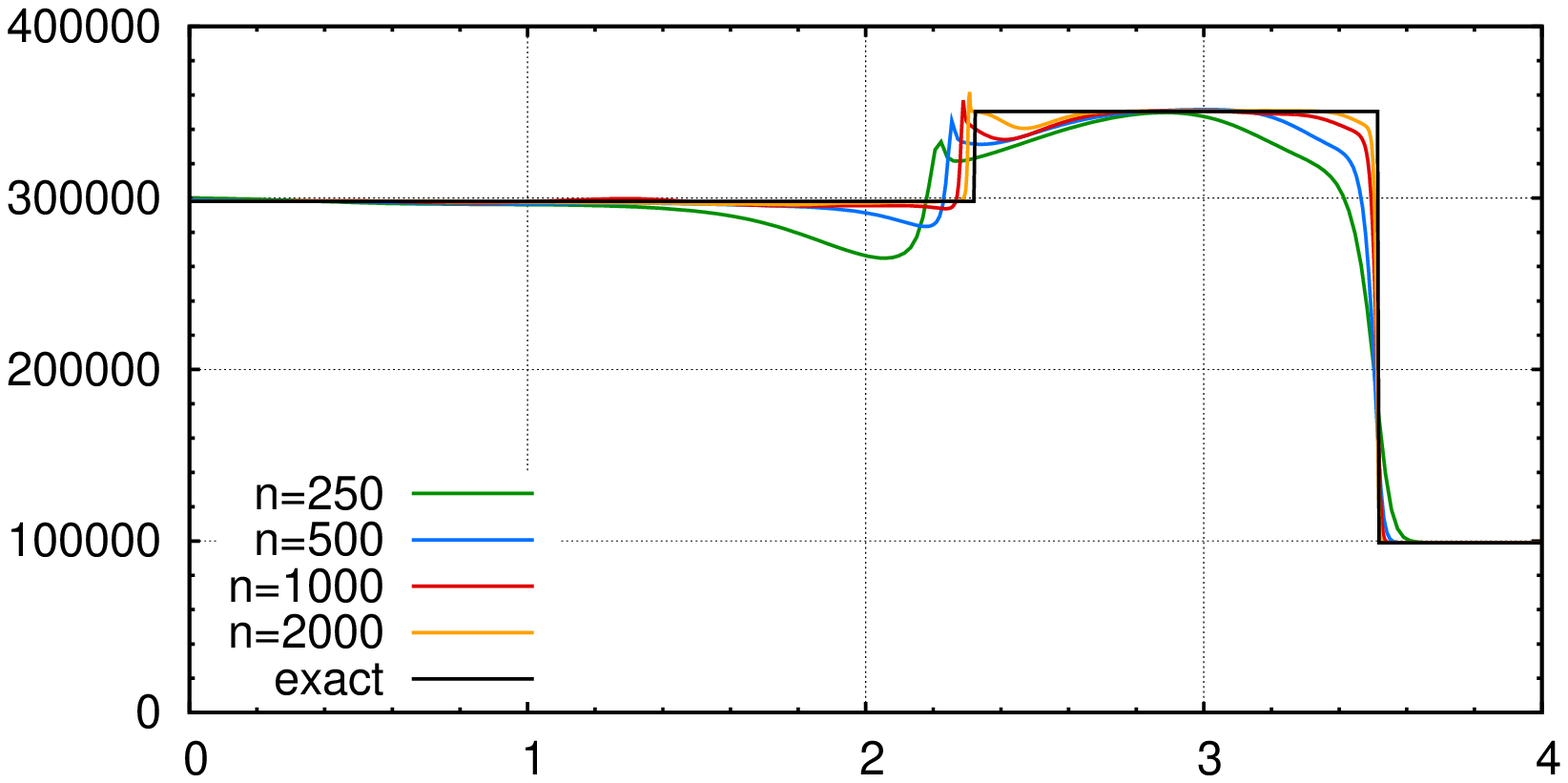,width=0.48\linewidth}\\
\epsfig{file=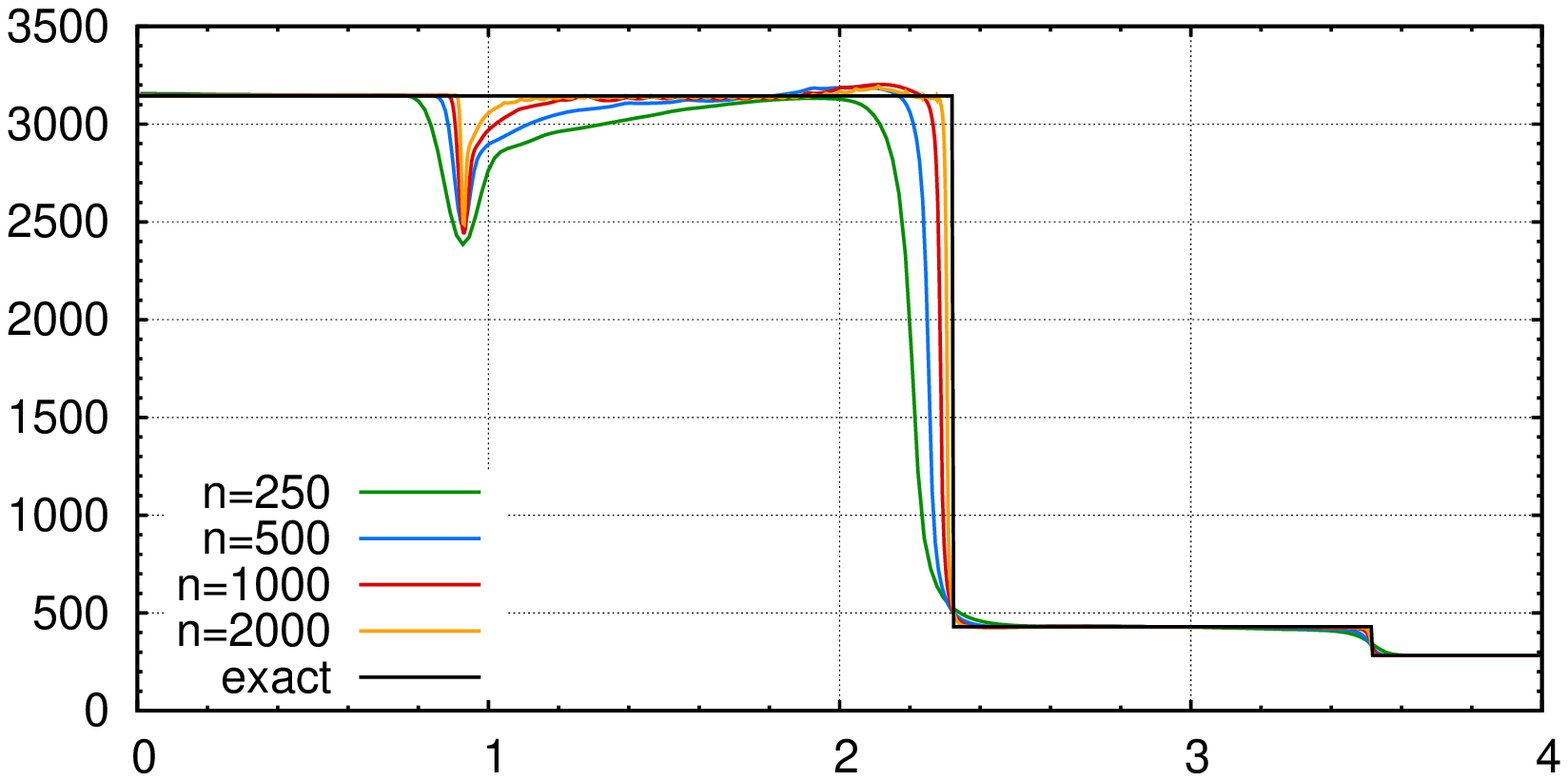,width=0.48\linewidth}
\epsfig{file=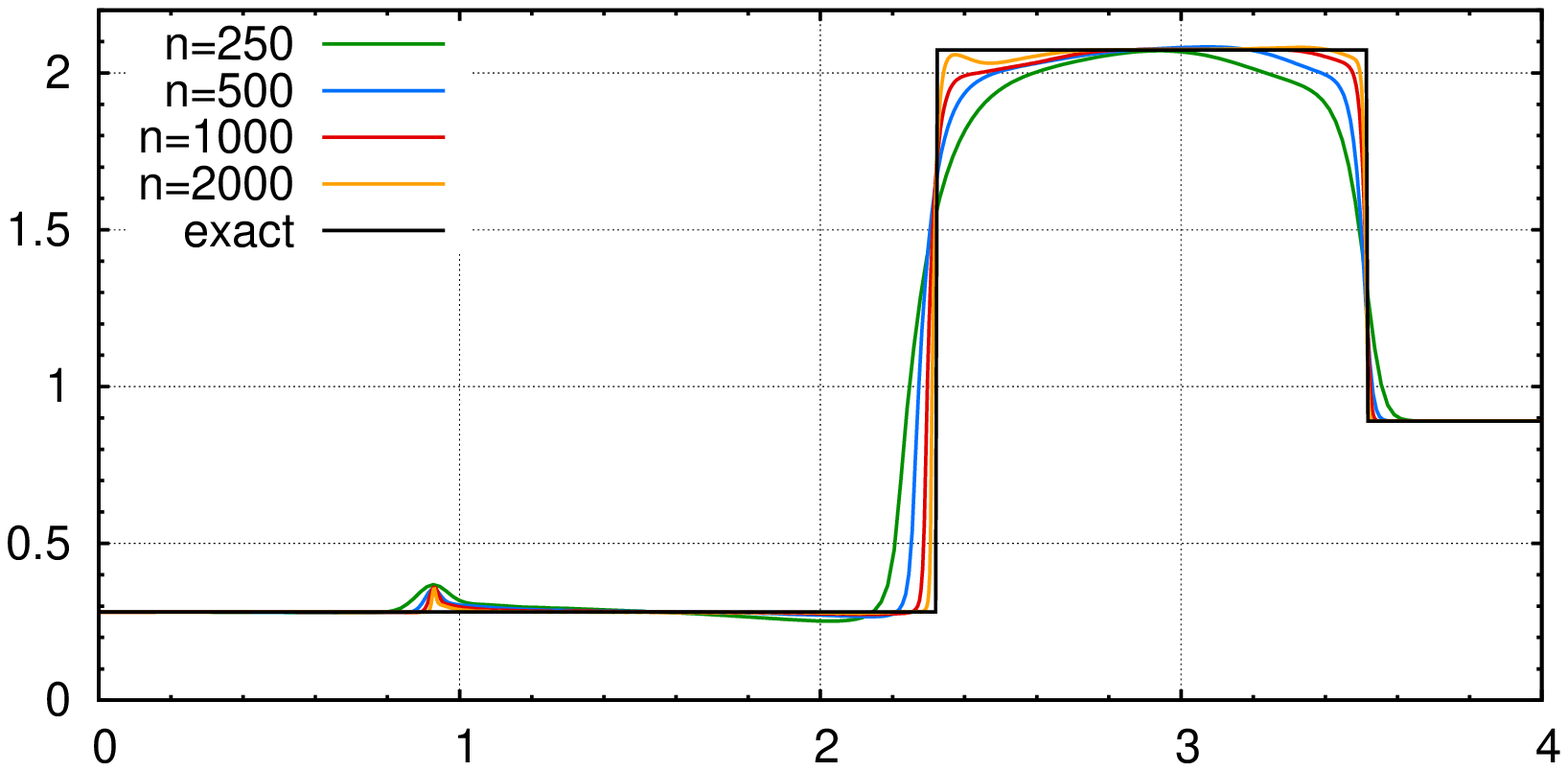,width=0.48\linewidth}
\end{minipage}}
\end{center}
\caption{Anti-diffusive scheme -- From top left to bottom right, fuel mass fraction, $G$, velocity, temperature and density at $t=0.005$, as a function of the space variable.
\label{fig:ad}}
\end{figure}

\begin{figure}[tb]
\begin{center}
\scalebox{0.85}{\begin{minipage}{\textwidth}
\epsfig{file=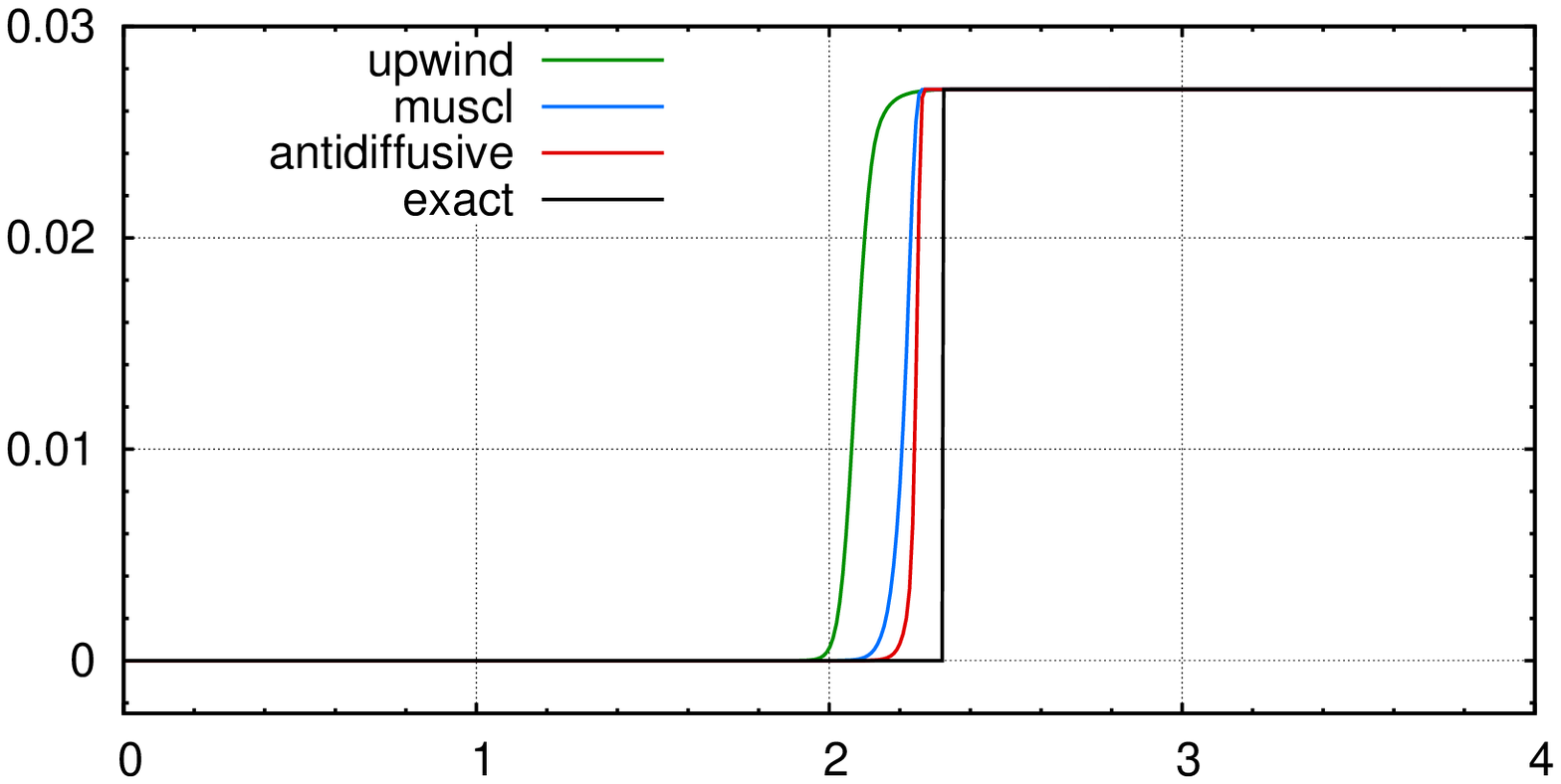,width=0.48\linewidth}
\epsfig{file=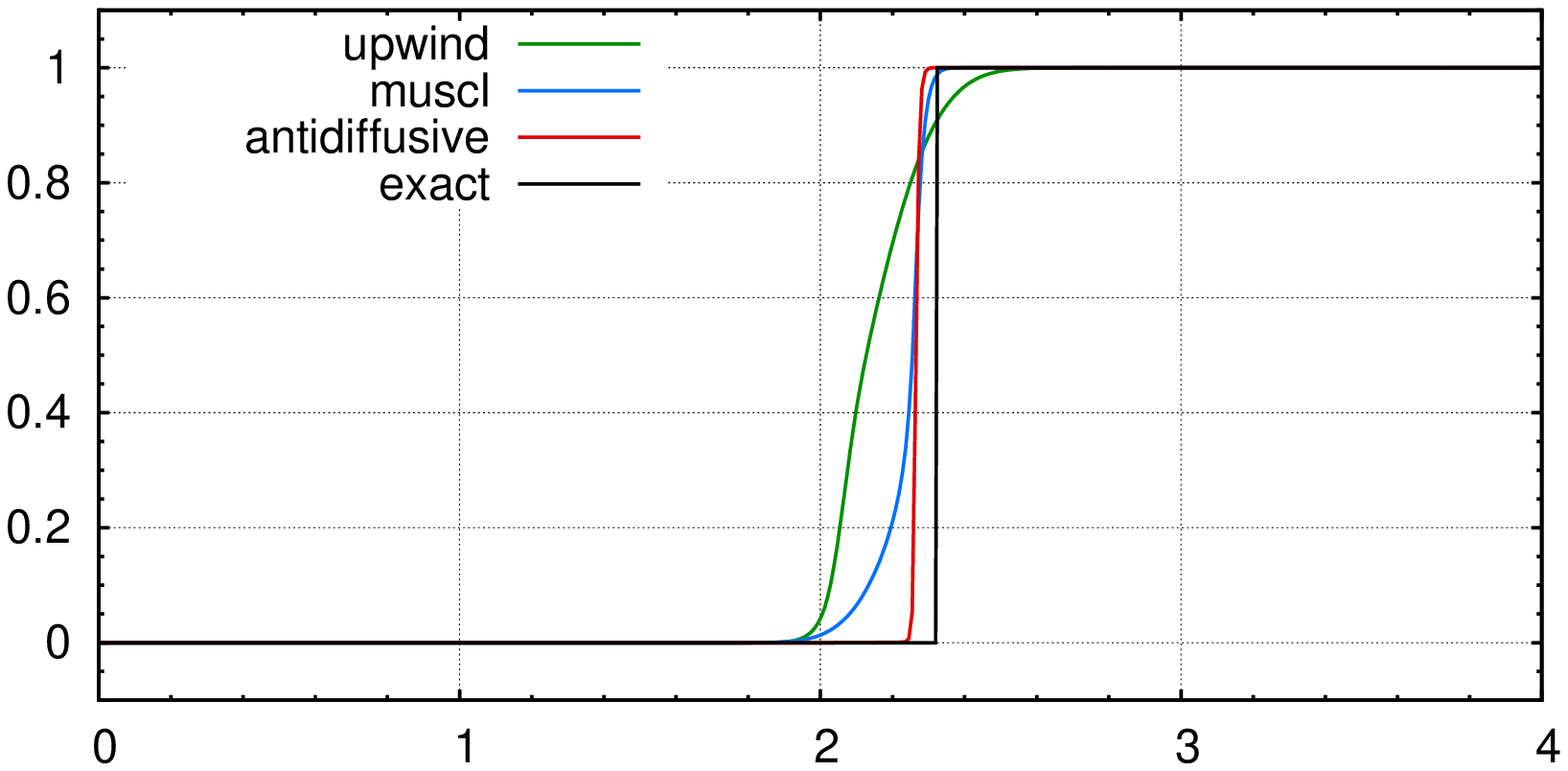,width=0.48\linewidth}\\
\epsfig{file=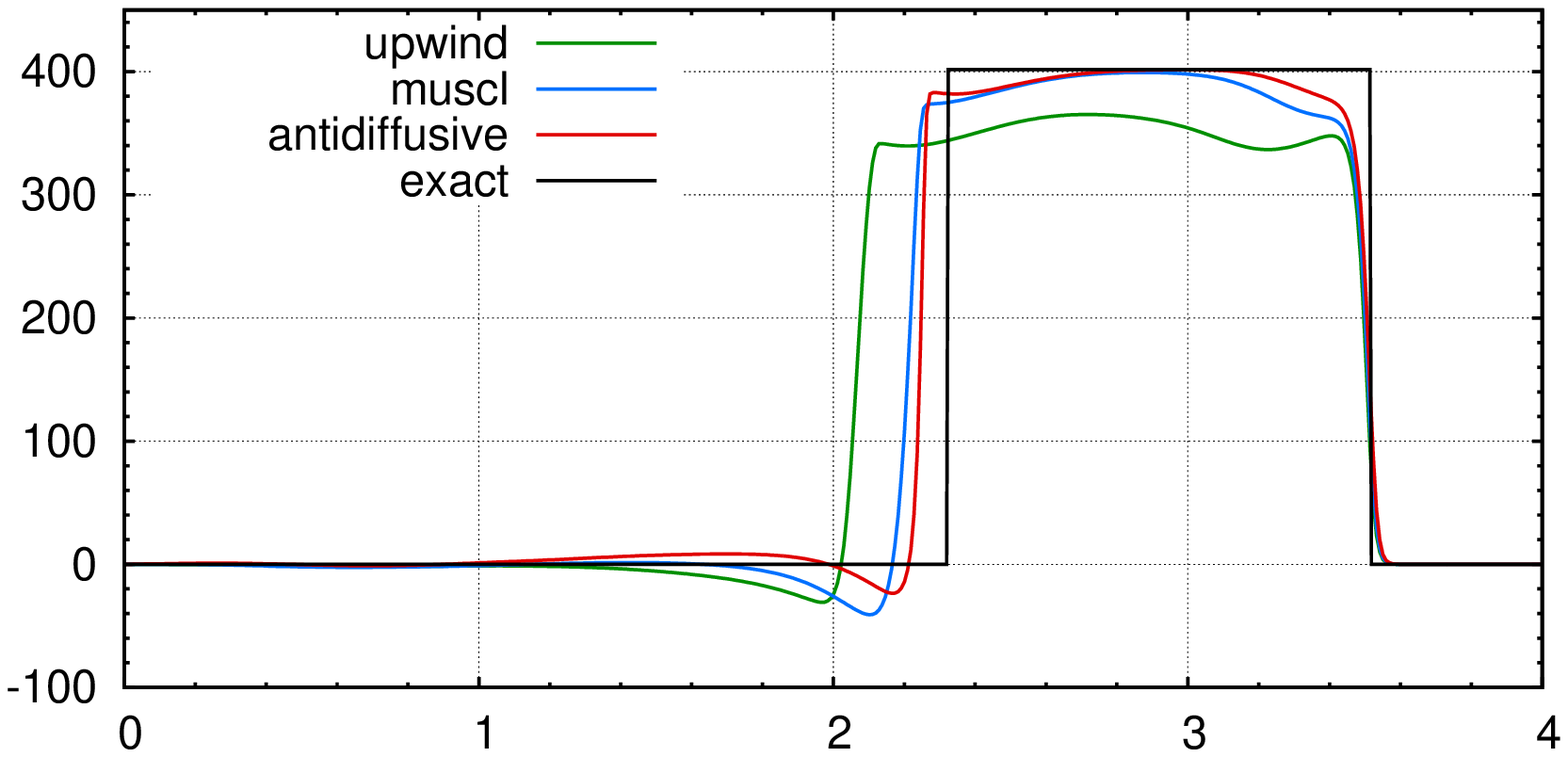,width=0.48\linewidth}
\epsfig{file=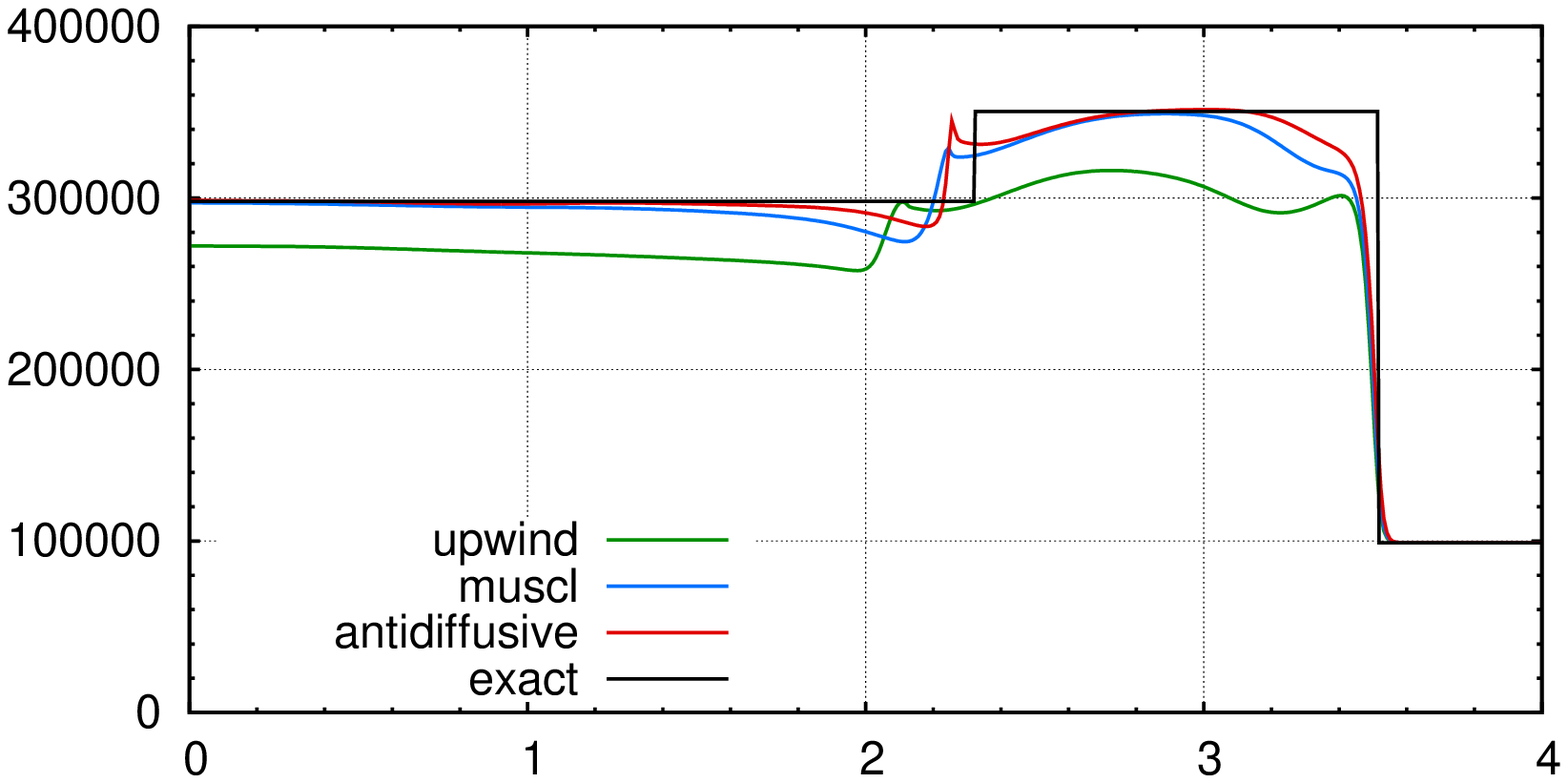,width=0.48\linewidth}\\
\epsfig{file=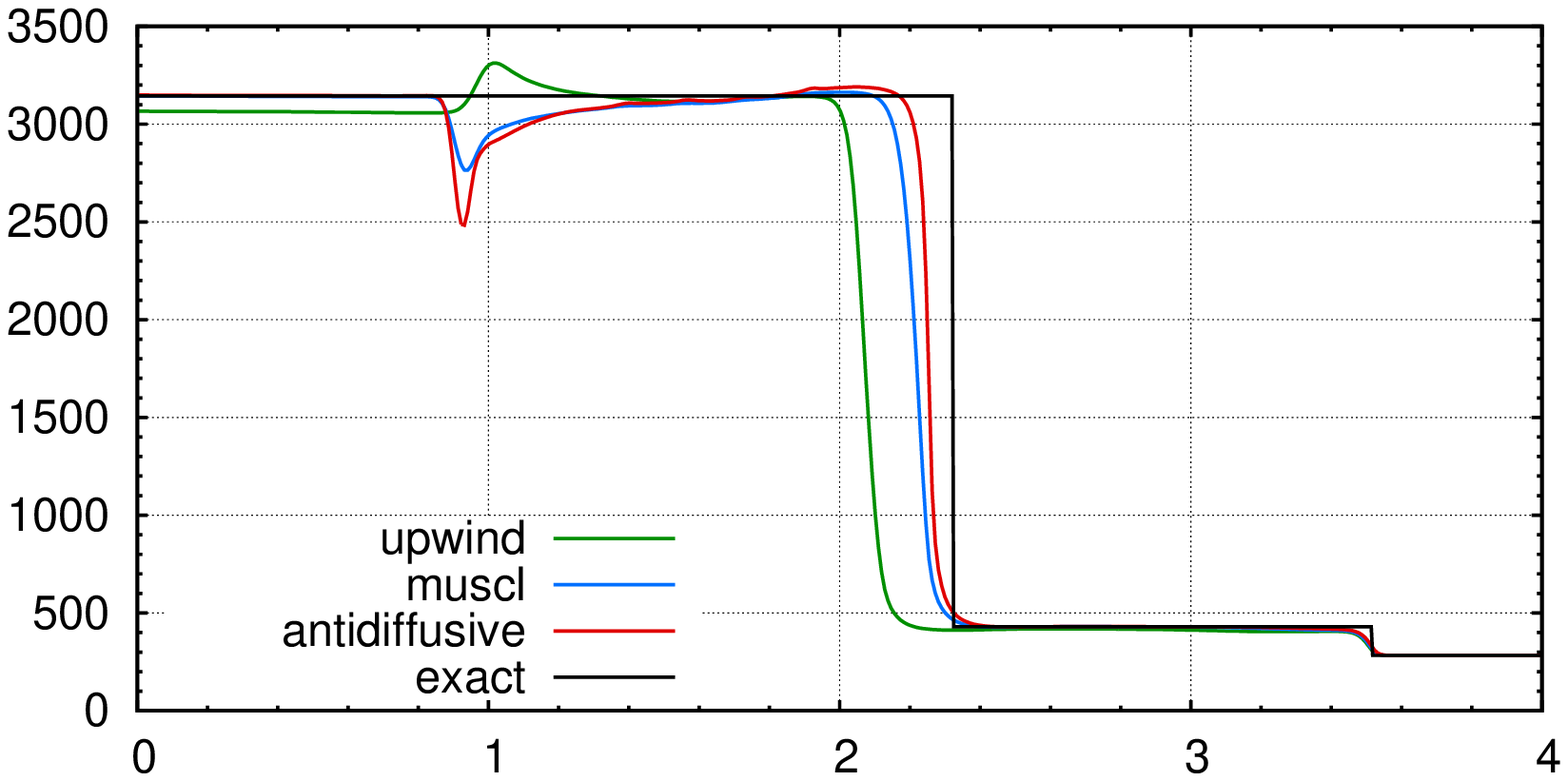,width=0.48\linewidth}
\epsfig{file=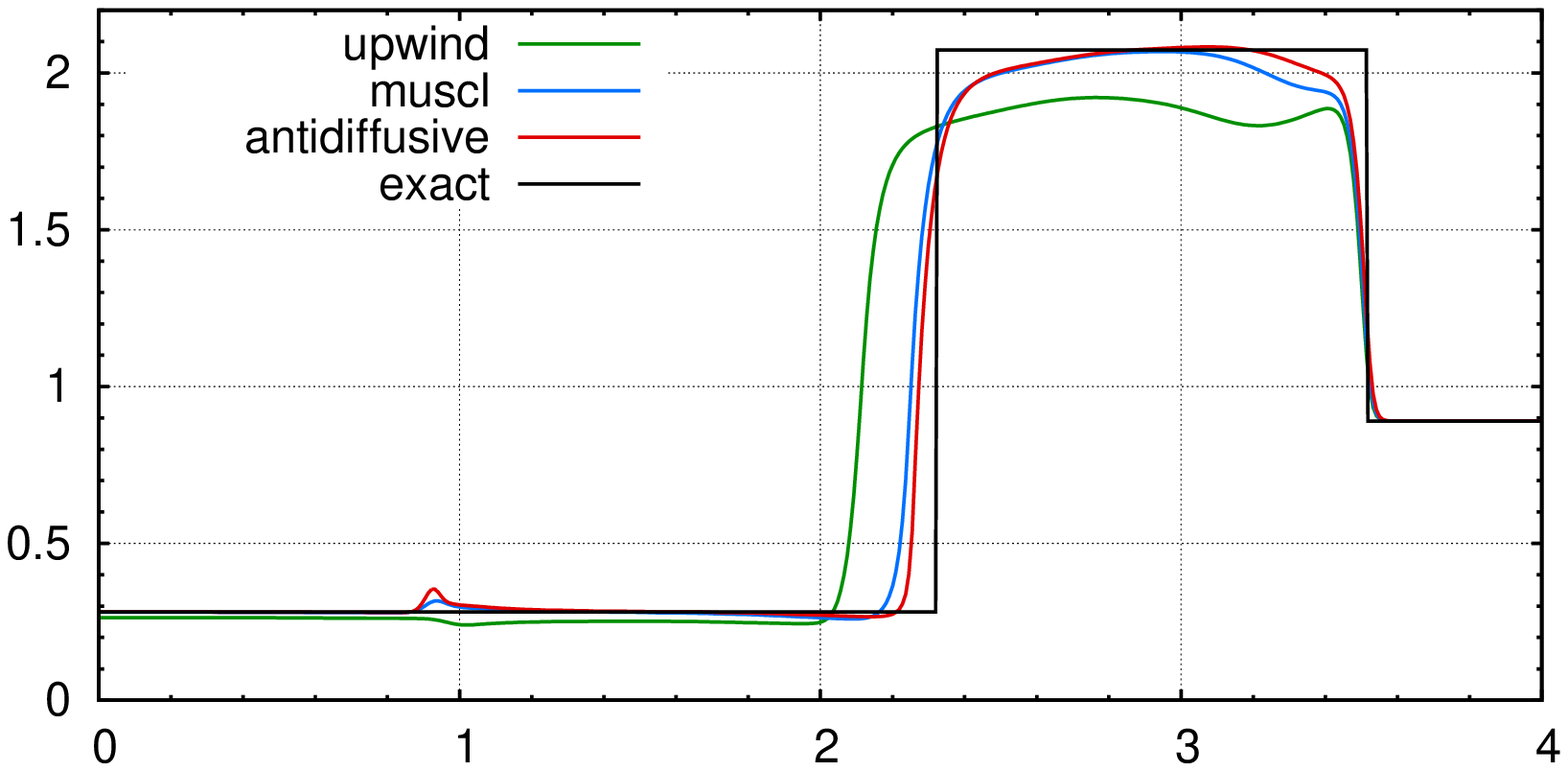,width=0.48\linewidth}
\end{minipage}}
\end{center}
\caption{Comparison of the solutions obtained with the upwind, MUSCL and anti-diffusive scheme -- From top to bottom, fuel mass fraction, $G$, velocity, temperature and density at $t=0.005$, as a function of the space variable.
Results obtained with a regular mesh composed of $n=500$ cells.
\label{fig:comp}}
\end{figure}

\medskip
Results obtained at $t=0.005$ with the upwind scheme, the MUSCL-like scheme and the anti-diffusive scheme, for increasingly refined meshes, are shown on Figure \ref{fig:upwind}, Figure \ref{fig:muscl} and Figure \ref{fig:ad} respectively, together with the analytical solution.
The expected convergence is indeed observed but, with the upwind discretization, the rate of convergence is poor.
This seems to be due to the interaction between the numerical diffusion of the upwind scheme, which artificially introduces unburnt reactive masses to the burnt zone, and the stiffness of the reaction term.
As expected in such a case, the results are significantly improved by the use of a less diffusive scheme for the chemical species balance equations.
Indeed, passing from the upwind to the MUSCL-like and to the anti-diffusive discretization improves the accuracy of the scheme, as may be observed in Figure \ref{fig:comp}, where the results obtained by the three discretizations for a regular mesh composed of 500 cells are plotted together with the continuous solution.
This observation is conforted by the measures, in $L^1$-norm, of the difference between the discrete and continuous solutions gathered in Table \ref{tab:errors}.
For every mesh and variable, the anti-diffusive scheme is the most accurate and the upwind one the least accurate.
The calculated orders of convergence are respectively close to 0.5 and 1 for the upwind scheme, on one part, and the MUSCL-like and anti-diffusive schemes, on the other part.

\begin{center}
\newcommand{\z}{\phantom{0}}
\begin{table}[htb]
\begin{tabular}{l||c|c|c||}
$h$       & $||p-p_{ex}||_{\mathrm L^1} \times 10^{-4}$   & $||\bfu-\bfu_{ex}||_{\mathrm L^1}\times10^{-2}$ & $ ||\rho-\rho_{ex}||_{\mathrm L^1} \times 10$  \\[0.2ex] \hline
      & upwind \qquad muscl   \qquad antidif     
          & $ 2.17\z \qquad \textcolor{bleuf}{1.56\z\z} \qquad \textcolor{rougec}{1.07\z\z} $     
          & $ 7.69   \qquad \textcolor{bleuf}{3.71\z}   \qquad \textcolor{rougec}{2.74\z}   $
\\
$h_0$     & $ 16.5\z \qquad \textcolor{bleuf}{7.26\z}   \qquad \textcolor{rougec}{4.59\z}   $  
          & $ 2.17\z \qquad \textcolor{bleuf}{1.56\z\z} \qquad \textcolor{rougec}{1.07\z\z} $     
          & $ 7.69   \qquad \textcolor{bleuf}{3.71\z}   \qquad \textcolor{rougec}{2.74\z}   $
\\
$h_0/2$   & $ 12.5\z \qquad \textcolor{bleuf}{3.88\z}   \qquad \textcolor{rougec}{2.43\z}   $  
          & $ 1.64\z \qquad \textcolor{bleuf}{0.787\z}  \qquad \textcolor{rougec}{0.579\z}  $     
          & $ 6.16   \qquad \textcolor{bleuf}{2.23\z}   \qquad \textcolor{rougec}{1.65\z}   $
\\
$h_0/4$   & $ \z9.66 \qquad \textcolor{bleuf}{2.05\z}   \qquad \textcolor{rougec}{1.38\z}   $  
          & $ 1.23\z \qquad \textcolor{bleuf}{0.471\z}  \qquad \textcolor{rougec}{0.371\z}  $
          & $ 4.73   \qquad \textcolor{bleuf}{1.26\z}   \qquad \textcolor{rougec}{0.913}    $
\\
$h_0/8$   & $ \z7.58 \qquad \textcolor{bleuf}{1.17\z}   \qquad \textcolor{rougec}{0.708}    $  
          & $ 0.958  \qquad \textcolor{bleuf}{0.263\z}  \qquad \textcolor{rougec}{0.175\z}  $
          & $ 3.63   \qquad \textcolor{bleuf}{0.691}    \qquad \textcolor{rougec}{0.476}    $
\\
$h_0/20$  & $ \z5.78 \qquad \textcolor{bleuf}{0.673}    \qquad \textcolor{rougec}{0.375}    $
          & $ 0.728  \qquad \textcolor{bleuf}{0.160\z}  \qquad \textcolor{rougec}{0.103\z}  $
          & $ 2.77   \qquad \textcolor{bleuf}{0.382}    \qquad \textcolor{rougec}{0.267}    $
\\
$h_0/40$  & $ \z4.31 \qquad \textcolor{bleuf}{0.414}    \qquad \textcolor{rougec}{0.194}    $
          & $ 0.543  \qquad \textcolor{bleuf}{0.0786}   \qquad \textcolor{rougec}{0.0458}   $
          & $ 2.03   \qquad \textcolor{bleuf}{0.201}    \qquad \textcolor{rougec}{0.134}    $     
\end{tabular}
\caption{Error in $L^1$ norm of the error between the discrete and continuous solutions for the various schemes - Black values correspond to the errors obtained with the upwind discretization, blue ones to the MUSCL discretization and the orange ones to the anti-diffusive scheme; we denote by $h_0=4.5/250$ the measure of the control volume of the least refined mesh.}
\label{tab:errors}\end{table}
\end{center}

\clearpage

\clearpage
\appendix
%
%
\section{From the Lagrange-projection scheme to the downwind-limited scheme}

\subsection{The one-dimensional constant-velocity problem}

\paragraph{\textbf{The scheme}}
We begin with the following model problem:
\[
\partial_t y + v\ \partial_x y =0, \qquad v \in \xR,\ v \geq 0.
\]
Let $y_-$, $y$ and $y_+$ be the value taken by the unknown in three successive cells, sorted from left to right.
We denote by $K$ the middle mesh and by $\edge$ the interface between $K$ and the right cell (see Figure \ref{fig:1D_notations}).
The numerical flux through $\edge$ outward $K$ reads $G_{K,\edge} = v\, y_\edge$ and our aim is to give a value for the face approximation at the face $y_\edge$.
We begin with the case where $K$ does not correspond to a (local) maximum (\ie\ we suppose $(y-y_-)\,(y_+ -y) \geq 0$ and either $y \neq y_-$ or $y \neq y_+$) and, without loss of generality, we assume $y_- \geq y$ and $y \geq y_+$, with $y_- > y_+$.

\begin{center}
\begin{figure}[ht]
\begin{tikzpicture}[scale=1]
\draw[thin, ->] (0,0) -- (7,0) node(xline)[below] {$\bfx$};
\draw (2,-0.5) -- (2,1.5);
\draw (5,-0.5) -- (5,1.5);
\draw[color=bleuf, thick] (0,1.2) -- (2,1.2) node[midway, above]{$y_-$};
\draw[color=bleuf, thick] (2,0.9) -- (5,0.9) node[midway, above]{$y$};
\draw[color=bleuf, thick] (5,0.2) -- (7,0.2) node[midway, above]{$y_+$};
\node at (3.5,0.2) {$K$};
\node at (5,1.7) {\textcolor{vertf}{$\edge$}};
\node at (2,1.7) {\textcolor{vertf}{$\edge'$}};
\end{tikzpicture}
\caption{Notations for the one-dimensional problem} 
\label{fig:1D_notations}
\end{figure}
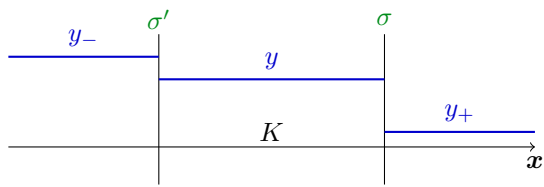
\end{center}

\medskip
The Lagrange-projection algorithm involves two steps: first, a reconstruction of the unknown in the cell $K$ supposing that the unknown is a step function taking the value $y_-$ on the left part of $K$ and the value $y_+$ on the right part, $y$ being the average of the unknown over $K$ (see Figure \ref{fig:rec_f}); second, the computation of $G_{K,\edge}$ as the integration over the time-step of the actual flux of the reconstructed function.

\begin{rmrk}[Lagrange-projection or method of characteristics?]
Introduced as such, the algorithm rather looks as a method of characteristics.
Another presentation, closer to the denomination ``Lagrange-projection", would be to consider the scheme as a three-steps algorithm: reconstruction, Lagrangian transport and projection by taking the mean value of each cell.
In fact, both computations lead to the same result (as, for instance, the Godunov scheme is equivalent to the upwind scheme), and the presentation chosen here offers the advantage to directly lead to a flux computation, which will be useful in the sequel.
\end{rmrk}

\begin{center}
\begin{figure}[h]
\begin{tikzpicture}[scale=1]
\draw[thin, ->] (0,0) -- (7,0) node(xline)[below] {$\bfx$};
\draw (2,-0.5) -- (2,1.5);
\draw (5,-0.5) -- (5,1.5);
\draw[color=bleuf, thick] (0,1.2) -- (2,1.2);
\draw[color=bleuf, thick] (2,0.9) -- (5,0.9);
\draw[color=bleuf, thick] (5,0.2) -- (7,0.2);
\draw[color=rougec, thick] (2,1.2) -- (4.1,1.2) -- (4.1,0.2) -- (5,0.2) node[midway, above]{$\hat y$};
\node at (3.5,0.2) {$K$};
\node at (5,1.7) {\textcolor{vertf}{$\edge$}};
\node at (2,1.7) {\textcolor{vertf}{$\edge'$}};
\end{tikzpicture}
\caption{Reconstruction step -- Blue: initial function. Red: reconstructed function in cell $K$} 
\label{fig:rec_f}
\end{figure}
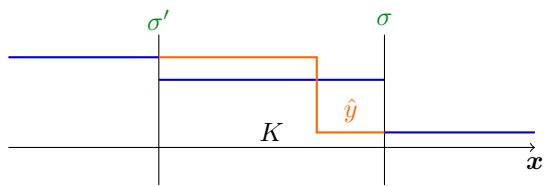
\end{center}

To alleviate the notations, we suppose that $K =(0,\delta x)$ and that the time interval under consideration is $(0,\delta t)$
The reconstructed function $\hat y$ over $K$, at the beginning of the time step so at $t=0$, reads:
\[
\hat y(x,0) = \left| \begin{array}{l} y_- \mbox{ if } x \leq \xi\, \delta x, \\[1ex]  y_+ \mbox{ otherwise,} \end{array} \right.
\]
where the constraint that the mean value of $\hat y$ over $K$ is $y$ yields:
\begin{equation}\label{eq:def_xi}
\xi = \frac{y_+-y}{y_+-y_-}.
\end{equation}
Let $\nu$ be a CFL number, defined by
\[
\nu = \frac{v \, \delta t}{\delta x}
\]
For $t \leq \delta x/v$, the transport of $\hat y$ at the velocity $v$ yields
\[
\hat y(\delta x,t) =
\left| \begin{array}{l} \displaystyle
y_+ \mbox{ if } t \leq \frac{(1-\xi)\, \delta x} v, 
\\[2ex]
y_- \mbox{ otherwise.}
\end{array} \right.
\]
Integrating over $(0,\delta t)$, we obtain that the 
numerical flux $G_{K,\edge}$ through $\edge$ outward $K$  satisfies:
\[
\delta t\ G_{K,\edge} = v\ 
\left| \begin{array}{l} \displaystyle
\delta t\ y_+ \mbox{ if } \delta t \leq \frac{(1-\xi)\, \delta x} v, 
\\[2ex] \displaystyle
\frac{(1-\xi)\ \delta x} v\ y_+  +  (\delta t - \frac{(1-\xi)\ \delta x} v)\ y_- \mbox{ if } \delta t \geq \frac{(1-\xi)\, \delta x} v.
\end{array} \right.
\]
Since, by definition of $y_\edge$, we have $ G_{K,\edge} = v\ y_\edge$, we get
\[
y_\edge = 
\left| \begin{array}{l} \displaystyle
y_+ \mbox{ if } \nu \leq 1-\xi, 
\\[2ex] \displaystyle
\bar y_\edge = \frac{1-\xi} \nu\ y_+  +  (1 - \frac{1-\xi} \nu)\ y_- \mbox{ if } \nu \geq 1-\xi.
\end{array} \right.
\]
Using \eqref{eq:def_xi}, we get:
\begin{equation} \label{eq:u_lim}
\bar y_\edge = y_- - \frac 1 \nu\ (y_- -y)
\end{equation}
Let us consider this latter expression as a function of $\nu$, for $\nu \in (0,1]$.
This function is increasing, from $-\infty$ ($\nu$ tending to zero) to $y$ for $\nu=1$.
In addition, thanks to \eqref{eq:def_xi}, $\nu=1-\xi$ implies $\bar y_\edge = y_+$.

\begin{rmrk}
These two results may be anticipated.
For the first one, from a Lagrangian point of view, the condition $\nu=1$ implies that the whole ``matter included in $K$" crosses $\edge$, so integrating the value of $\hat y$ at $\edge$ over $(0,\delta t)$ amounts to integrating the initial value of $\hat y$ over $(0,\delta x)$.
For the second one, the equality $\bar y_\edge = y_+$ for $\nu=1-\xi$ follows from the fact that, by construction, the expression of the flux is continuous with respect to $\delta t$.
\end{rmrk}

The behaviour of $\bar y_\edge(\nu)$ as a function of $\nu$ thus shows the equivalence between the two conditions $\nu \leq (1-\xi)$ and $\bar y_\edge(\nu) \leq y_+$, and the expression of $y_\edge$ may be recast as
\[
y_\edge = 
\left| \begin{array}{l} \displaystyle
y_+ \mbox{ if } \bar y_\edge \leq y_+, 
\\[2ex] \displaystyle
\bar y_\edge  \mbox{ otherwise,}
\end{array} \right.
\]
or just simply
\[
y_\edge = \max (y_+, \bar y_\edge).
\]
This latter expression may itself be obtained by the following two-step computation:
\begin{list}{-}{\itemsep=0.5ex \topsep=0.5ex \leftmargin=1.cm \labelwidth=0.3cm \labelsep=0.5cm \itemindent=0.cm}
\item[(i)] choose $y_+$ as tentative value for $y_\edge$,
\item[(ii)] define
\[
I_\edge = [y + \frac{1- \nu} \nu\ (y- y_-),\ y]
\]
as admissible interval and project the tentative value $y_+$ over $I_\edge$.
\end{list}
It is a straightforward exercise to check that the same conclusion may be drawn in the case where the unknown is increasing (\ie\ $y_- \leq y \leq y_+$ with $y_- < y_+$).
Finally, in the case where $y$ is a local maximum (\ie\ $(y-y_-)\,(y_+ -y) < 0$ or $y_- = y = y_+$), the definition of the reconstructed function $\hat y$ in not possible, and the propose two-step computation yields $y_\edge = y$ (\ie\ the upwind choice), which is a reasonable choice (and, in fact, the only one which ensures a discrete maximum principle, as we will see later).
We have thus recast the Lagrange-projection algorithm as a downwind scheme with a suitable limitation, which yields a discrete maximum principle, as we check in the following.

%
%
\medskip
\paragraph{\textbf{Checking the discrete maximum principle}}
The discrete maximum principle directly follows from the Lagrange-projection viewpoint on the scheme.
It may also be checked by showing, as usual, that the updated values of the unknown are convex combinations of the begining-of-step values.
Indeed, denoting by $\edge'$ the left face of the cell $K$, the above construction shows that
\[
y_{\edge'} - y = \alpha_{\edge'}\ (y_- -y) \mbox{ and } y_\edge - y = \alpha_\edge (y_- - \frac 1 \nu\ (y_- -y) -y) = \alpha_\edge \frac{1 - \nu} \nu (y - y_-)
\]
with $\alpha_\edge$ and $\alpha_{\edge'}$ in the interval $[0,1]$.
The first relation is due to the fact that $y_{\edge'}$ is obtained by projection of $y$ on an interval containing $y_-$, so is a convex combination of these two values.
The second relation is obtained by writing that $y_\edge$ is a convex combination of the two bounds of $I_\edge$.
We thus have, dropping for short the time index $n$:
\[
y^{n+1} = y - \nu (y_\edge - y_{\edge'}) = y - \nu \bigl[ (y_\edge - y) - (y_{\edge'} -y) \bigr]
= y -\nu \bigl[ \alpha_\edge \frac{1 - \nu} \nu (y - y_-) - \alpha_{\edge'}\ (y_- -y) \bigr].
\]
Hence,
\[
y^{n+1} = \bigl[ 1 - \alpha_\edge (1-\nu) -\alpha_{\edge'} \nu \bigr]\, y + \bigl[ \alpha_\edge (1-\nu) +\alpha_{\edge'} \nu \bigr]\, y_-,
\]
which yields that $y^{n+1}$ is a convex combination of $y$ and $y_-$ if $\nu \leq 1$.

%
%
\medskip
\paragraph{\textbf{Transport of step functions}}
It may be checked that the scheme exactly transports the Heaviside function (so, by an easy extension, any step function) in the following sense: at each time, the discrete solution in a cell corresponds to the mean value of the exact solution.
The proof of this result is obtained by induction, using the Lagrange-projection viewpoint on the scheme: suppose that the result is true at step $n$; then there is only one cell where the solution differs from 0 or 1, the reconstruction is exact in this cell and, finally, so are the fluxes.

\medskip
This feature of the scheme implies the following weaker property (P1): let us consider the case where $y_-=1$, $y \in (0,1)$ and $y_+=0$; then the expression \eqref{eq:u_lim} of $\bar y_\edge$ ensures that $y^{n+1}=1$ for $\nu$ close to 1.
The property (P1) will be used in the following to propose an extension of \eqref{eq:u_lim} to the case of a transport equation with a non-constant (in space and time) velocity.
%
%
\subsection{A one-dimensional transport problem}

\paragraph{\textbf{A tentative scheme}} - Let us now consider the problem
\begin{equation} \label{eq:pb1D_conv}
\partial_t(\rho y) + \partial_x(\rho y v) =0,
\end{equation}
where $\rho$ and $v$ satisfy
\begin{equation} \label{eq:pb1D_mass}
\partial_t \rho  + \partial_x(\rho v) =0.
\end{equation}
Under this latter assumption, the solution $y$ satisfies the transport equation
\[
\partial_t y + v\ \partial_x y =0,
\]
and $y$ thus obeys to a maximum principle.
For the solution of Equation \eqref{eq:pb1D_conv}, a derivation of the scheme based on the (exact) transport of a reconstructed function for the unknown $y$ is not straightforward, since the velocity is not constant and thus the transport would necessitate a reconstruction of the velocity itself.
We thus implement an alternate strategy, which consists in writing the scheme as a ``generic" downwind limited scheme, and then tuning the limitation to obtain a non-diffusive approximation of step functions, in the sense that the transition from one plateau to the another one is captured in only one cell.

\medskip
We use the same notations as in the previous section except that we now denote by $y_K$ the value of the unknown associated to the cell $K$ (instead of $y$ in the previous section) in order to prepare for the multidimensional setting (see below); we denote by $\rho_K$ the density in $K$, the discretization of the mass balance equation \eqref{eq:pb1D_mass} reads (skipping the exponent $~^n$ at time $t_n$ for the sake of simplicity):
\begin{equation} \label{eq:pb1D_mass_D}
\frac{|K|}{\delta t}\ (\rho^{n+1}_K - \rho_K) + F_{K,\edge} - F_{K,\edge'}=0.
\end{equation}
The discretization of Equation \eqref{eq:pb1D_conv} takes the following form:
\begin{equation} \label{eq:pb1D_conv_D}
\frac{|K|}{\delta t}\ (\rho_K^{n+1}\,y_K^{n+1} - \rho_K\,y_K) + F_{K,\edge}\,y_\edge - F_{K,\edge'}\,y_{\edge'}=0.
\end{equation}
Multiplying Equation \eqref{eq:pb1D_mass_D} by $y_K$ and subtracting to Equation \eqref{eq:pb1D_conv_D} yields the discrete form of the transport equation satisfied by $y$:
\[
\frac{|K|}{\delta t}\ \rho_K^{n+1}\,(y_K^{n+1} - y_K) + F_{K,\edge}\,(y_\edge-y_K) - F_{K,\edge'}\,(y_{\edge'}-y_K)=0,
\]
or, equivalently:
\begin{equation}\label{eq:y=}
y_K^{n+1}= y_K - \mathrm{sign}(F_{K,\edge})\ \nu (y_\edge-y_K) + \mathrm{sign}(F_{K,\edge'})\ \nu' (y_{\edge'}-y_K),
\end{equation}
with
\[
\nu = \frac{\delta t\ |F_{K,\edge}|}{|K|\ \rho_K^{n+1}},
\quad \nu'=\frac{\delta t\ |F_{K,\edge'}|}{|K|\ \rho_K^{n+1}}.
\]
To design the limitation process, we consider a specific case.
We suppose that $F_{K,\edge} > 0$, $F_{K,\edge'} >0$ and, in all the cells at the left (resp. right) side of $K$, the unknown is equal to $y_- >0$ (resp. is equal to $0$), with $0 < y_K < y_-$; in this case, we expect the scheme to allow $y_K^{n+1}$ to take the value $y_-$ (just) before limitation (this is the property (P1) introduced above).
In addition, we suppose that the scheme imposes $y_{\edge'}=y_-$ (and we will check a posteriori that it is indeed the case) and that it is a downwind-limited scheme, and thus that $y_\edge=0$ if the limitation is not active.
Equation \eqref{eq:y=} yields:
\[
y_K^{n+1}= y_K - \nu (y_\edge-y_K) + \nu' (y_- -y_K),
\]
so that we have, at the point where limitation becomes active:
\[
y_- = y_K - \nu (\bar y_\edge-y_K) + \nu' (y_- - y_K).
\]
\ie
\[
\bar y_\edge = y_K + \frac{1-\nu'}{\nu}(y_K-y_-).
\]
By the same arguments as in the previous section (monotonicity of $\bar y_\edge$ with respect to the time step, and $\bar y_\edge$ tends to $-\infty$ when the time step tends to zero), this result suggest that $y_\edge$ my be obtained by the projection of $y_+=0$ over the interval
\begin{equation}\label{eq:def_I_2}
I_\edge = [y_K + \frac{1- \nu'} \nu\ (y_K- y_-),\ y_K]
\end{equation}

\medskip
This computation suggests the following scheme, given here in the multi-dimensional setting for the sake of generality.
Let $\edge=K|L$ an internal edge of the mesh.
Without loss of generality, we suppose that $K$ is the upwind cell, \ie\ $F_{K,\edge} > 0$, so that $y_L = y_+$ in the preceding setting.
We denote by $op(\edge)$ the opposite site of $\edge$ with respect to $K$ (the definition of which is clear for structured meshes and may be extended to quadrangles or hexahedra via the mapping linking the actual and the reference element).
We define $\nu$ and $\nu'$ by
\[
\nu = \frac{\delta t\ |F_{K,\edge}|}{|K|\ \rho_K^{n+1}}, \qquad \nu' = \frac{\delta t\ |F_{K,op(\edge)}|}{|K|\ \rho_K^{n+1}}.
\]
Let $I_\edge$ be defined by Relation \eqref{eq:def_I_2}.
Then we define $y_\edge$ by the projection of $y_L$ onto $I_\edge$.
In addition, denoting by $M_{op(\edge)}$ the cell such that $op(\edge)=K|M_{op(\edge)}$, then for any $\edge \in \edges$, the face value $y_\edge$ satisfies the following two properties:
\begin{equation}\label{eq:props}
\begin{array}{ll}
(i) &
y_\edge - y_K = \alpha_{\edge,1}\ (y_L - y_K),
\\[1ex]
(ii) & \displaystyle
y_\edge - y_K = \alpha_{\edge,2}\ \xi\ (y_K-y_{M_{op(\edge)}}),\quad \mbox{with } \xi=\frac{1-\nu'}{\nu}=\frac{|K|\ \rho_K^{n+1} - \delta t\ |F_{K,op(\edge)}|}{\delta t\ |F_{K,\edge}|},
\end{array}
\end{equation}
with $\alpha_{\edge,1}, \alpha_{\edge,2} \in [0,1]$.

\medskip
By construction, this scheme satisfies the property (P1).

%
%
\bigskip
\paragraph{\textbf{Checking the discrete maximum principle}} - Let us recast \eqref{eq:y=} in the multi-dimensional setting, for $K$ a generic cell of the mesh:
\[
y_K^{n+1} = y_K - \frac{\delta t}{|K|\ \rho_K^{n+1}} \sum_{\edge \in \edges(K)} F_{K,\edge}\ (y_\edge - y_K)
\]
Using the properties $(i)$ and $(ii)$ of Equation \eqref{eq:props}, we get
\begin{multline}\label{eq:check_dmp}
y_K^{n+1} = y_K + \frac{\delta t}{|K|\ \rho_K^{n+1}}\ \Bigl[
\sum_{\edge \in \edges(K),\ F_{K,\edge} \geq 0} |F_{K,\edge}|\ \alpha_{\edge,2}
\ \frac{|K|\ \rho_K^{n+1} - \delta t\ |F_{K,op(\edge)}|}{\delta t\ |F_{K,\edge}|}\ (y_{M_{op(\edge)}} - y_K)
\\
+ \sum_{\edge \in \edges(K),\ F_{K,\edge} \leq 0} |F_{K,\edge}|\ \alpha_{\edge,1}\ (y_L - y_K)
\Bigr].
\end{multline}
The coefficient $c_K$ multiplying $y_K$ at the right-hand side of this relation reads:
\[
c_K = 1
- \sum_{\edge \in \edges(K),\ F_{K,\edge} \geq 0} \alpha_{\edge,2}\ \bigl[ 1 - \frac{\delta t\ |F_{K,op(\edge)}|}{|K|\ \rho_K^{n+1}}\bigr]
- \sum_{\edge \in \edges(K),\ F_{K,\edge} \leq 0} \alpha_{\edge,1}\ \frac{\delta t\ |F_{K,\edge}|}{|K|\ \rho_K^{n+1}}
\]
This relation shows that, without modification, the scheme cannot satisfy a discrete maximum principle (since we are not able to guarantee that $c_K \geq 0$).
For instance, for a one-dimensional problem, if the cell $K$ only has outward mass fluxes (\ie\ $F_{K,\edge} >0$ and $F_{K,\edge'} >0$, with $\edge$ and $\edge'$ the two faces of $K$), the coefficient $c_K$ may take values close to $-1$.
This problem may be traced back to the fact that the limit slope $(1-\nu')/\nu$ may blow up.
A possible modification of the scheme, which ensures that the scheme satisfies a discrete maximum principle at least with small CFL numbers, consists in giving an upper bound to this slope, that is to define the admissible interval $I_\edge$ as:
\begin{equation}\label{eq:def_I_3}
I_\edge = [y_K + \zeta\, (y_K - y_-),\ y_K], \qquad \zeta = \min(\frac{1- \nu'} \nu, S_{\rm max}),
\end{equation}
where $S_{\rm max}$ is a user-defined parameter.
With this choice, the following property, which is in some sense stronger than Property $(ii)$ of Equation \eqref{eq:props}, is guaranteed:
\begin{equation}\label{eq:props2}
(ii-b) \qquad y_\edge - y_K = \alpha_{\edge,2}\ S_{\rm max}\ (y_K-y_{M_{op(\edge)}}).
\end{equation}
Using now $(ii-b)$ instead of $(ii)$ in \eqref{eq:check_dmp}, we obtain for the coefficient $c_K$
\[
c_K = 1 - \frac{\delta t\ |F_{K,\edge}|}{|K|\ \rho_K^{n+1}}\ \Bigl[
\sum_{\edge \in \edges(K),\ F_{K,\edge} \geq 0} \alpha_{\edge,2}\ S_{\rm max}
+ \sum_{\edge \in \edges(K),\ F_{K,\edge} \leq 0} \alpha_{\edge,1} \Bigr],
\]
so $c_K \geq 0$ for $\delta t$ small enough (this relation showing that the stability time step decreases when the user-defined parameter $S_{\rm max}$ increases).

\medskip
\begin{rmrk}[On alternate direction versions of the Lagrange-projection step for incompressible flows]
For incompressible flows on structured meshes (typically, VOF applications in square or cubic domains), the one-dimensional Lagrange-projection scheme is often used combined with a direction-by direction strategy (sometimes, with a reconstruction step using a slightly smoother function than a step function).
When working direction-by-direction, the convection velocity, which is obtained by setting to zero $d-1$ components of the flow velocity and keeping the remaining one unchanged, is not divergence-free.
The only strategy to preserve a dicrete maximum principle for the convected function thus seems to be to switch to a discrete transport equation, which precisely what is done in this section.
Thus, the material presented here essentially applies to this case also.
\end{rmrk}
%
%
\bibliography{ereac}
\bibliographystyle{plain}
\end{document}